\documentclass[a4paper,12pt,reqno]{amsart}
\usepackage[utf8x]{inputenc}
\usepackage{amssymb,amsmath,latexsym}
\usepackage{epsfig}
\usepackage{mathrsfs}
\usepackage{hyperref}
\usepackage{color}

\usepackage{epic,eepic,wrapfig,color, ifthen}
\usepackage{url}
\usepackage{hyperref}

\newtheorem{thm}{Theorem}[section]
\newtheorem{cor}[thm]{Corollary}
\newtheorem{lemma}[thm]{Lemma}
\newtheorem{prop}[thm]{Proposition}

\newtheorem{remark}[thm]{Remark}
\newtheorem{example}[thm]{Example}
\numberwithin{equation}{section}

\newcommand{\formula}[2][nolabel]
{\ifthenelse{\equal{#1}{nolabel}}
 {\begin{align*} #2 \end{align*}}
 {\ifthenelse{\equal{#1}{}}
  {\begin{align} #2 \end{align}}
  {\begin{align} \label{#1} #2 \end{align}}
 }
}

\def\qed{{\hfill $\Box$ \bigskip}}
\def\R{{\mathbb R}}
\def\P{{\mathbb P}}
\def\E{{\mathbb E}}
\def\1{{\bf 1}}

\newcommand{\cal}[1]{\mathcal{#1}}
  
 \def\sE {{\cal E}} \def\sF {{\cal F}}

\def\sM {{\cal M}}  
  
 \def\sT {{\cal T}}

\def\R {{\mathbb R}}
\def\N {{\mathbb N}}

\def\nn{\nonumber}

\def\wt{\widetilde}

\def\E{{\mathbb E}}
\def\P{{\mathbb P}}
\def\eps{{\varepsilon}}

\def\bea{\begin{align*}}
\def\eea{\end{align*}}
\def\bee{\begin{equation}}
\def\eee{\end{equation}}

\makeatletter
\@addtoreset{equation}{section}

\makeatother

\begin{document}
\bibliographystyle{plain}

\title[Estimates on transition densities of subordinators]
{ \bf Estimates on transition densities of subordinators with  jumping density decaying in mixed polynomial orders}

\author
{Soobin Cho \quad and \quad Panki Kim}

\address[Cho]{Department of Mathematical Sciences,
	Seoul National University,
	Building 27, 1 Gwanak-ro, Gwanak-gu
	Seoul 08826, Republic of Korea}
\curraddr{}
\email{soobin15@snu.ac.kr}

\address[Kim]{Department of Mathematical Sciences and Research Institute of Mathematics,
	Seoul National University,
	Building 27, 1 Gwanak-ro, Gwanak-gu
	Seoul 08826, Republic of Korea}\thanks{This research is  supported by the National Research Foundation of Korea(NRF) grant funded by the Korea government(MSIP) (No. 2016R1E1A1A01941893).
}
\curraddr{}
\email{pkim@snu.ac.kr}

 \date{}

\maketitle

\begin{abstract}
In this paper, we discuss estimates on transition densities  for subordinators, which are global in time.
We establish the sharp two-sided estimates on the transition densities for subordinators whose L\'evy measures are absolutely continuous and decaying in mixed polynomial orders. Under a weaker assumption on  L\'evy measures,  we also obtain a precise asymptotic behaviors of the transition densities  at infinity. Our results cover geometric stable subordinators, Gamma subordinators and much more. 
\end{abstract}

\bigskip


\bigskip

\section{Introduction and Main results}
Since there are only a few known examples of stochastic processes for which the transition density can be computed explicitly, estimates and asymptotic behaviors of transition densities of stochastic processes  are extremely important and have studied a lot. 
When the process is symmetric and has a strong Markov property,  there are many beautiful results on this topic (see \cite{BKKL, BBCK, BGK, BGR14, Chen, CKK09, CKK, CK03, CK08, CKW, GHH17, GHH18, GHL14, HL, Mi} and references therein for estimates for symmetric jump processes).
But when the process is a non-symmetric jump process, estimates and asymptotic hehavior of its transition density are known much less. See \cite{CHXZ,
CZ16, CZ17, CZ18, GS1, GLT, GS, KS, KL, KR, KSV, Pi} and references therein. In this paper, we discuss estimates on transition densities  for  a large class of  non-decreasing L\'evy processes on $\R$. 

Let $S=(S_t)_{t\ge0}$ be a subordinator, that is, a non-decreasing L\'evy process on $\R$ with $S_0=0$. The process $S$ is characterized by its Laplace exponent $\phi$ which is given by
$$
\E e^{-\lambda S_t}=e^{-t \phi(\lambda)} \quad  \text{   for all } \; t, \lambda \ge 0.
$$
It is well known that $\phi$ is a Bernstein function with $\phi(0)=0$ and there exists a unique constant $a \ge 0$ and a Borel measure $\nu$ on $(0,\infty)$ satisfying $\int_0^{\infty} \text{min}\{1,s\} \nu(ds)<\infty$ such that
\begin{equation}\label{defphi}
\phi(\lambda) = a \lambda + \int_0^{\infty} (1-e^{-\lambda s}) \nu(ds).
\end{equation}
The constant $a$ is called the drift and $\nu$ is called the L\'evy measure of $S$ in the literature.

 The main objective of this paper is to obtain two-sided estimates on the transition density for a large class of subordinators. 
 Note that except a few special cases (see \cite{BKKK}), 
 the transition probability density of subordinators can not be computed explicitly along side an expression for the  L\'evy measure.
 Through subordination and inverse  subordination, the sharp estimates of the transition density of subordinators provide
 the sharp estimates of heat kernel of subordinate Markov process and 
 two-sided estimates for the fundamental solution, respectively. See \cite[Section 5]{GLT},  \cite[Section 4]{CKKK}. and \cite{BS, S}.
 
 Our assumptions are quite general; imposing only  (mixed) polynomially decaying conditions
  locally at zero (or infinity) on the density of the  L\'evy measure. This paper is a continuation of the authors' previous work in \cite{CK}. In \cite{CK}, we studied tail probabilities of subordinators under various decaying conditions on the tail of the L\'evy measure. In this paper, we concentrate on the case when the density of L\'evy measure is (mixed) polynomially decaying. (cf. conditions {\bf (S.Poly.)} and {\bf (L.Poly.)} in \cite{CK}.)
 
 Recently in \cite{GLT}, estimates for transition density  of subordinators have been also studied (see \cite[Theorem A and Theorem 4.15]{GLT}). In our context, their main assumptions on subordinators can be interpreted as our condition {\bf (S)} or {\bf (G)} holds with $0<\alpha_1 \le \alpha_2<1$. (See, \eqref{Slow} and \eqref{Shigh} below.) In this paper, by imposing scaling conditions on the L\'evy density directly, we allow the upper scaling index $\alpha_2$ at zero to be bigger than $1$. (cf. \cite[Remark 1.3(1)]{CK}.) Moreover, we establish the large time counterpart of that result. In this situation, we even allow that the lower scaling index at zero can be negative (see {\bf (L-3)} below). Hence, our results cover geometric stable subordinators. (See Example \ref{example} and Section 4.1. below.)

 In analysis on distributions of subordinators, by considering the subordinator $\wt{S}_t=S_t-at$, we may assume that $a=0$ without loss of generality. Hence, we always assume that $a=0$ in this paper. Moreover, {\it we always assume that the L\'evy measure $\nu$ has a density function $\nu(x)$ and the following Hartman-Wintner type condition holds throughout this paper.}
 
 \vspace{2mm}
 \noindent{\bf (E)} There exists a constant $T_0 \in [0, \infty)$ such that 
 $$
 \liminf_{x \to 0} x\nu(x) = 1/T_0,
 $$
 with a convention that $1/0=\infty$. 
 
 \smallskip
 
 In particular, the condition {\bf (E)} implies that $\nu(0,\infty) = \infty$ and hence the subordinator $S_t$ is not a compounded Poisson process. Moreover, as a consequence of this condition, we obtain the existence and boundedness of the transition density function.

\begin{prop}\label{p:existence}
For all $t>T_0$, the transition density $p(t,x)$ of the subordinator $S_t$ exists and is a continuous bounded function on $(0, \infty)$ as a function of $x$. 
\end{prop}
\proof
According to \cite[(64) and (74)]{HW}, (see also \cite[$({\rm HW}_{1/t})$]{KS13},) it suffices to show that
\begin{align*}
\liminf_{|\xi| \to \infty}\frac{\text{Re}\;\phi(i \xi)}{\log(1+ |\xi|)} \ge \frac{1}{T_0}.
\end{align*}

We first assume that $T_0>0$. Fix an arbitrary $\eps>0$. Then, by the assumption {\bf (E)}, there exists a constant $\delta>0$ such that $\nu(x) \ge (1-\eps)T_0^{-1}x^{-1}$ for $x \in (0, \delta)$. On the other hand, since a Gamma subordinator, whose Laplace exponent is $\log(1+\lambda)$, has the L\'evy density $s^{-1}e^{-s}$, we get the following equalities:
$$
\log\sqrt{1+\xi^2}= \text{Re}\;\log(1 + i \xi) = \text{Re}\; \int_0^\infty (1-e^{-i \xi s})s^{-1}e^{-s} ds = \int_0^\infty (1-\cos (\xi s)) s^{-1}e^{-s}ds.
$$
It follows that
\begin{align*}
\liminf_{|\xi| \to \infty}\frac{\text{Re}\;\phi(i \xi)}{\log(1+ |\xi|)} & = \liminf_{\xi \to \infty}\frac{\int_0^\infty (1-\cos (\xi s))\nu(s)ds }{\int_0^\infty (1-\cos (\xi s)) s^{-1}e^{-s}ds} \\
&\ge \frac{1-\eps}{T_0}\liminf_{\xi \to \infty}\frac{\int_0^\delta (1-\cos (\xi s))s^{-1}e^{-s}ds }{\int_0^\infty (1-\cos (\xi s)) s^{-1}e^{-s}ds} \\
& \ge \frac{1-\eps}{T_0}\left(1-\limsup_{\xi \to \infty}\frac{\int_\delta^\infty s^{-1}e^{-s}ds }{\log(1+ \xi)}\right) = \frac{1-\eps}{T_0}.
\end{align*}
Hence, we get the result by letting $\eps\to 0$. 

Then, we also deduce the result for the case when $T_0=0$ by letting $T_0 \to 0$.
\qed

Now, we enumerate our other main assumptions for the L\'evy measure $\nu$.

\vspace{1mm}

\noindent{\bf (S-1)} There are constants $c_1>0$, $R_1 \in (0, \infty]$ and $\alpha_1>0$ such that
\begin{align}\label{Slow}
\frac{\nu(r)}{\nu(R)} \ge c_1 \left(\frac{R}{r}\right)^{1+\alpha_1}    \quad \text{for all} \;\; 0<r \le R < R_1;
\end{align}

\noindent{\bf (S-2)} There are constants $c_2>0$, $R_1 \in (0, \infty]$ and $\alpha_2>0$ such that
\begin{align}\label{Shigh}
\frac{\nu(r)}{\nu(R)} \le c_2 \left(\frac{R}{r}\right)^{1+\alpha_2}    \quad \text{for all} \;\; 0<r \le R < R_1;
\end{align}

\noindent{\bf (S-3)} There are constants $c_3>0$ and $R_1 \in (0, \infty]$ such that
\begin{align}\label{S1}
\sup_{r \ge R_1} \nu(r) \le c_3,
\end{align}
with a convention that $\sup \emptyset = 0$;

\smallskip
\noindent{\bf (S-3*)} There are constants $c_4, c_5>0$ and $R_1 \in (0, \infty]$ such that
\begin{align}\label{S3*}
c_4\sup_{u \ge r} \nu(u) \le \nu(r) \quad \text{and} \quad c_5\nu(r) \le \nu(2r) \quad \text{for all} \;\; r \ge R_1/2.
\end{align}

\noindent{\bf (S)} There exist a common constant $R_1 \in (0, \infty]$ and constants $c_1,c_2,c_3>0$, $\alpha_2 \ge \alpha_1>0$ such that {\bf (S-1)}, {\bf (S-2)} and {\bf (S-3)} hold.

\vspace{3mm}

\noindent{\bf (L-1)} There are constants $c_6>0$, $R_2>0$ and $\alpha_3>0$ such that
\begin{align}\label{Llow}
\frac{\nu(r)}{\nu(R)} \ge c_6 \left(\frac{R}{r}\right)^{1+\alpha_3}    \quad \text{for all} \;\; R_2 \le r \le R < \infty;
\end{align}

\noindent{\bf (L-2)} There are constants $c_7>0$, $R_2>0$ and $\alpha_4>0$ such that
\begin{align}\label{Lhigh}
\frac{\nu(r)}{\nu(R)} \le c_7 \left(\frac{R}{r}\right)^{1+\alpha_4}    \quad \text{for all} \;\; R_2 \le r \le R < \infty;
\end{align}

\noindent{\bf (L-3)} There are constants $c_8,c_9>0$ and $R_3>0$ such that
\begin{align}\label{L-1}
\frac{\nu(r)}{\nu(R)} \ge c_8 \left(\frac{R}{r}\right)^{-c_9}  \;\;\quad \text{for all} \;\; 0<r \le R < R_3;
\end{align}

\noindent{\bf (L)} There exist a common constant $R_2 \in (0, \infty)$ and constants $c_6,c_7,c_8, c_9, R_3>0$, $\alpha_4 \ge \alpha_3>0$ such that {\bf (L-1)}, {\bf (L-2)} and {\bf (L-3)} hold.

\vspace{3mm}

\noindent{\bf (G)} The condition {\bf (S)}  holds with $R_1 = \infty$.

\vspace{1mm}

\begin{remark}\label{r:remark1}
	{\rm (1) The condition {\bf (S-1)} implies the condition {\bf (E)} with $T_0=0$ and the condition {\bf (L-3)} with $R_3 = R_1$.
		
		\noindent 
		(2) The constant $\alpha_1$ in the condition {\bf (S-1)} should be less than $1$. Indeed, since we have
		$$
		\infty>\int_0^{r} s\nu(s)ds \ge c_1\nu(r)r^{-1-\alpha_1} \int_0^{r}s^{-\alpha_1}ds \quad \text{ for all} \;\;  r \in (0, R_1),
		$$
		it must hold that $\alpha_1<1$.
		
		\noindent
		(3) A truncated $\alpha$-stable subordinator, whose L\'evy measure $\nu(ds)$ is given by 
		$$\nu(ds) = s^{-1-\alpha}\1_{(0,1)}(s)ds \quad (0<\alpha<1),$$
		satisfies the condition {\bf (S)} with $R_1=1$.

		\noindent
		(4) Clearly, the condition {\bf (S-3*)} implies the condition {\bf (S-3)} with the same constant $R_1$. Indeed, we get $\sup_{r \ge R_1}\nu(r) \le c_4^{-1} \nu(R_1)$ under the condition  {\bf (S-3*)}.

			\noindent 
		(5) 
		Let $0<\alpha_1 \le \alpha_2<1$ and  $m$ be a finite measure on $[\alpha_1, \alpha_2]$.
		Let $S$ be a subordinator without drift whose L\'evy measure $\nu(dx)$ is given by 
		$$\nu(dx) = \left(\int_{\alpha_1}^{\alpha_2} \frac{\beta}{\Gamma(1-\beta)x^{\beta+1}} m(d\beta)\right)dx.$$
		 Then, we can see that the subordinator $S$ satisfies the condition {\bf (G)}. Note that if $\alpha_1=\alpha_2=\alpha \in (0,1)$ and $m$ is a Dirac measure on $\alpha$, then $S$ is a $\alpha$-stable subordinator.
		
		\noindent 
	   (6) A geometric stable subordinator, whose Laplace exponent is $\log(1+\lambda^\alpha)$ for $\alpha \in (0,1)$, has the L\'evy density $\nu(x)$ such that $c_1^{-1}x^{-1} \le \nu(x) \le c_1x^{-1}$ for $x \in (0,1)$ while $c_1^{-1}x^{-1-\alpha}\le \nu(x) \le c_1x^{-1-\alpha}$ for $x \in [1,\infty)$, for some constant $c_1>1$. Hence it satisfies the condition {\bf (L)} while not satisfy {\bf (S)}. Note that the condition {\bf (E)} is satisfied with $T_0=1/\alpha>0$.
		
		\noindent 
		(7) The condition {\bf (L-3)} is very mild. For instance, if the L\'evy density is almost decreasing, then it holds trivially. Therefore, every subordinator whose Laplace exponent is a complete Bernstein function satisfies that assumption since its L\'evy measure has a completely monotone density. (See \cite[Chapter 16]{SSV} for examples of complete Bernstein functions.) 
}
\end{remark}

\smallskip

Following \cite{JP}, we let $H(\lambda)=\phi(\lambda)-\lambda \phi'(\lambda)$ and we define
\begin{align*}
b(t)&=(\phi'\circ H^{-1})(1/t) \qquad \text{and} \qquad w(r)=\nu(r,\infty).
\end{align*}

The function $H$ has an important role in estimates for the distributions of the subordinators. (see, e.g. \cite{JP, Mi}.)
Also, the function $b$ is used in authors' previous paper \cite{CK} (the definition of $b$-function in \cite{CK} is the same as $tb(t)$ in this paper) to describe a displacement with the highest probability of given subordinator at time $t$. We can see that $b(t)$ is strict increasing and $b(t) < \phi'(0)$.

From the definitions, we see that for every $\lambda>0$,
 \begin{equation}\label{basicH1}
H(\lambda) \ge \int_0^{1/\lambda} (1-e^{-\lambda s}-\lambda s e^{-\lambda s}) \nu(s) ds \ge \frac{1}{2e}\lambda^2 \int_0^{1/\lambda} s^2 \nu(s)ds
\end{equation}
and
 \begin{equation}\label{basicH2}
H(\lambda) \ge \int_{1/\lambda}^\infty (1-e^{-\lambda s}-\lambda s e^{-\lambda s}) \nu(s) ds \ge \frac{e-2}{e}w(1/\lambda).
\end{equation}
These inequalities follow from the facts that $1-e^{-x}-xe^{-x} =e^{-x}(e^x-1-x)\ge (2e)^{-1}x^2$ for all $0\le x\le 1$ and $1-e^{-x}-xe^{-x} \ge e^{-1}(e-2)$ for all $x \ge 1$. In particular, \eqref{basicH2} implies that for all $t>0$,
\begin{equation}\label{wHinverse}
w^{-1}(2e/t) \le w^{-1}(e(e-2)^{-1}/t) \le H^{-1}(1/t)^{-1}.
\end{equation}

\smallskip

Let 
$$
D(t):= t\max_{s \in [w^{-1}(2e/t), H^{-1}(1/t)^{-1}]} sH(s^{-1}).
$$
Then we define a function $\theta:(0, \infty) \times [0, \infty) \to (0, \infty)$ by
\begin{equation}\label{theta}
\theta(t,y):=
\begin{cases}
H^{-1}(1/t)^{-1}  &\text{ if} \;\; y \in [0, H^{-1}(1/t)^{-1}),\\
\min\left\{s \in \big[w^{-1}(2e/t), H^{-1}(1/t)^{-1}\big] : t s H(s^{-1}) = y\right\} &\text{ if} \;\; y \in [H^{-1}\left(1/t\right)^{-1}, D(t)], \\
w^{-1}(2e/t)  &\text{ if} \;\; y \in (D(t),\infty).
\end{cases}
\end{equation}
Note that $\theta(t,y) \in [w^{-1}(2e/t), H^{-1}(1/t)^{-1}]$ for all $t>0$ and $y \ge 0$. In particular, for each fixed $y \ge 0$, we have $\lim_{t \to 0} \theta(t,y)=0$ and $\lim_{t \to \infty} \theta(t,y) = \infty$. However, neither $t \mapsto \theta(t,y)$ nor $y \mapsto \theta(t,y)$ is a monotone function in general.

\smallskip

Following \cite{JP}, for $t>0$ and $x \in (0, t\phi'(0))$, we abbreviate
\begin{equation}\label{defsigma}
\sigma=\sigma(t,x):=(\phi')^{-1}(x/t).
\end{equation}
(This function is denoted by $\lambda_t$ in \cite{JP}.) Since $\phi'$ is non-increasing, $\sigma$ is a non-increasing function on $x$ for each fixed $t$ and a non-decreasing function on $t$ for each fixed $x$. From the definitions and the monotonicities of $H, \sigma$ and $b$, we have that for every $t>0$,
$$
tH(\sigma) > 1 \quad \text{for all} \;\; x \in (0,tb(t)) \qquad \text{and} \qquad tH(\sigma)|_{(t,x)=(t,tb(t))} =1.
$$

 Hereinafter, we denote $a\wedge b:=\min\{a,b\}$ and $a\vee b:=\max\{a,b\}$. The following theorems are the main results of this paper.
 
\begin{thm}\label{t:main1}
	Let $S$ be a subordinator satisfying {\bf (S)}. Then, for every $T>0$, there exist constants $c_1,c_2, c_3, c_5>1$ and $c_4>0$ such that the following estimates hold for all $t \in (0, T]$.
	
	\noindent
	\noindent {\rm (1)} It holds that for all $x \in (0, t b(t)]$,
	\begin{align}\label{e:main1}
	 \frac{c_1^{-1}}{\sqrt{t(-\phi''(\sigma))}}\exp\big(-tH(\sigma)\big) \le 	p(t,x) \le  \frac{c_1}{\sqrt{t(-\phi''(\sigma))}}\exp\big(-tH(\sigma)\big),
	\end{align}
	where $\sigma$ is defined as \eqref{defsigma}.  In particular, it holds that for all $x \in (0, tb(t)]$,
	\begin{align}\label{otherform}
	 c_2^{-1}H^{-1}(1/t)\exp\big(-2tH(\sigma)\big) \le 	p(t,x) \le  c_2H^{-1}(1/t)\exp\big(-\frac{t}{2}H(\sigma)\big)
	\end{align}

		\noindent
	{\rm (2)} It holds that for all $y \in [0, R_1/2)$,
	\begin{align}\label{e:main2}
	&c_3^{-1} H^{-1}(1/t) \min\left\{1, \frac{t \nu(y)}{H^{-1}(1/t)} + \exp \big(-\frac{c_4y}{\theta(t, y/(8e^2))} \big) \right\} \nn\\
	&\qquad \le p(t, tb(t)+y) \le c_3 H^{-1}(1/t) \min\left\{1, \frac{t \nu(y)}{H^{-1}(1/t)} + \exp \big(-\frac{y}{8\theta(t, y/(8e^2))} \big) \right\},
	\end{align}
	where $\theta(t,y)$ is defined as \eqref{theta}. In particular, for all $y \in (D(t), R_1/2)$,
	\begin{align}\label{pure}
 c_5^{-1} t \nu(y)	\le p(t, tb(t)+y) \le c_5 t \nu(y).
	\end{align}
	
	 Moreover, if  $S$ also satisfies  the condition {\bf (S-3*)}, then \eqref{e:main2} holds for all $y \in [0, \infty)$ and \eqref{pure} holds for all $y \in (D(t), \infty)$.
	
	\end{thm}

\begin{thm}\label{t:main2}
	Let $S$ be a subordinator satisfying {\bf (E)} and {\bf (L)}. 
	
	\noindent
\noindent {\rm (1)} There exist constants $T_1>T_0$, $c_1,c_2,c_3,c_5>1$ and $c_4>0$ such that for all $t \in [T_1, \infty)$, \eqref{e:main1} holds for all $x \in (0,tb(t)]$, \eqref{otherform} holds for all $x \in [tb(T_1), tb(t)]$, \eqref{e:main2} holds for all $y \in [0, \infty)$ and \eqref{pure} holds for all $y \in (D(t), \infty)$. 	

\noindent
{\rm (2)} If $T_0=0$ in the condition {\bf (E)}, then for every $T>0$, there are comparison constants such that for all $t \in [T, \infty)$, \eqref{e:main1} holds for all $x \in (0,tb(t)]$, \eqref{otherform} holds for all $x \in [tb(T),tb(t)]$, \eqref{e:main2} holds for all $y \in [0, \infty)$ and \eqref{pure} holds for all $y \in (D(t), \infty)$.
\end{thm}

\begin{cor}\label{c:main1}
	Let $S$ be a subordinator satisfying {\bf (G)}. Then, there exist constants $c_1,c_2,c_3,c_5>1$ and $c_4>0$ such that for all $t \in (0, \infty)$, \eqref{e:main1} and \eqref{otherform} hold for all $x \in (0,tb(t)]$, \eqref{e:main2} holds for all $y \in [0, \infty)$ and \eqref{pure} holds for all $y \in (D(t), \infty)$.
\end{cor}
\vspace{3mm}

Our main theorems also cover the cases when $\alpha_3 \le 1$ and $\alpha_4 \ge 2$. In such cases, the exponential term in the right tail estimates may have an efficient effect on estimates at specific times while have no role in other time values. (See, Section 4.2.)
Note that since the condition {\bf (E)} guarantees the existence of a continuous bounded transition density function $p(t,x)$ only for $t>T_0$, we should choose the constant $T_1$ bigger than $T_0$ in Theorem \ref{t:main2}.

\smallskip

If we impose additional conditions on decaying orders of the density of L\'evy measure, then we can simplify the right tail estimates in our theorems. Consider the following further conditions:

\vspace{2mm}

\noindent{\bf (S.Pure)} Condition {\bf (S)} holds with $\alpha_2<2$.

\vspace{1mm}

\noindent{\bf (L.Pure)} Condition {\bf (L)} holds with $\alpha_4<2$.

\vspace{1mm}

\noindent{\bf (L.Mixed)} Condition {\bf (L)} holds with $\alpha_3>1$.

\vspace{2mm}

\begin{remark}\label{r:remark2}
	{\rm (1) Since $\alpha_1$ should be less than $1$, (see Remark \ref{r:remark1}(2),) there is no analogous condition to {\bf (L.Mixed)} concerning the condition {\bf (S)}. 
	
	\noindent 
(2) We have $\phi'(0)<\infty$ under the condition {\bf (L.Mixed)}. Indeed, we see that
$$
\phi'(0)=\int_0^{R_3}s\nu(s)ds + \int_{R_3}^\infty s\nu(s)ds \le c + c_6^{-1} \nu(R_3)R_3^{1+\alpha_3} \int_{R_3}^\infty s^{-\alpha_3}ds <\infty.
$$

}
\end{remark}

\smallskip

Under either of the conditions {\bf (S.Pure)} or {\bf (L.Pure)}, we obtain pure jump type estimates on the right tails of $p(t,x)$. 

Recall that $\sigma=(\phi')^{-1}(x/t)$ for $t>0$ and $x \in (0, t\phi'(0))$. In the following corollary, we let $\sigma=0$ for $t>0$ and $x \ge t\phi'(0)$ so that $x \mapsto \sigma$ is a non-increasing function on $(0, \infty)$ for each fixed $t>0$. We use the notation $z_+ = \max\{z, 0\}$ for $z \in \R$.

\begin{cor}\label{cor1}
	Let $S$ be a subordinator satisfying {\bf (S.Pure)}. Then, for every $T>0$, there exists a constant $c_1>1$ such that for all $t \in (0, T]$ and $y \in [0, R_1/2)$,
	\begin{equation}\label{e:pj}
	c_1^{-1} \left( H^{-1}(1/t)\wedge t\nu(y) \right) \le p(t, tb(t)+y) \le c_1 \left( H^{-1}(1/t)\wedge t\nu(y) \right).
	\end{equation}
	Therefore, there exists a constant $c_2>1$ such that for all $t \in (0, T]$ and $x \in (0, R_1/2)$,
	\begin{align}\label{e:pj'}
	&c_2^{-1}\min\bigg\{H^{-1}(1/t) \exp\big(-2tH(\sigma)\big), t \nu\big((x-tb(t))_+\big)  \bigg\} \nn\\
	&\qquad \le p(t,x) \le c_2 \min\bigg\{H^{-1}(1/t) \exp\big(-\frac{t}{2}H(\sigma)\big), t \nu\big((x-tb(t))_+\big)  \bigg\}.
	\end{align}

	 Moreover, if $S$ also satisfies the condition {\bf (S-3*)}, then  \eqref{e:pj} holds for all $t \in (0, T]$ and  $y \in [0, \infty)$, and \eqref{e:pj'} holds for all $t \in (0,T]$ and $x \in (0, \infty)$.	
\end{cor}

\begin{cor}\label{cor2}
	Let $S$ be a subordinator satisfying {\bf (E)} and {\bf (L.Pure)}. Then, there exist constants $T_1>T_0$ and $c_1>1$ such that \eqref{e:pj} holds for all $t \in [T_1,\infty)$ and $y \in [0, \infty)$, and \eqref{e:pj'} holds for all $t \in [T_1, \infty)$ and $x \in [tb(T_1), \infty)$.
	
Moreover, if $T_0=0$ in the condition {\bf (E)}, then for every $T>0$, there are  comparison constants such that \eqref{e:pj} holds for all $t \in [T, \infty)$ and $y \in [0, \infty)$, and \eqref{e:pj'} holds for all $t \in [T, \infty)$ and $x \in [tb(T), \infty)$.
\end{cor}

\vspace{2mm}

Under the condition {\bf (L.Mixed)}, we can find a monotone function which is easy to compute and can play the same role as the function $\theta$.
Define
$$
\mathscr{H}(r):= \inf_{s \ge r} \frac{1}{sH(s^{-1})} \qquad \text{and} \qquad \mathscr{H}^{-1}(u):= \sup\{r \in \R : \mathscr{H}(r) \le u\}.
$$
Recall that under the condition {\bf (L.Mixed)}, $\phi'(0)$ is finite.
See \eqref{push} and a line below.

\begin{cor}\label{cor3}
	Let $S$ be a subordinator satisfying {\bf (E)} and {\bf (L.Mixed)}. Then, there exist constants $T_1>T_0$, $c_1>1$ and $c_2,c_3>0$ such that for all $t \in [T_1, \infty)$ and $y \in [0, \infty)$,
	\begin{align*}
&c_1^{-1} H^{-1}(1/t) \min\left\{1, \frac{t \nu(y)}{H^{-1}(1/t)} + \exp \big(-\frac{c_2y}{\mathscr{H}^{-1}(t/y)} \big) \right\} \\
&\qquad \le p(t, t\phi'(0)+y) \le c_1 H^{-1}(1/t) \min\left\{1, \frac{t \nu(y)}{H^{-1}(1/t)} + \exp \big(-\frac{c_3y}{\mathscr{H}^{-1}(t/y)} \big) \right\}.
\end{align*}

Moreover, if $T_0=0$ in the condition {\bf (E)}, then for every $T>0$, there are comparison constants such that the above estimates hold for all $t \in [T, \infty)$ and $y \in [0,\infty)$.
\end{cor}
The above corollary may be considered as a counterpart of 
\cite[Theorem 1.5(2)]{BKKL} where a similar result was obtained for symmetric jump processes. (See, Section \ref{s:R}.)

In this paper, we also discuss the precise asymptotical properties of densities of subordinators. (cf, \cite{DR, GLT}.)
The asymptotic expressions are given in terms of $\phi$ and its derivatives.
Under the condition {\bf (L-3)} we show in Corollary \ref{exactconv} that the density of the subordinator is asymptotically  equal to $(2\pi t(-\phi''(\sigma)))^{-1/2} \exp\big(-tH(\sigma)\big)$ as $t \to \infty$. If, in addition,  the constant $T_0=0$ in the condition {\bf (E)} is zero then the same result holds as $x \to 0$. In Example \ref{example}, we apply Corollary \ref{exactconv} to geometric stable subordinators and get the exact asymptotic behavior of the transition density  of  geometric stable subordinators as $t \to \infty$.
Up to authors' knowledge,  since Pillai introduced the series formula \eqref{exactden}  of the transition density  of  geometric stable subordinator in \cite{Pi} in 1990, its exact asymptotic behaviors given in 
\eqref{exactasymp0} and
\eqref{generalform}--\eqref{exactasymp2}
 have been unknown.

\vspace{3mm}

{\bf Notations}: In this paper, the positive constants $T_0, \alpha_1,\alpha_2,\alpha_3, \alpha_4, R_1,R_2$ and $R_3$ will remain the same. Lower case letters $c$'s without subscripts denote strictly positive constants whose values are unimportant and which may change even within a line, while values of lower case letters with subscripts $c_i$, $i=0,1,2,...$ are fixed in each statement and proof, and the labeling of these constants starts anew in each proof.

We use the symbol ``$:=$'' to denote a definition, 
which is read as ``is defined to be.''  Recall that  $a\wedge b:=\min\{a,b\}$ and $a\vee b:=\max\{a,b\}$.

The notation $f(x) \asymp g(x)$ means that there exist comparison constants $c_1,c_2>0$ such that $c_1g(x)\leq f (x)\leq c_2 g(x)$ for the specified range of the variable $x$. On the other hand, the notation $f(x) \simeq g_1(x) + g_2(x)h(cx)$ means that there exist comparison constants $c_3,c_4,c_5,c_6>0$ such that $c_3(g_1(x)+g_2(x)h(c_4x))\le f(x) \le c_5(g_1(x)+g_2(x)h(c_6x))$ for the specified range.

\vspace{3mm}

\section{Auxiliary functions and basic estimates}

Recall that
\begin{align*}
&H(\lambda)= \phi(\lambda)-\lambda\phi'(\lambda) = \int_0^\infty (1- e^{-\lambda s} - \lambda s e^{-\lambda s})\nu(s)ds,\\
&b(t)=(\phi' \circ H^{-1})(1/t) = \int_0^{\infty} se^{-H^{-1}(1/t)s} \nu(s) ds,  \qquad w(r)= \nu(r, \infty).
\end{align*}

Since $\phi$ is a Bernstein function with $\phi(0)=0$, we see that $H(0)=0$ and $H$ is strictly increasing on $(0, \infty)$. Also, it is easy to see that
$H(R) \le (R/r)^2 H(r)$ for all $R \ge r >0$.
Recall that $H$ satisfies 
\eqref{basicH1} and \eqref{basicH2}.
Moreover, it holds that
\begin{align}\label{e:CK21}
e^{-1}\lambda^2 \int_0^{1/\lambda} s w(s)ds \le H(\lambda) \le 5\lambda^2 \int_0^{1/\lambda} s w(s)ds \quad \text{for all } \lambda>0.
\end{align}
Indeed, by  the Fubini's theorem, $\phi(\lambda)/\lambda = \int_0^\infty \int_0^{s}e^{-\lambda u}  \nu(s)du ds= \int_0^\infty e^{-\lambda u}\int_u^\infty \nu(s)ds du = \int_0^\infty e^{-\lambda u}w(u)du$ for all $\lambda>0$. It follows that
$$
\frac{H(\lambda)}{\lambda^2} = -\left(\frac{\phi(\lambda)}{\lambda}\right)' = \int_0^\infty e^{-\lambda s} sw(s)ds \quad \text{for all} \; \lambda>0.
$$
Since $\int_0^{1/\lambda} e^{-\lambda s} sw(s)ds \ge e^{-1}\int_0^{1/\lambda} sw(s)ds \ge e^{-1}w(1/\lambda) \int_0^{1/\lambda} sds =2^{-1}e^{-1}\lambda^{-2}w(1/\lambda) $ and 
$\int_{1/\lambda}^\infty e^{-\lambda s} sw(s)ds \le w(1/\lambda) \int_{1/\lambda}^\infty se^{-\lambda s} ds = 2e^{-1} \lambda^{-2}w(1/\lambda)$, we get 
$$ \int_0^{1/\lambda} e^{-\lambda s} sw(s)ds \le \int_0^\infty e^{-\lambda s} sw(s)ds \le 5 \int_0^{1/\lambda} e^{-\lambda s} sw(s)ds \quad \text{for all} \; \lambda>0.
$$
Therefore, since $e^{-1} \le e^{-\lambda s} \le 1$ for $s \in [0, 1/\lambda]$, we get \eqref{e:CK21}.

\smallskip

Recall that the condition {\bf (L-3)} is weaker than the condition {\bf (S-1)}. (See, Remark \ref{r:remark1}(1).) Hence, the following lemma also hold under the condition {\bf (S-1)}.

We denote $\phi^{(n)}$ the $n$-th derivative of the function $\phi$.

\begin{lemma}\label{l:L1diff}
	Suppose that {\bf (L-3)} holds. 
	
	\noindent {\rm (1)} For every $\lambda_0>0$, there are constants $c_n>1$, $n=1,2,...$ such that
	\begin{align*}
	e^{-1}\int_0^{1/\lambda} s^n \nu(s)ds \le |\phi^{(n)}(\lambda)| \le c_n\int_0^{1/\lambda} s^n \nu(s)ds \qquad \text{for all} \;\; \lambda \ge \lambda_0 \;\; \text{and} \;\; n \ge 1.
	\end{align*}
	\noindent {\rm (2)} For every $\lambda_0>0$, there are constants $c_n'>1$, $n=1,2,...$ such that
	\begin{align*}
	c_n'^{-1}|\phi^{(n)}(2\lambda)| \le |\phi^{(n)}(\lambda)| \le c_n'|\phi^{(n)}(2\lambda)| \qquad \text{for all} \;\; \lambda \ge \lambda_0 \;\; \text{and} \;\; n \ge 1.
	\end{align*}
	
	\noindent {\rm (3)} For every $\lambda_0>0$, there are constants $c_n''>0$, $n=1,2,...$ such that
	$$
	|\lambda \phi^{(n+1)}(\lambda)| \le c_n''|\phi^{(n)}(\lambda)| \qquad \text{for all} \;\; \lambda \ge \lambda_0 \;\; \text{and} \;\; n \ge 1.
	$$
\end{lemma}
\proof (1) First, we see that for all $\lambda>0$ and $n \ge 1$,
\begin{equation}\label{difflow}
|\phi^{(n)}(\lambda)|  \ge \int_0^{1/\lambda} s^n e^{-\lambda s} \nu(s) ds \ge e^{-1} \int_0^{1/\lambda} s^n  \nu(s) ds.
\end{equation}

On the other hand, we have that for all $\lambda \ge 2R_3^{-1}$ and $n \ge 1$,
\begin{align}\label{diffbound}
|\phi^{(n)}(\lambda)| &= \int_0^{1/\lambda} s^n e^{-\lambda s} \nu(s) ds + \int_{1/\lambda}^{R_3} s^{n+c_9} e^{-\lambda s} s^{-c_9}\nu(s) ds + \lambda^{-n}\int_{R_3}^\infty \lambda^ns^{n} e^{-\lambda s} \nu(s) ds\nn \\
& \le \int_0^{1/\lambda} s^n e^{-\lambda s} \nu(s) ds + c \lambda^{c_9}\nu(1/\lambda) \int_{1/\lambda}^{R_3} s^{n+c_9} e^{-\lambda s} ds+c\lambda^{-n} \int_{R_3}^\infty e^{-\lambda s/2} \nu(s) ds \nn \\
& \le\int_0^{1/\lambda} s^n \nu(s) ds+ c\lambda^{-n-1} \nu(1/\lambda)\int_1^\infty u^{n+c_9}e^{-u}du + c \lambda^{-n}e^{-\lambda R_2/2} w(R_3),
\end{align}
where $c_9>0$ is the constant in \eqref{L-1}.
We used the assumption \eqref{L-1} and the fact that for every $n \ge 1$, there exists a constant $c>0$ such that $x^n \le ce^{x/2}$ for all $x \ge 0$ in the first inequality and the change of the variables $u=\lambda s$ in the second inequality.

Using the assumption \eqref{L-1}  (twice) and  the inequality $x^n \le ce^{x/2}$ again,  it also hold that for all $\lambda \ge 2R_3^{-1}$,
\begin{align}\label{scaleint}
\int_0^{1/\lambda} s^n \nu(s) ds &\ge \nu(1/\lambda)\int_{1/(2\lambda)}^{1/\lambda} s^n\frac{\nu(s)}{\nu(1/\lambda)}ds \ge
c\nu(1/\lambda)\int_{1/(2\lambda)}^{1/\lambda} s^n (s\lambda)^{c_9}ds \nn\\ &\ge c \lambda^{-n-1} \nu(1/\lambda) 
\ge c \lambda^{-n-1-c_9} R_3^{-c_9}\nu(R_3) \ge  c\lambda^{-n}e^{-\lambda R_3/2}w(R_3),
\end{align}
We deduce from \eqref{diffbound} and \eqref{scaleint} that $|\phi^{(n)}(\lambda)| \le c \int_0^{1/\lambda} s^n\nu(s)ds$ for all $\lambda \ge 2R_3^{-1}$ and $n \ge 1$. Then, by considering the constants $\inf_{ \lambda \in [\lambda_0, 2R_3^{-1}]} \big(|\phi^{(n)}(\lambda)|^{-1} \int_0^{1/\lambda}s^n\nu(s)ds\big)$, we get the desired result.

\smallskip
\noindent
(2) By (1), the change of the variables and the assumption \eqref{L-1}, we have that for $\lambda \ge 2R_3^{-1}$ and $n \ge 1$, 
$$
|\phi^{(n)}(2\lambda)| \asymp \int_0^{1/(2\lambda)}s^n \nu(s)ds = 2^{-n-1}\int_0^{1/\lambda}s^n \nu(s/2)ds \ge c \int_0^{1/\lambda}s^n \nu(s)ds \asymp |\phi^{(n)}(\lambda)|.
$$
By (1), we also get 
$$
|\phi^{(n)}(\lambda)| \ge e^{-1}  \int_0^{1/(2\lambda)}s^n \nu(s)ds
\asymp |\phi^{(n)}(2\lambda)|.
$$
By considering $\inf_{ \lambda \in [\lambda_0, 2R_3^{-1}]} |\phi^{(n)}(\lambda)/\phi^{(n)}(2\lambda)|$ and $\sup_{ \lambda \in [\lambda_0, 2R_3^{-1}]} |\phi^{(n)}(\lambda)/\phi^{(n)}(2\lambda)|$, we get the result.

\smallskip
\noindent
(3) By (1), we have that for all $\lambda \ge \lambda_0$ and $n \ge 1$,
\begin{align*}
|\lambda \phi^{(n+1)}(\lambda)| \asymp \int_0^{1/\lambda} (\lambda s)s^n \nu(s)ds \le \int_0^{1/\lambda} s^n \nu(s)ds \asymp |\phi^{(n)}(\lambda)|,
\end{align*}
which yields the result. \qed

\begin{lemma}\label{l:L2diff}
	Suppose that {\bf (L-1)} holds. 
	
	\noindent {\rm (1)} For every $\lambda_0>0$, there are constants $c_n>1$, $n=1,2,...$ such that
	\begin{align*}
	e^{-1}\int_0^{1/\lambda} s^n \nu(s)ds \le |\phi^{(n)}(\lambda)| \le c_n\int_0^{1/\lambda} s^n \nu(s)ds \qquad \text{for all} \;\; 0<\lambda \le \lambda_0 \;\; \text{and} \;\; n \ge 1.
	\end{align*}
	
	\noindent {\rm (2)} For every $\lambda_0>0$, there are constants $c_n'>1$, $n=1,2,...$ such that
	\begin{align*}
	c_n'^{-1}|\phi^{(n)}(2\lambda)| \le |\phi^{(n)}(\lambda)| \le c_n'|\phi^{(n)}(2\lambda)| \qquad \text{for all} \;\; 0<\lambda \le \lambda_0 \;\; \text{and} \;\; n \ge 1.
	\end{align*}
	
	\noindent {\rm (3)} For every $\lambda_0>0$, there are constants $c_n''>0$, $n=1,2,...$ such that
	$$
	|\lambda \phi^{(n+1)}(\lambda)| \le c_n''|\phi^{(n)}(\lambda)| \qquad \text{for all} \;\; 0<\lambda \le \lambda_0 \;\; \text{and} \;\; n \ge 1.
	$$
\end{lemma}
\proof (1) By \eqref{difflow}, it remains to prove the upper bounds. By  \eqref{Llow} and the first line in \eqref{scaleint}, we have that for all $0<\lambda \le R_2^{-1}$ and $n \ge 1$,
\begin{align*}
|\phi^{(n)}(\lambda)| &\le \int_0^{1/\lambda} s^n \nu(s)ds + \int_{1/\lambda}^\infty s^n e^{-\lambda s}\nu(s)ds \le \int_0^{1/\lambda} s^n \nu(s)ds + c\nu(1/\lambda)\int_{1/\lambda}^\infty s^n e^{-\lambda s}ds \\
&\le \int_0^{1/\lambda} s^n \nu(s)ds + c\lambda^{-n-1}\nu(1/\lambda) \le c \int_0^{1/\lambda} s^n \nu(s)ds.
\end{align*}
Again, by considering the constants $\inf_{ \lambda \in [R_2^{-1}, \lambda_0]} \big(|\phi^{(n)}(\lambda)|^{-1} \int_0^{1/\lambda}s^n\nu(s)ds\big)$, we get the result.

\noindent
(2) As in the proof of (1), it suffices to prove for $0<\lambda \le R_2^{-1}$ and $n \ge 1$. By (1), the change of the variables and  \eqref{Llow},
\begin{align*}
|\phi^{(n)}(2\lambda)| &\asymp \int_0^{2R_2}s^n \nu(s)ds + \int_{R_2}^{1/(2\lambda)}s^n \nu(s)ds \\
&=\int_0^{2R_2}s^n \nu(s)ds  +  2^{-n-1}\int_{2R_2}^{1/\lambda}s^n \nu(s/2)ds\\
& \asymp \int_0^{2R_2}s^n \nu(s)ds  +  \int_{2R_2}^{1/\lambda}s^n \nu(s)ds \asymp |\phi^{(n)}(\lambda)|.
\end{align*}

\noindent
(3) We get the result by the same proof as the one for Lemma \ref{l:L1diff}(3). \qed

\begin{lemma}\label{l:smallH}
Suppose that {\bf (S-1)} holds.

\noindent {\rm (1)} There are constants $c_1,c_2>0$ such that
$$
w(r) \ge c_1 r \nu(2r) \quad \text{and} \quad w(r) \le c_2 r \nu(r) \qquad \text{for all} \;\; 0<r < R_1/2.
$$
In particular, if we further assume that {\bf (S-2)} holds, then $w(r) \asymp r \nu(r)$ for all $r \in (0, R_1/2)$.

\noindent {\rm (2)} There is a constant $c_3>0$ such that
$$ \frac{w(r)}{w(R)} \ge c_3 \left(\frac{R}{r}\right)^{\alpha_1} \qquad \text{for all} \;\; 0 < r \le R < R_1/2.
$$

\noindent {\rm (3)} For every $r_0>0$, there is a constant $c_4>0$ such that
$$
\frac{H(R)}{H(r)} \ge c_4 \left( \frac{R}{r} \right)^{\alpha_1}   \qquad \text{for all} \;\; r_0 \le r \le R<\infty.
$$
In particular, 
$$
\frac{H^{-1}(t)}{H^{-1}(s)} \le c_4^{-1/\alpha_1} \left( \frac{t}{s} \right)^{1/\alpha_1}   \qquad \text{for all} \;\; H(r_0) \le s \le t<\infty.
$$
\noindent {\rm (4)} 
For every $\lambda_0>0$, there are comparison constants such that
$$
H(\lambda)  \asymp \lambda^2 \int_0^{1/\lambda} s^2\nu(s)ds  \asymp \lambda^2(-\phi''(\lambda))  \qquad \text{for all} \;\; \lambda_0 \le \lambda<\infty.
$$
\end{lemma}
\proof
(1) 
By \eqref{Slow}, we have that 
for all $r \in (0, R_1/2)$,
$$
w(r) \ge \int_{r}^{2r} \nu(s) ds= \nu(2r)\int_r^{2r} \frac{\nu(s)}{\nu(2r)} ds \ge c_1 \nu(2r)\int_r^{2r} (2r/s)^{1+\alpha_1} ds  \ge c_1r \nu(2r) 
$$
and
\begin{align*}
w(r) &= \int_r^{R_1} \nu(s)ds + w(R_1) \le \frac{\int_{R_1/2}^{R_1} \nu(s)ds + w(R_1)}{\int_{R_1/2}^{R_1} \nu(s)ds} \int_r^{R_1} \nu(s)ds \\
& \le c\nu(r) \int_r^{R_1} \frac{\nu(s)}{\nu(r)}ds \le  cc_1^{-1} r^{1+\alpha_1}\nu(r) \int_r^{R_1} s^{-1-\alpha_1} ds  \le cc_1^{-1}\alpha_1^{-1}r\nu(r).
\end{align*}

Moreover, if {\bf (S-2)} holds, then there is a constant $c>0$ such that $\nu(2r) \ge c\nu(r)$ for all $0<r < R_1/2$. Hence, we obtain $w(r) \asymp r\nu(r)$ for all $r \in (0, R_1/2)$.

\noindent
(2) By (1) and \eqref{Slow}, for all $0<2r \le R<R_1/2$,
$$
\frac{w(r)}{w(R)} \ge c\frac{r}{R}\frac{ \nu(2r)}{ \nu(R)} \ge c \frac{r}{R} \left( \frac{R}{2r} \right)^{1+\alpha_1} = c2^{-1-\alpha_1} \left( \frac{R}{r} \right)^{\alpha_1}.
$$
On the other hand, for all $0<r \le R \le (R_1/2) \wedge (2r)$, by the monotonicity of $w$,
$$
\frac{w(r)}{w(R)} \ge 1 \ge 2^{-\alpha_1}\left( \frac{R}{r} \right)^{\alpha_1}.
$$
Hence,  we get the result in both cases.

\noindent
(3) By \eqref{e:CK21}, the change of the variables and (2), we have that for all $2R_1^{-1} < r \le R$,
\begin{align*}
\frac{H(R)}{H(r)} \ge \frac{1}{5e} \left(\frac{R}{r} \right)^2 \frac{\int_0^{1/R} sw(s)ds}{\int_0^{1/r} sw(s)ds} =  \frac{1}{5e} \frac{\int_0^{1/r} sw(rs/R)ds}{\int_0^{1/r} sw(s)ds}  \ge \frac{c_3}{5e} \left(\frac{R}{r} \right)^{\alpha_1} .
\end{align*}
Moreover, for all $r_0 \le r \le 2R_1^{-1} < R$, we see that
$$
\frac{H(R)}{H(r)} \ge  \frac{H(R)}{H(2R_1^{-1})} \ge \frac{c_3}{5e}  \left(\frac{R}{2R_1^{-1}} \right)^{\alpha_1} \ge \frac{c_3}{5e}  \left(\frac{r_0}{2R_1^{-1}} \right)^{\alpha_1} \left(\frac{R}{r} \right)^{\alpha_1}.
$$
Lastly, for all $r_0 \le r \le R \le 2R_1^{-1}$, we get $H(R)/H(r)  \ge 1 \ge (r_0R_1/2)^{\alpha_1} (R/r)^{\alpha_1}$. Hence, the assertion holds.

\noindent

\noindent
 (4) As in the proof of Lemma \ref{l:L1diff}, we may and do assume that $\lambda_0 > 2R_1^{-1}$. By (1), we get 
 $$
 \lambda^2\int_0^{1/\lambda}s^2 \nu(s) ds \ge c_2^{-1}\lambda^2\int_0^{1/\lambda}s w(s) ds  \ge c_2^{-1}\lambda^2 w(1/\lambda)\int_0^{1/\lambda}s  ds = (2c_2)^{-1}w(1/\lambda).
 $$
 Hence, from the definition of $H$, we get
 $$
 H(\lambda)  \le \lambda^2\int_0^{1/\lambda}s^2 \nu(s) ds + w(1/\lambda) \le (2c_2+1)\lambda^2\int_0^{1/\lambda}s^2 \nu(s) ds.
 $$
 Therefore, by combining with \eqref{basicH1}, we obtain the first comparison. Then, according to \cite[(5.6) and (5.7) in Lemma 5.1]{JP}, 
 \eqref{basicH2} and the first comparison imply  $H(\lambda) \asymp \lambda^2(-\phi''(\lambda))$. This completes the proof. 
\qed

Similar results to Lemma \ref{l:smallH} hold under the condition {\bf (L-1)}.

\begin{lemma}\label{l:largeH}
	Suppose that {\bf (L-1)} holds.
	
	\noindent {\rm (1)} There are constants $c_1, c_2>0$ such that
	$$
	w(r) \ge c_1 r \nu(2r) \quad \text{and} \quad w(r) \le c_2 r \nu(r) \qquad \text{for all} \;\; r \ge R_2.
	$$
	In particular, if we further assume that {\bf (L-2)} holds, then $w(r) \asymp r \nu(r)$ for all $r \in [R_2, \infty)$.

	\noindent {\rm (2)} There is a constant $c_3>0$ such that
	$$
\frac{w(r)}{w(R)} \ge c_3 \left(\frac{R}{r}\right)^{\alpha_3} \qquad \text{for all} \;\; R_2 \le r \le R < \infty.
	$$
	
	\noindent {\rm (3)} For every $r_0>0$, there is a constant $c_4>0$ such that
	$$
	 \frac{H(R)}{H(r)} \ge c_4 \left( \frac{R}{r} \right)^{\alpha_3 \wedge (3/2)}   \qquad \text{for all} \;\; 0 < r \le R \le r_0.
	$$
	In particular, 
	$$
	\frac{H^{-1}(t)}{H^{-1}(s)} \le c_4^{-1/(\alpha_3 \wedge (3/2))} \left( \frac{t}{s} \right)^{(1/\alpha_3) \vee (2/3)}   \qquad \text{for all} \;\;  0<  s \le t \le H(r_0).
	$$
	
	\noindent {\rm (4)} For every $\lambda_0>0$, there are comparison constants such that
	$$
	H(\lambda)  \asymp \lambda^2 \int_0^{1/\lambda} s^2\nu(s)ds \asymp \lambda^2(-\phi''(\lambda))   \qquad \text{for all} \;\; 0 < \lambda \le \lambda_0.
	$$
\end{lemma}
\proof
(1) For all $r \ge R_2$, we see from \eqref{Llow} that $w(r) \ge \int_r^{2r} \nu(s)ds \ge cr\nu(2r)$ and $w(r) \le c r^{1+\alpha_3} \nu(r) \int_r^{\infty}s^{-1-\alpha_3} ds \le cr\nu(r)$. Moreover, if {\bf (L-2)} holds further, then $\nu(r) \asymp \nu(2r)$ for all $r \ge R_2$ and hence $w(r) \asymp r\nu(r)$ for all $r \ge R_2$.

\noindent
(2) This follows from \eqref{Llow} and (1). (See, the proof of Lemma \ref{l:smallH}(2).)

\noindent
(3) By the proof of Lemma \ref{l:smallH}(3), using the monotoncity of  $H$, it suffices to show that 
\begin{equation}\label{Hlarge}
\frac{H(R)}{H(r)} \ge c_3 \left( \frac{R}{r} \right)^{\alpha_3 \wedge (3/2)} \qquad \text{for all} \;\; 0 < r \le R \le (2R_2)^{-1}. 
\end{equation}
According to \eqref{e:CK21} and (1), for all $0<\lambda \le (2R_2)^{-1}$, 
\begin{align*}
 \lambda^{-2} H(\lambda) \asymp\int_0^{1/\lambda} s w(s)ds \le  c\int_0^{R_2} s w(s)ds + c\int_{R_2}^{1/\lambda} s^2 \nu(s)ds \le c\int_{R_2}^{1/\lambda} s^2 \nu(s)ds.
 \end{align*}
The  holds since $\int_{R_2}^{1/\lambda} s^2 \nu(s)ds \ge \int_{R_2}^{2R_2} s^2 \nu(s)ds \ge c \ge c\int_{0}^{R_2} s w(s)ds$ where we have used \eqref{Llow}. Therefore, in view of \eqref{basicH1}, we get
\begin{equation}\label{e:dw1}
H(\lambda) \asymp \lambda^2 \int_{R_2}^{1/\lambda} s^2 \nu(s) ds \qquad \text{for all} \;\; 0<\lambda \le (2R_2)^{-1}.
\end{equation}

Let $\alpha_3' = \alpha_3 \wedge (3/2) \in (0,2)$. Observe that by \eqref{Llow}, for all $0<r \le R \le (2R_2)^{-1}$, 
\begin{align*}
&\int_{1/R}^{1/r} s^2 \nu(s)ds = \int_{1/R}^{1/r}  s^{1-\alpha_3'}s^{1+\alpha_3'}\nu(s)ds \le c R^{-1-\alpha_3'} \nu(1/R) \int_{1/R}^{1/r}  s^{1-\alpha_3'}ds \\
&\;\; \le  c r^{\alpha_3'-2}R^{-1-\alpha_3'} \nu(1/R) = c r^{\alpha_3'-2}R^{2-\alpha_3'} R^{-3} \nu(1/R)  \le  c (R/r)^{2-\alpha_3'} \int_{1/(2R)}^{1/R} s^2\nu(s)ds.
\end{align*}
Thus
$$
\int_{R_2}^{1/r} s^2 \nu(s)ds \le \int_{R_2}^{1/R} s^2 \nu(s)ds+ c (R/r)^{2-\alpha_3'} \int_{1/(2R)}^{1/R} s^2\nu(s)ds \le
c (R/r)^{2-\alpha_3'}   \int_{R_2}^{1/R} s^2 \nu(s)ds.
 $$
It follows that by \eqref{e:dw1}, for all $0<r \le R \le (2R_2)^{-1}$,
\begin{align*}
\frac{H(R)}{H(r)} &\ge c \left(\frac{R}{r}\right)^{2} \frac{ \int_{R_2}^{1/R} s^2 \nu(s)ds}{\int_{R_2}^{1/r} s^2 \nu(s)ds} \ge c \left(\frac{R}{r}\right)^{2} \frac{ \int_{R_2}^{1/R} s^2 \nu(s)ds}{(R/r)^{2-\alpha_3'}\int_{R_2}^{1/R} s^2 \nu(s)ds} \ge c \left(\frac{R}{r}\right)^{\alpha_3'}.
\end{align*}
This proves \eqref{Hlarge}.

\noindent
(4) The first comparison follows from \eqref{basicH1} and \eqref{e:dw1}. Then, 
using  \cite[(5.6) and (5.7) in Lemma 5.1]{JP}, 
we obtain the second comparison from
the first comparison
and  \eqref{basicH2}.
 \qed

\vspace{2mm}

\vspace{2mm}

Now, we give some basic properties of the $b$-function. It is easy to verify that $b$ is strictly increasing, $\lim_{t \to 0}b(t)=0$ and $\lim_{t \to \infty} b(t) = \phi'(0) \in (0, \infty]$. Also, by \cite[Lemma 2.4]{CK},
\begin{equation}\label{estib}
\frac{1}{\phi^{-1}(1/t)} \le tb(t) \le \frac{1}{\phi^{-1}(c_*/t)}, \quad c_*=\frac{e-2}{e^2-e}, \;\;\qquad \text{for all} \;\; t>0.
\end{equation}

The following lemma is useful when $\phi$ is not comparable to the function $H$.

\begin{lemma}\label{l:b}
For every $a_2 \ge a_1 >0$ and $a_3>0$, it holds that for all $t>0$,
\begin{align*}
&tb(t/a_1)-tb(t/a_2) \le 2e\frac{a_2}{H^{-1}(a_2/t)} + e^{-1}\frac{tw(H^{-1}(a_2/t)^{-1})}{H^{-1}(a_1/t)}\le \frac{2e^2-4e+1}{e-2} \frac{a_2}{ H^{-1}(a_1/t)},\\
&tb(t/a_3)-tb(t/(4a_3))  \ge \frac{1}{2}\frac{tH^{-1}(4a_3/t)^2 (\phi''\circ H^{-1})(4a_3/t)}{H^{-1}(4a_3/t)}.
\end{align*}

In particular, if the condition {\bf (S-1)} holds, (resp. {\bf (L-1)} holds,) then for every $a_2 \ge a_1>0$, $a_3>0$ and $T>0$, there exist $c_1,c_2>0$ such that for all $t \in (0, T]$, (resp. $t \in [T, \infty)$,)
\begin{align*}
tb(t/a_1)-tb(t/a_2) \le c_1H^{-1}(1/t)^{-1} \quad\;\; \text{and} \quad\;\; tb(t/a_3)-tb(t/(4a_3)) \ge c_2 H^{-1}(1/t)^{-1}.
\end{align*}
\end{lemma}
\proof
By the mean value theorem and \eqref{basicH2}, we have
\begin{align*}
b(t/a_1)&-b(t/a_2) =\int_0^\infty\left(e^{-sH^{-1}(a_1/t)}-e^{-sH^{-1}(a_2/t)}\right) s \nu(s) ds \\
& \le  H^{-1}(a_2/t) \int_0^{H^{-1}(a_2/t)^{-1}}s^2 \nu(s) ds + \int_{H^{-1}(a_2/t)^{-1}}^\infty se^{-sH^{-1}(a_1/t)}  \nu(s) ds \\
& \le \frac{2e H(H^{-1}(a_2/t))}{ H^{-1}(a_2/t)} + \frac{e^{-1}}{H^{-1}(a_1/t)} \int_{H^{-1}(a_2/t)^{-1}}^\infty \nu(s) ds\\
&=\frac{2ea_2}{t H^{-1}(a_2/t)} + \frac{e^{-1}}{H^{-1}(a_1/t)} w(H^{-1}(a_2/t)^{-1})\\
& \le \frac{2ea_2}{t H^{-1}(a_1/t)} + \frac{H(H^{-1}(a_2/t))}{(e-2) H^{-1}(a_1/t)} = \frac{2e^2-4e+1}{e-2} \frac{a_2}{t H^{-1}(a_1/t)}.
\end{align*}
In the second inequality, we used the fact that for every $\lambda>0$, the map $s \mapsto se^{-\lambda s}$ on $[0,\infty)$ has the maximum value $\lambda^{-1}e^{-1}$.

On the other hand, we also have that by the mean value theorem,
\begin{align*}
&b(t/a_3)-b(t/(4a_3)) =\int_0^\infty\left(e^{-sH^{-1}(a_3/t)}-e^{-sH^{-1}(4a_3/t)}\right) s \nu(s) ds \\
& \ge (H^{-1}(4a_3/t)-H^{-1}(a_3/t)) \int_0^{\infty} e^{-sH^{-1}(4a_3/t)}s^2 \nu(s) ds  \ge \frac{1}{2}H^{-1}(4a_3/t) (\phi''\circ H^{-1})(4a_3/t).
\end{align*}
In the last inequality, we used the fact that $H^{-1}(\kappa^2 \lambda) \ge \kappa H^{-1}(\lambda)$ for all $\kappa \ge 1$ and $\lambda \ge 0$ which follows from that $H(\kappa \lambda) \le \kappa^2 H(\lambda)$ for all $\kappa \ge 1$ and $\lambda \ge 0$.

Then, by Lemmas \ref{l:smallH} and \ref{l:largeH}, we can see that the last assertion holds.
\qed

\vspace{2mm}

\section{Main results}

In this section, we give proofs for our main results. Recall that we always assume that $S=(S_t)$ is a subordinator without drift, whose L\'evy measure has a density function satisfying the condition {\bf (E)} with the constant $T_0 \in [0, \infty)$.

\subsection{Estimates on left tail probabilities} 
In this subsection, we study estimates on $p(t,x)$ when $x$ is small.
We first present a result established in \cite{GLT}, which holds under the condition {\bf (S-1)}. Recall from \eqref{defsigma} that we use the abbreviation 
$\sigma=\sigma(t,x)=(\phi')^{-1}(x/t)$ for   $t>0$ and  $0<x < t \phi'(0).$

\begin{prop}\label{lefttail}
	Suppose that the condition {\bf (S-1)} holds. Then, for every $T>0$, there exists a constant $M_0>0$ such that for all $t \in (0, T]$ and $x \in (0, t b(t/M_0)]$,
	\begin{align}\label{e:lefttail}
	p(t,x) \asymp \frac{1}{\sqrt{t(-\phi''(\sigma))}}\exp(-tH(\sigma)).
	\end{align}
\end{prop}
\proof
According to Lemma \ref{l:smallH}(3) and (4), we can see that for every $x_0>0$, the condition $-\phi'' \in \text{WLSC}(\alpha_1-2,c,x_0)$ in \cite[Theorem 3.3]{GLT} is satisfied with some constant $c>0$.  Since $x \mapsto \sigma$ decreases for each fixed $t$, we have that for $t \in (0, T]$ and $x \in (0, tb(t/M_0)]$,
$$
\sigma \ge ((\phi')^{-1} \circ b)(t/M_0) = H^{-1}(M_0/t) \ge H^{-1}(M_0/T).
$$
Also, by the above inequality and Lemma \ref{l:smallH}(4), there exists a constant $c_1>0$ such that
$$
t\sigma^2 (-\phi''(\sigma)) \ge c_1tH(\sigma) \ge c_1tH(H^{-1}(M_0/t)) = c_1M_0.
$$
Hence, the result follows from \cite[Theorem 3.3]{GLT}. \qed

Now, we establish left tail probabilities under the conditions {\bf (L-1)} and {\bf (L-3)}.  Since subordinators can not decrease, if $x$ is small compare to $t$, then left tail probabilities mainly depend on small jumps of subordinators. This is why we impose the assumption {\bf (L-3)} on small jumps in the condition {\bf (L)}.

Define a function $\sM:(0,\infty) \times (0, \infty) \times (-\infty, \infty) \to \mathbb{C}$ by
\begin{equation}\label{defM}
\sM(s,z,u):=\phi(z+\frac{iu}{\sqrt{s(-\phi''(z))}}) - \phi(z) -  \phi'(z)\frac{iu}{\sqrt{s(-\phi''(z))}}.
\end{equation}
In the settings of \cite{GLT}, the Laplace exponent $\phi$ should satisfy a lower weak scaling condition at infinity (i.e., the lower Matuszewska index (at infinity) of the function $\phi(\lambda)\1_{\{\lambda \ge 1\}}$ should be strictly bigger than $0$.) It follows that a map $u \mapsto e^{-t\sM(t,\sigma,u)}$  for fixed $t>0$ decreases at least subexponentially. This property has an important role in the proof of \cite[Theorem 3.3]{GLT}. 
Unlike \cite{GLT}, in our settings, the Laplace exponent $\phi$ can be slowly varying at infinity so that the map $u \mapsto e^{-t\sM(t,\sigma,u)}$ can decays in polynomial order. Therefore, we should  bound the integral $\int_{-\infty}^\infty e^{-t\sM(t,\sigma,u)}du$ more carefully in the following proposition.

\begin{prop}\label{lefttail:L}
	Suppose that the conditions {\bf (L-1)} and {\bf (L-3)} holds. Then, there exist $T_1>T_0$ , $M_0>0$ and comparison constants such that \eqref{e:lefttail} holds for all $t \in [T_1, \infty)$ and $x \in (0, t b(t/M_0)]$.
	
		Moreover, if $T_0=0$ in the condition {\bf (E)}, then for every $T>0$, there exist $M_0>0$ and comparison constants such that \eqref{e:lefttail} holds for all $t \in [T, \infty)$ and $ x\in (0, tb(t/M_0)]$.
\end{prop}
\proof
Recall that   $\sM$ is defined in \eqref{defM}. Since $\phi'(\sigma)=x/t$,
by the Fourier-Mellin inversion formula (see, e.g. \cite[(4.3)]{Sc},) and the change of the variables, we have
\begin{align}\label{Mellin}
p(t,x) &=\frac{e^{-t\phi(\sigma)+\sigma x}}{2\pi} \int_{-\infty}^\infty \exp \left(-t\big(\phi(\sigma + iu) - \phi(\sigma)  \big) + iux \right) du \nn\\
 &=\frac{e^{-t(\phi(\sigma)-\sigma \phi'(\sigma))}}{2\pi} \int_{-\infty}^\infty \exp \left(-t\big(\phi(\sigma + iu) - \phi(\sigma)   - iu\phi'(\sigma)\big) \right) du \nn\\
&=\frac{e^{-tH(\sigma)}}{2\pi\sqrt{t(-\phi''(\sigma))}} \int_{-\infty}^\infty e^{-t \sM(t,\sigma,u)} du,
\end{align}
whenever the integral converges. Note that if $|e^{-t \sM(t, \sigma,u)}|$ is integrable on $\R$ with respect to $u$, then
\begin{align*}
\int_{-\infty}^\infty e^{-t \sM(t,\sigma,u)} du = \int_{0}^\infty (e^{-t \sM(t,\sigma,u)}+e^{-t \sM(t,\sigma,-u)}) du.
\end{align*}
Since the complex conjugate of $\sM(t,\sigma,u)$ is given by
\begin{align*}
\overline{\sM(t,\sigma,u)} &= \int_0^{\infty}  \bigg(1-\exp\big(-(\sigma+\frac{-iu}{\sqrt{t(-\phi''(\sigma))}})s\big)\bigg)\nu(s)ds -\phi(\sigma) - \phi'(\sigma)\frac{-iu}{\sqrt{t(-\phi''(\sigma))}}\\
&= \sM(t,\sigma,-u),
\end{align*}
we have that $e^{-t \sM(t,\sigma,u)}+e^{-t \sM(t,\sigma,-u)} \in \R$ for all $t,\sigma>0$ and $u \in \R$. Hence, $p(t,x)$ is a real number whenever $|e^{-t \sM(t, \sigma,u)}|$ is integrable on $\R$ with respect to $u$.

\smallskip

 Let $T>T_0$ be a constant which will be chosen later and fix any $\delta>0$ such that $T \ge T_0+\delta$.
 We claim that the integral in \eqref{Mellin} converges for all $t \ge T$ and $\sigma>0$. Indeed, by a similar proof to that of  Proposition \ref{p:existence}, for $\eps=\delta/(2T_0+2\delta)$, there are constants $\sigma_0>0, \xi_0>1$ such that 
\begin{align}\label{e:sig1}
\nu(s) \ge \frac{(1-\eps/2)}{T_0}s^{-1} \quad \text{for all} \;\; s \in (0, -\sigma_0^{-1}\log(1-\eps/2))
\end{align}
and
\begin{align}\label{e:sig2}
\int_0^{-\log(1-\eps/2)} (1 - \cos( \xi s))s^{-1} ds \ge \frac{1-\eps}{1-\eps+\eps^2/4} \log(1+ \xi) \quad \text{for all} \;\;  \xi \ge \xi_0.
\end{align}
It follows that for all $|u|>\xi_0 (\sigma_0 \vee \sigma)\sqrt{t(-\phi''(\sigma))}$, 
\begin{align}\label{logdecay}
\text{Re}& \; t\sM(t,\sigma,u)  = t\int_0^\infty e^{-\sigma s}\big(1-\cos \frac{us}{\sqrt{t(-\phi''(\sigma))}}\big) \nu(s)ds\nn\\
& \ge \frac{(1-\eps/2)^2t}{T_0}\int_0^{\frac{-\log (1-\eps/2) }{\sigma_0 \vee \sigma}} \big(1-\cos \frac{us}{\sqrt{t(-\phi''(\sigma))}}\big) s^{-1}ds  \nn\\
& =  \frac{(1-\eps/2)^2t}{T_0}\int_0^{-\log (1-\eps/2)} \big(1-\cos \frac{us}{(\sigma _0\vee \sigma)\sqrt{t(-\phi''(\sigma))}}\big) s^{-1}ds \nn\\
& \ge \frac{(1-\eps)t}{T_0} \log \big( 1+ \frac{u }{(\sigma_0 \vee \sigma)\sqrt{t(-\phi''(\sigma))}} \big).
\end{align}
Since $(1-\eps)t/T_0 \ge (1-\eps)T/T_0  \ge (1-\eps)(T_0+\delta)/T_0=1+\delta/(2T_0)>1$, we see from \eqref{logdecay} that $|e^{-t \sM(t, \sigma,u)}| = e^{-t \text{Re}\; \sM(t,\sigma,u)}$ is integrable on $\R$ with respect to $u$. This yields that \eqref{Mellin} holds.

\smallskip
Next, we will show that there exists a constant $M_0>1$ such that for all $t \in [T, \infty)$ and $x \in (0, tb(t/M_0)]$,
\begin{equation}\label{claimM}
\sqrt{\pi} \le \int_{-\infty}^\infty e^{-t \sM(t,\sigma,u)} du \le 2\sqrt{\pi},
\end{equation}
which implies \eqref{e:lefttail} in view of \eqref{Mellin}. 

Define
$$
\sT_0=\sT_0(t,\sigma):=(\sigma_0 \vee \sigma)\sqrt{t(-\phi''(\sigma))} \qquad \text{and} \qquad \sT=\sT(t,\sigma):=\sigma\sqrt{t(-\phi''(\sigma))}.
$$
Clearly, we have $\sT_0 \ge \sT$. For $\sigma > \sigma_0$, we see from \eqref{e:sig1} that $\sT^2 \ge \sigma^2 t \int_0^{1/\sigma} s^2 e^{-\sigma s} \nu(s) ds \ge c_1\sigma^2 t  \int_0^{1/\sigma } s^2 s^{-1}ds = c_1t/2 \ge c_1T/2$. On the other hand, for $\sigma \le \sigma_0$, we see from Lemma \ref{l:largeH}(4) and the monotonicity of the function $\sigma$ that $\sT^2 \ge c_2 tH(\sigma) \ge c_2 t(H \circ (\phi')^{-1} \circ b)((t/M_0)) = c_2M_0$. It follows that
\begin{equation}\label{sTsize}
\sT_0^2 \ge \sT^2 \ge (c_1T/2) \wedge (c_2M_0).
\end{equation}
We claim that
\begin{equation}\label{CLT}
\lim_{\sT \to \infty} \int_{-\infty}^\infty e^{-t\sM(t, \sigma, u)}du = \int_{-\infty}^\infty e^{-\frac{1}{2}u^2}du = \sqrt{2\pi},
\end{equation}
which yields the desired result. Indeed, if \eqref{CLT} is true, then there exists a constant $c_3>0$ such that \eqref{claimM} holds for $\sT \ge c_3$. By choosing $T=2c_1^{-1}c_3^2$ and $M_0=c_2^{-1}c_3^2$, we get the result from \eqref{sTsize}. Now, we prove \eqref{CLT}.

First, we note that according to \eqref{logdecay},
\begin{align}\label{e:eq1}
\left|\int_{|u| > \xi_0 \sT_0} e^{-t \sM(t,\sigma,u)} du\right| \le 2\int_{\xi_0 \sT_0}^\infty \big( 1+ \frac{u}{\sT_0} \big)^{-(2T_0+\delta)/(2T_0)} du.
\end{align}

On the other hand, by Taylor's theorem, we have
\begin{align*}
\bigg|t\sM(t,\sigma,u )-\frac{1}{2}u^2\bigg| &= \bigg|t\big(\phi(\sigma + \frac{i\sigma }{\sT}u) - \phi(\sigma) - \phi'(\sigma)\frac{i\sigma }{\sT}u\big)-\frac{1}{2}u^2\bigg| \\
& \le \frac{1}{2}u^2\sup_{z \in [-|u|,|u|]}  \bigg|(-\phi''(\sigma + \frac{i\sigma}{\sT}z))\frac{t\sigma^2}{\sT^2} - 1\bigg|\\
& = \frac{1}{2}\frac{u^2}{(-\phi''(\sigma))}\sup_{z \in [-|u|,|u|]}  \bigg|-\phi''(\sigma + \frac{i\sigma}{\sT}z) + \phi''(\sigma)\bigg|.
\end{align*}
Note that
\begin{align*}
&\sup_{z \in [-|u|,|u|]}  \bigg|-\phi''(\sigma + \frac{i\sigma}{\sT}z) + \phi''(\sigma)\bigg| \le \sup_{z \in [-|u|,|u|]} \int_0^\infty s^2e^{-\sigma s}\bigg|\cos(\frac{ \sigma zs}{\sT})-1-i \sin(\frac{ \sigma zs}{\sT})\bigg|  \nu(s) ds \\
&\quad = 2\sup_{z \in [-|u|,|u|]} \int_0^\infty s^2e^{-\sigma s} \big|\sin (\frac{\sigma z s}{2 \sT}) \big| \nu(s) ds \le \frac{\sigma |u|}{\sT} \int_0^\infty s^3e^{-\sigma s}   \nu(s) ds = \frac{\sigma |u|}{\sT} \phi'''(\sigma).
\end{align*}
We used the fact that $|\sin x| \le |x|$ for all $x\in \R$ in the second inequality. Hence, we get
$$
\bigg|t\sM(t,\sigma,u )-\frac{1}{2}u^2\bigg| \le \frac{\sigma \phi'''(\sigma)}{2\sT(-\phi''(\sigma))}|u|^3.
$$
Then, combining with the fact that $|e^z-1| \le |z|e^{|z|}$ for $z \in \mathbb{C}$, it follows that for all $u \in \R$,
\begin{align}\label{taylor}
\left|e^{-t \sM(t,\sigma,u)}-e^{-\frac{1}{2}u^2} \right| &= e^{-\frac{1}{2}u^2} \left|\exp\big(\frac{1}{2}u^2-t \sM(t,\sigma,u)\big)-1 \right| \nn\\
&\le \frac{\sigma \phi'''(\sigma)}{2\sT(-\phi''(\sigma))}|u|^3 \exp\big(-\frac{1}{2}{u^2} +  \frac{\sigma \phi'''(\sigma)}{2\sT(-\phi''(\sigma))}|u|^3\big).
\end{align}
Below, we consider the cases $\sigma >\sigma_0$ and $\sigma \le \sigma_0$, separately.

\smallskip
\noindent
(Case 1): Assume that $\sigma>\sigma_0$. By Lemma \ref{l:L1diff}(3), there exists a constant $c_4>0$ such that $\sigma \phi'''(\sigma) \le c_4(-\phi''(\sigma))$. Let $\xi_1 = (2c_4)^{-1} \wedge \xi_0$. Then, according to \eqref{taylor}, 
\begin{align}\label{e:eq2}
&\left|\int_{|u| \le \xi_1  \sT} (e^{-t \sM(t,\sigma,u)}-e^{-\frac{1}{2}u^2}) du \right|\le \frac{c_4}{\sT} \int_0^{ \xi_1  \sT} u^3\exp\big(-(\frac{1}{2} -  \frac{c_4}{2\sT}u) u^2\big) du \nn\\
& \le \frac{c_4}{\sT} \int_0^{ \xi_1  \sT} u^3\exp\big(-(\frac{1}{2} -  \frac{c_4 \xi_1}{2}) u^2\big) du \le \frac{c_4}{\sT} \int_0^{\infty} u^3\exp\big(-\frac{1}{4}u^2\big) du \le \frac{c_5}{\sT}.
\end{align}

On the other hand, note that $\sigma |u|/\sT > \sigma_0 \xi_1$ for $|u|>\xi_1 \sT$. Hence, by Lemma \ref{l:L1diff}(1),  for all $|u|>\xi_1 \sT$,
\begin{align}\label{midrange}
&\text{Re} \; t\sM(t,\sigma,u)  \ge  t\int_0^{\sT/(\sigma |u|)} \big(1-\cos \frac{\sigma us}{\sT}\big)e^{-\sigma s} \nu(s)ds \nn\\
& \ge  t\frac{\cos 1}{2} \frac{\sigma^2u^2}{\sT^2}e^{-\sT/|u|}\int_0^{\sT/(\sigma |u|)}  s^2\nu(s)ds \ge c_6e^{-1/\xi_1} t \frac{\sigma^2u^2}{\sT^2} |\phi''(\sigma |u|/\sT)|.
\end{align}
In the first inequality above, we used the fact that $1-\cos x \ge \frac{\cos 1}{2} x^2$ for all $|x|\le 1$.
It follows that
\begin{align}\label{e:eq31}
\left|\int_{\xi_1\sT < |u| \le \xi_0\sT_0} e^{-t \sM(t,\sigma,u)}du \right| &\le 2 \xi_0 \sT \max_{\xi_1\sT < u \le \xi_0\sT}\exp\big(-c_6e^{-1/\xi_1} t \frac{\sigma^2u^2}{\sT^2} |\phi''(\sigma u/\sT)|\big)\nn\\
&\le 2\xi_0 \sT \exp\big(-c_6e^{-1/\xi_1}t\xi_1^2 \sigma^2  |\phi''(\sigma\xi_0)|\big) \nn\\
& \le 2\xi_0 \sT \exp \big(-c_7 t\sigma^2 |\phi''(\sigma)|\big) = 2\xi_0 \sT \exp \big(-c_7 \sT^2\big).
\end{align}
We used the fact that $\sT_0=\sT$ under the assumption $\sigma>\sigma_0$ in the first inequality and Lemma \ref{l:L1diff}(2) in the third inequality.

\smallskip

Finally, by the triangle inequality and inequalities \eqref{e:eq1}, \eqref{e:eq2} and \eqref{e:eq31}, we obtain
\begin{align}\label{tri}
&\left|\int_{\R} (e^{-t \sM(t,\sigma,u)}-e^{-\frac{1}{2}u^2}) du \right| \nn\\
&\le \left|\int_{|u| \le \xi_1 \sT} (e^{-t \sM(t,\sigma,u)}-e^{-\frac{1}{2}u^2}) du \right|  + \left|\int_{|u|>\xi_1 \sT} e^{-t \sM(t, \sigma, u)} du \right| +  \left|\int_{|u|>\xi_1 \sT} e^{-\frac{1}{2}u^2} du \right| \nn\\
& \le \frac{c_5}{\sT} + 2\xi_0 \sT \exp\big(-c_7\sT^2\big) + 2\int_{\xi_0 \sT}^\infty \big( 1+ \frac{u}{\sT} \big)^{-(2T_0+\delta)/(2T_0)} du  +  2\int_{\xi_1 \sT}^\infty e^{-\frac{1}{2}u^2} du \nn\\
& \to 0 \qquad \text{as} \quad \sT \to \infty.
\end{align}
This proves \eqref{CLT}.

\smallskip
\noindent
(Case 2): Assume that $\sigma \le \sigma_0$. We follow the proof given in (Case 1). First, using Lemma \ref{l:L2diff}(3) instead of Lemma \ref{l:L1diff}(3),  \eqref{e:eq2} still hold with possibly different constants $\xi_1$ and $c_5$. 

Next, note that $\sigma |u|/\sT\le \xi_0\sigma_0$ for $|u| \le \xi_0\sT_0$ in this case. Hence, by Lemma \ref{l:L2diff}(1), we see that \eqref{midrange} holds for all $|u| \le \xi_0\sT_0$ with a different constant $c_6$. Also, by Lemma \ref{l:largeH}(4), we have that $\sT^2 \asymp tH(\sigma)$ and $\sigma^2u^2 \sT^{-2} |\phi''(\sigma|u|/\sT)| \asymp H(\sigma |u|/\sT)$ for all $|u| \le \xi_0/\sT_0$.

Then, by \eqref{midrange} and Lemma \ref{l:largeH}(3) and (4), we have that for $\alpha_3':=\alpha_3 \wedge (3/2)$, 
\begin{align}\label{e:eq32}
&\left|\int_{\xi_1\sT < |u| \le \xi_0\sT_0} e^{-t \sM(t,\sigma,u)}du \right| \le 2 \int_{\xi_1 \sT}^{\xi_0\sT_0} \exp \big(- c_8 t H(\sigma u/\sT)\big)du\nn\\
& \le 2 \int_{\xi_1 \sT}^{\xi_0\sT_0} \exp\big(-c_9 t H(\sigma)u^{\alpha_3'}\sT^{-\alpha_3'}\big)du = 2 \sT \int_{\xi_1}^{\xi_0\sT_0/\sT} \exp\big(-c_9 t H(\sigma)u^{\alpha_3'}\big)du  \nn\\
&\le c_{10}\sT\exp\big(-\frac{c_9}{2} t H(\sigma)\xi_1^{\alpha_3'}\big)\int_{\xi_1}^{\infty} (tH(\sigma)u^{\alpha_3'})^{-2/\alpha_3'}du \le c_{11} \sT^{1-4/\alpha_3'}\exp \big( - c_{12} \sT^2 \big).
\end{align}
We used the change of the variables in the first equality and the fact that there exists a constant $c>0$ such that $e^{-x/2} \le cx^{-2/\alpha_3'}$ for all $x>\xi_1$ in the third inequality.

Using \eqref{e:eq32} instead of \eqref{e:eq1}, we see that \eqref{tri} still valid . Hence, we obtain \eqref{CLT}.

\vspace{2mm}

To complete the proof, we further assume that $T_0=0$. Choose any $0<T \le 2c_1^{-1}c_3$. To prove the second assertion, it suffices to show that there exist constants $c_{12}>1$ and $M_0 \ge c_2^{-1}c_3^2$ such that for all $t \in [T, 2c_1^{-1}c_3]$ and $x \in (0, tb(t/M_0)]$, 
\begin{equation}\label{claimdrag}
c_{12}^{-1} \le \int_{-\infty}^{\infty}e^{-t \sM(t, \sigma, u)} du \le c_{12},
\end{equation}
in view of \eqref{Mellin}. 
Note that \eqref{logdecay} is still valid with possibly different constants $\eps, \sigma_0$ and $\xi_0$. Hence, we have $\int_{-\infty}^{\infty}e^{-t \sM(t, \sigma, u)} du \in \R$ for all $t \in [T, 2c_1^{-1}c_3]$ and $\sigma>0$. Also, since inequalities \eqref{e:eq1}, \eqref{e:eq2}, \eqref{e:eq31} and \eqref{e:eq32} still work, by a similar argument to \eqref{tri}, we see that there exists a constant $c_{13}>0$ such that if $\sT=\sigma \sqrt{t(-\phi''(\sigma))} \ge c_{13}$, then \eqref{claimdrag} holds. Hence, it remains to prove that for a set  $A:=\{(t,\sigma) : t \in [T, 2c_1^{-1}c_3], \sigma>0, \sT < c_{13} \}$,
\begin{equation}\label{compact}
\inf_{(t,\sigma) \in A} \int_{-\infty}^{\infty}e^{-t \sM(t, \sigma, u)} du \asymp \sup_{(t,\sigma) \in A}\int_{-\infty}^{\infty}e^{-t \sM(t, \sigma, u)} du \asymp 1.
\end{equation}

Recall that $\sT^2 \ge c_2M_0$ if $\sigma \le R_2^{-1}$. By taking $M_0=c_2^{-1}c_{13}^2$, we have
$A \subset [T,2c_1^{-1}c_3] \times [R_2^{-1},\infty).$ On the other hand, since $T_0=0$, we have 
\begin{equation}\label{t0}
\lim_{\sigma \to \infty} \sigma^2(-\phi''(\sigma)) \ge \lim_{\sigma \to \infty} e^{-1}\sigma^2\int_0^{1/\sigma} s (s\nu(s)) ds \ge (2e)^{-1} \liminf_{\sigma \to \infty} \inf_{0<s<1/\sigma} (s\nu(s)) = \infty.
\end{equation}
Thus, there exists a constant $\sigma_1>0$ such that $\sT^2 \ge T\sigma^2(-\phi''(\sigma)) \ge c_{13}^2$ for all $\sigma >\sigma_1$ and hence $A \subset [T,2c_1^{-1}c_3] \times [R_2^{-1},\sigma_1]=:A_0$. Clearly, $(t, \sigma) \mapsto \int_{-\infty}^{\infty} e^{-t\sM(t,\sigma,u)}du$ is a continuous function on $A_0$. Therefore, we deduce \eqref{compact} from the extreme value theorem.
\qed

As a consequence, we obtain the following corollary.

\begin{cor}\label{exactconv}
	Suppose that the condition {\bf (L-3)} holds. Then, for every $N>0$,
	\begin{equation}\label{tlimit}
	\lim_{t \to \infty}p(t,x)  \sqrt{t(-\phi''(\sigma))}\exp\big(tH(\sigma)\big) = (2\pi)^{-1/2} \quad \text{uniformly in} \;\; x \in (0, N].
	\end{equation}
	
	If we further assume that the constant $T_0=0$ in the condition {\bf (E)}, then for every $N>0$,
	\begin{equation}\label{xlimit}
	\lim_{x \to 0}p(t,x)  \sqrt{t(-\phi''(\sigma))}\exp\big(tH(\sigma)\big) = (2\pi)^{-1/2} \quad \text{uniformly in} \;\; t \in [N, \infty).
	\end{equation}
\end{cor}
\proof
 Let $\sT=\sigma\sqrt{t(-\phi''(\sigma))}.$ Fix the constant $\sigma_0$ satisfying \eqref{e:sig1} and \eqref{e:sig2} with $\delta=1$.
 
 Observe that $\sigma=(\phi')^{-1}(x/t) \to \infty$ as $t \to \infty$. Hence, there exists a constant $t_N>T_0+1$ such that $\sigma> \sigma_0$ for all $t>t_N$ and $x \in (0,N]$.
 As we observed 
in the the proof of Proposition \ref{lefttail:L},  
 $\sT^2 \ge c_1t/2$ if $\sigma>\sigma_0$. 
 Since
 Lemma \ref{l:L1diff} holds under the condition {\bf (L-3)} only, we can use it and follow (Case 1) in the proof of Proposition \ref{lefttail:L}. Thus, \eqref{tri} holds for all $x \in (0,N]$ if $t>t_N$ (so that $\sigma>\sigma_0$). Since $\lim_{t \to \infty} \sT \ge \lim_{t \to \infty} (c_1t/2)^{1/2} = \infty$, \eqref{tlimit} follows from \eqref{Mellin} and \eqref{CLT}. 
 
 Now, we further assume that $T_0=0$. Since $\sigma$ also go to infinity as $x \to 0$, there exists a constant $x_N>0$ such that $\sigma> \sigma_0$ for all $t \ge N$ and $x \in (0, x_N)$. Hence, \eqref{tri} holds for all $t \ge N$ if $x<x_N$. Moreover, by \eqref{t0}, we get that $\lim_{x \to 0} \sT = \infty$ uniformly in $t \in [N, \infty)$ since $t \mapsto \sigma$ is increasing. Therefore, we also deduce \eqref{xlimit} from \eqref{Mellin} and \eqref{CLT}.
 \qed

A similar result to Corollary \ref{exactconv} is obtained in \cite[Section 3]{GLT}. 
Note that since the condition {\bf (L-3)} is very mild and do not require any lower scaling assumptions, our result covers geometric stable subordinators and Gamma subordinators which are not covered in \cite[Corollary 3.6]{GLT}.

\begin{example}\label{example}
	{\rm
Let $0<\alpha \le 1$ and $S_t$ be a geometric $\alpha$-stable subordinator whose Laplace exponent is given by $\log(1+\lambda^\alpha)$. When $\alpha=1$, $S_t$ is called a Gamma subordinator in the literature. It is known that the density of the L\'evy measure $\nu(x)$ and the transition density $p(t,x)$ of $S_t$ are equal to
\begin{align}\label{exactden}
\nu(x)=\frac{\alpha}{x} \sum_{n=0}^\infty (-1)^n \frac{x^{\alpha n}}{\Gamma(1+\alpha n)} \quad \text{and} \quad p(t,x)=\frac{x^{\alpha t-1}}{\Gamma(t)}\sum_{n=0}^\infty (-1)^n \frac{\Gamma(t+n)x^{\alpha n}}{n! \Gamma( \alpha t+ \alpha n)},
\end{align}
where $\Gamma(t)=\int_0^\infty y^{t-1}e^{-y}dy$ is the Gamma function. (See, \cite[Example 5.11]{BBKRSV}, \cite[Section 9.2]{GKMR} and \cite[Theorem 4.2]{Pil}.)
Then, we can see that the subordinator $S_t$ satisfies the conditions {\bf (E)} with $T_0=1/\alpha$, and {\bf (L-3)}. Thus, we can apply \eqref{tlimit}. (We can not expect that \eqref{xlimit} holds since $T_0>0$.) Below, we get the exact asymptotic behavior of $p(t,x)$ given in \eqref{exactden} as $t \to \infty$, from \eqref{tlimit}.

\smallskip

Observe that $\phi'(r) = \alpha r^{\alpha - 1}/(1+ r^{\alpha})$. Hence, for every $\lambda<\phi'(0)$, we can see that
\begin{align}\label{phiinverse}
(\phi')^{-1}(\lambda) + (\phi')^{-1}(\lambda)^{1-\alpha} = \frac{\alpha}{\lambda}.
\end{align}

Fix $N>0$. By \eqref{phiinverse}, for every $t>N/\phi'(1)$ and $x \in (0,N]$, $\sigma \in (1, \infty)$ is determined by 
\begin{align}\label{exactsigma}
 \sigma + \sigma^{1-\alpha} = \frac{\alpha  t}{x}.
\end{align}
We claim that
\begin{equation}\label{unif1}
\lim_{t \to \infty} \sigma x/t  = \alpha, \qquad \text{uniformly in} \;\; x\in (0,N].
\end{equation}
Indeed, \eqref{exactsigma} implies that $\sigma \le \alpha t/x$ and hence
$\sigma^{1-\alpha}x/t \le \alpha^{1-\alpha} x^\alpha t^{-\alpha} \le \alpha^{1-\alpha}N^\alpha t^{-\alpha}$. It follows that
$|\sigma x /t - \alpha| = |\sigma^{1-\alpha} x/t| \to 0$ as $t \to \infty$ uniformly in $x \in (0,N]$.
Then, since
\begin{align*}
-\phi''(\sigma) = \frac{\alpha(1-\alpha)\sigma^{\alpha-2}}{1+\sigma^\alpha} + \frac{\alpha^2 \sigma^{2\alpha-2}}{(1+\sigma^\alpha)^2}=  \frac{(1-\alpha)}{\sigma}\phi'(\sigma) + \phi'(\sigma)^2 = \frac{x^2}{t^2} \left(\frac{(1-\alpha)t}{x \sigma} +1 \right),
\end{align*}
we get from \eqref{unif1} that
\begin{align}\label{phi''}
\lim_{t \to \infty} \frac{ x^2/(\alpha t)}{t(-\phi''(\sigma))} = \lim_{t \to \infty} \frac{ 1/\alpha}{(1-\alpha)t/(x\sigma)+1}= \lim_{t \to \infty} \frac{ 1/\alpha}{(1-\alpha)/\alpha+1} =1, \;\;\; \text{uniformly in} \;\; x\in (0,N].
\end{align}
On the other hand, we get from \eqref{exactsigma} that
\begin{align}\label{expH}
&\exp\big(tH(\sigma)\big)= \exp\big(t\log(1+\sigma^\alpha) - t\sigma\phi'(\sigma)\big) = (1+\sigma^\alpha)^{t}e^{-\sigma x}   \nn\\
&\quad = (\sigma^{1-\alpha}+\sigma)^{\alpha t} (\sigma^{1-\alpha}+\sigma)^{(1-\alpha) t}\sigma^{-(1-\alpha)t}e^{-\sigma x}  =(\alpha t/x)^{\alpha t} (1 +\sigma^{-\alpha})^{(1-\alpha)t}e^{-\sigma x}.
\end{align} 
Therefore, according to \eqref{tlimit}, \eqref{phi''} and \eqref{expH},  it holds that
\begin{align*}
\lim_{t \to \infty} p(t,x)x^{1-\alpha t} (\alpha t)^{\alpha t-1/2}(1 +\sigma^{-\alpha})^{(1-\alpha)t}e^{-\sigma x}= (2\pi)^{-1/2}, \qquad \text{uniformly in} \;\; x\in (0,N],
\end{align*}
which is equivalent to (by Stirling's formula,)
\begin{align}\label{exactasymp0}
\lim_{t \to \infty} p(t,x)x^{1-\alpha t} \Gamma(\alpha t)(1 +\sigma^{-\alpha})^{(1-\alpha)t}e^{\alpha t-\sigma x} = 1,  \qquad \text{uniformly in} \;\; x\in (0,N],
\end{align}
where $\sigma =(\phi')^{-1}(x/t)\in (1, \infty)$ is determined by \eqref{exactsigma}. In other words, according to \eqref{exactden},
\begin{align*}
\lim_{t \to \infty} (1 +\sigma^{-\alpha})^{(1-\alpha)t}e^{\alpha t-\sigma x}\sum_{n=0}^\infty (-1)^n \frac{ \Gamma(t+n)\Gamma(\alpha t)x^{\alpha n}}{n!\Gamma(t) \Gamma( \alpha t+ \alpha n)}  = 1, \qquad \text{uniformly in} \;\; x\in (0,N].
\end{align*}

\smallskip

 Now, we express  \eqref{exactasymp0} in terms of $t$ and $x$ only.
  Let $\zeta = (x/(\alpha t))^\alpha$ and define
  $$\eta_{1}(\lambda) = 1 - \lambda^{1/\alpha}(\phi')^{-1}(\alpha \lambda^{1/\alpha}), \quad \eta_{2}(\lambda) = 1 + (\phi')^{-1}(\alpha \lambda^{1/\alpha})^{-\alpha} \quad \text{for} \;\; \alpha\lambda^{1/\alpha} \le \phi'(1).$$ 
   Observe that for $(t,x) \in (N/\phi'(1), \infty) \times (0,N]$, we have $\zeta \le (\phi'(1)/\alpha)^\alpha = 2^{-\alpha}$, $\eta_1(\zeta)= 1 - \sigma x/(\alpha t)$ and $\eta_2(\zeta) = 1 + \sigma^{-\alpha}$. Since $\phi$ is a Bernstein function, $\eta_1$ and $\eta_2$ are infinitely differentiable. 
   Moreover, according to \eqref{phiinverse}, it holds that for all $\alpha \lambda^{1/\alpha} \le \phi'(1)$,
   $$
   \lambda^{-1/\alpha}(1-\eta_1(\lambda)) + [\lambda^{-1/\alpha}(1-\eta_1(\lambda))]^{1-\alpha} = \alpha/(\alpha \lambda^{1/\alpha}),
   $$
   which is equivalent to
   $$
   (1-\eta_1(\lambda))+ \lambda (1-\eta_1(\lambda))^{1-\alpha}=1.
   $$
   Thus, $\eta_1(0)=0$ and
   $
   \eta_1'(\lambda) = [1-\eta_1(\lambda)]/[(1-\eta_1(\lambda))^\alpha + (1-\alpha) \lambda]
   $
   by the implicit function theorem.  Moreover, one can show the following by induction: for every $j \ge 1$, there exists an infinitely differentiable function $F_j$ such that $F_j(0), F_j'(0)<\infty$ and for all $\alpha \lambda^{1/\alpha} \le \phi'(1)$,
   $$
   \eta_1^{(j)}(\lambda) = \frac{F_j(\lambda)}{\left((1-\eta_1(\lambda))^\alpha + (1-\alpha) \lambda\right)^{2^j}},
   $$
   where $\eta_1^{(j)}$ denotes the $j$-th derivative of $\eta_1$. In particular $\eta_1^{(j)}(0)<\infty$ for all $j \ge 0$.
   
    Besides, we also get from \eqref{phiinverse} that
   $
   \eta_2(\lambda) (\phi')^{-1}(\alpha \lambda^{1/\alpha}) = \alpha/(\alpha \lambda^{1/\alpha})
   $
   and hence
   $$
   \eta_2(\lambda) = \frac{1}{\lambda^{1/\alpha}(\phi')^{-1}(\alpha \lambda^{1/\alpha})} = \frac{1}{1-\eta_1(\lambda)}.
   $$
   In particular, we can see that $\eta_2^{(j)}(0)<\infty$ for all $j \ge 0$. 
   
    Then, by the Taylor series expansion, there exist sequences of bounded functions $\{\eps_{j,1}\}_{j \ge 1}$ and $\{\eps_{j,2}\}_{j \ge 1}$ on $(0, 2^{-\alpha}]$ such that for every $k \ge 1$ and $(t,x) \in (N/\phi'(1),\infty) \times (0,N]$,
   \begin{align}\label{taylor1}
   \alpha t - \sigma x&= \alpha t \sum_{j=1}^k \frac{\eta_1^{(j)}(0)}{j!} \left(\frac{x}{\alpha t}\right)^{j\alpha} + \alpha t\eps_{k+1,1}(\zeta)\left(\frac{x}{\alpha t}\right)^{(k+1)\alpha}, \\
   1+\sigma^{-\alpha}&= 1+ \sum_{j=1}^k \frac{\eta_2^{(j)}(0)}{j!} \left(\frac{x}{\alpha t}\right)^{j\alpha} + \eps_{k+1,2}(\zeta)\left(\frac{x}{\alpha t}\right)^{(k+1)\alpha}. \nn
   \end{align}
   Moreover, from the uniqueness of the Taylor series, there exist a unique sequence $\{\delta_j\}_{j \ge 1}$ and a sequence of bounded functions $\{\eps_{j,3}\}_{j \ge 1}$ on $(0, 2^{-\alpha}]$ such that for every $k \ge 1$ and $(t,x) \in (N/\phi'(1),\infty) \times (0,N]$,
   \begin{align}\label{taylor2}
   1 + \sigma^{-\alpha}=\left(1+ \eps_{k+1,3}(\zeta)\left(\frac{x}{\alpha t}\right)^{(k+1)\alpha}\right) \exp \left( \sum_{j=1}^k \delta_j \left(\frac{x}{\alpha t}\right)^{j\alpha}  \right).
   \end{align}
   For example, we can calculate that $\eta_2'(0)=1$ and $\eta_2''(0)=2 \alpha$. It follows that
   \begin{align*}
   1 + \sigma^{-\alpha}&= 1 + \left(\frac{x}{\alpha t}\right)^{\alpha} + \alpha \left(\frac{x}{\alpha t}\right)^{2\alpha} + \eps_{3,2}(\zeta)\left(\frac{x}{\alpha t}\right)^{3\alpha}\\
   & = 1 +  \delta_1 \left(\frac{x}{\alpha t}\right)^{\alpha} +  (\delta_1^2/2 + \delta_2) \left(\frac{x}{\alpha t}\right)^{2\alpha}  + O\left(\left(\frac{x}{\alpha t}\right)^{3\alpha}\right)\\
   &=\left(1+ \eps_{3,3}(\zeta)\left(\frac{x}{\alpha t}\right)^{3\alpha}\right) \exp \left( \delta_1 \left(\frac{x}{\alpha t}\right)^{\alpha} +  \delta_2 \left(\frac{x}{\alpha t}\right)^{2\alpha}  \right),  
   \end{align*}
   and hence we see that $\delta_1= 1$ and $\delta_2 = \alpha -1/2$.

Note that for every $k \ge 1$ and $(t,x) \in (N/\phi'(1),\infty) \times (0,N]$,
$$\big|\alpha t \eps_{k+1,1}(\zeta) \zeta^{k+1}\big| \le \alpha^{1-(k+1)\alpha}N^{(k+1)\alpha}\lVert \eps_{k+1,1} \rVert_{\infty}t^{1-(k+1)\alpha}$$
 and 
$$
\left|\left(1+ \eps_{k+1,3}(\zeta)\zeta^{k+1}\right)^{(1-\alpha)t}-1 \right| \le c (1-\alpha)\alpha^{-(k+1)\alpha}N^{(k+1)\alpha} \lVert \eps_{k+1,3} \rVert_{\infty}t^{1-(k+1)\alpha}.
$$
Therefore, from \eqref{taylor1} and \eqref{taylor2}, we get that 
$$
\lim_{t \to \infty} (1+\sigma^{-\alpha})^{(1-\alpha)t}e^{\alpha t - \sigma x} \exp \left(-t \sum_{j=1}^{\lfloor 1/\alpha \rfloor} \bigg(\frac{\alpha\eta_1^{(j)}(0)}{j!}+(1-\alpha)\delta_j \bigg)\left(\frac{x}{\alpha t}\right)^{j\alpha} \right)=1,
$$
uniformly in $x \in (0,N]$.

Finally, according to \eqref{exactasymp0}, we conclude that 
\begin{align}\label{generalform}
&\lim_{t \to \infty} p(t,x)x^{1-\alpha t}\Gamma(\alpha t)\exp \left(t \sum_{j=1}^{\lfloor 1/\alpha \rfloor} \bigg(\frac{\alpha\eta_1^{(j)}(0)}{j!}+(1-\alpha)\delta_j \bigg)\left(\frac{x}{\alpha t}\right)^{j\alpha} \right)\nn\\
=&\lim_{t \to \infty} \exp \left(t \sum_{j=1}^{\lfloor 1/\alpha \rfloor} \bigg(\frac{\alpha\eta_1^{(j)}(0)}{j!}+(1-\alpha)\delta_j \bigg)\left(\frac{x}{\alpha t}\right)^{j\alpha} \right)\sum_{n=0}^\infty (-1)^n \frac{ \Gamma(t+n)\Gamma(\alpha t)x^{\alpha n}}{n!\Gamma(t) \Gamma( \alpha t+ \alpha n)}  = 1,
\end{align}
uniformly in $x \in (0,N]$.

In particular, we can check that $\eta_1'(0)=1$, $\eta_1''(0)=-2(1-\alpha)$, $\delta_1=1$ and $\delta_2=\alpha-1/2$. From these calculations, we obtain the following two special results:
if $\alpha \in (1/2,1]$, then
\begin{align}\label{exactasymp1}
&\lim_{t \to \infty} p(t,x)x^{1-\alpha t}\Gamma(\alpha t)\exp \bigg( t\left(\frac{x}{\alpha t}\right)^{\alpha}  \bigg) \nn\\
=&\lim_{t \to \infty} \exp \bigg( t\left(\frac{x}{\alpha t}\right)^{\alpha}  \bigg)\sum_{n=0}^\infty (-1)^n \frac{ \Gamma(t+n)\Gamma(\alpha t)x^{\alpha n}}{n!\Gamma(t) \Gamma( \alpha t+ \alpha n)}  = 1, \quad \text{uniformly in} \;\; x\in (0,N],
\end{align}
and if $\alpha \in (1/3, 1/2]$, then
\begin{align}\label{exactasymp2}
&\lim_{t \to \infty} p(t,x)x^{1-\alpha t}\Gamma(\alpha t)\exp \left( t\left(\frac{x}{\alpha t}\right)^{\alpha} - \frac{1-\alpha}{2}t \left(\frac{x}{\alpha t}\right)^{2\alpha} \right)\nn\\
=&\lim_{t \to \infty} \exp \left( t\left(\frac{x}{\alpha t}\right)^{\alpha} - \frac{1-\alpha}{2}t \left(\frac{x}{\alpha t}\right)^{2\alpha} \right)\sum_{n=0}^\infty (-1)^n \frac{ \Gamma(t+n)\Gamma(\alpha t)x^{\alpha n}}{n!\Gamma(t) \Gamma( \alpha t+ \alpha n)}  = 1,
\end{align}
uniformly in $x \in (0,N]$.

\smallskip

Since for each fixed $n \ge 0$, 
$
\lim_{t \to \infty} (\alpha t)^{\alpha n} \Gamma(\alpha t)/\Gamma(\alpha t + \alpha n)=\lim_{t \to \infty} t^n \Gamma(t)/\Gamma(t+n) = 1,
$
one may expect that for all sufficiently large $t$,
\begin{align*}
\sum_{n=0}^\infty (-1)^n \frac{\Gamma(t+n)\Gamma(\alpha t)x^{\alpha n}}{n! \Gamma(t) \Gamma( \alpha t+ \alpha n) }   \sim \sum_{n=0}^\infty (-1)^n \frac{t^n}{(\alpha t)^{\alpha n}}  \frac{x^{\alpha n}}{n!} = \exp \bigg(-  t\left(\frac{x}{\alpha t}\right)^\alpha \bigg).
\end{align*}
However, \eqref{generalform} says that this heuristic only works when $\alpha>1/2$. \qed

\vspace{2mm}

}
\end{example}

\vspace{3mm}

\subsection{Estimates on the transition density near the maximum value} In this subsection, we obtain maximum estimates on $p(t,x)$. Then, we extend the left tail estimates obtained in Section 3.1 as a corollary.

\begin{lemma}\label{fourier}
	Let $a \in [0,\infty)$, $\beta, c_1>0$ be constants and $f$ be a non-negative, non-decreasing function. Assume that 
	\begin{equation}\label{assu}
	\frac{f(R)}{f(r)} \ge c_1\left(\frac{R}{r}\right)^\beta \qquad \text{for all} \;\; a \le r \le R \;\; (\text{resp.} \;\; 0 \le r \le R \le a).
	\end{equation}
	Then, for every $c_2>0$, there exists a constant $c_3>0$ such that
	\begin{align*}
	&\int_a^\infty \exp(-c_2tf(\xi))d\xi \le c_3f^{-1}(1/t) \quad \text{for all} \;\; t \in (0,1/f(a)), \\
	(\text{resp.} \;\; &\int_0^a \exp(-c_2tf(\xi))d\xi \le c_3f^{-1}(1/t) \quad \text{for all} \;\; t \in [1/f(a),\infty),
	\end{align*}
	where $f^{-1}(s):= \inf \{r \ge 0: f(r) > s\}$ with a convention that $\inf \emptyset = \infty$.
\end{lemma}
\proof 
We first assume that \eqref{assu} holds for $a \le r \le R$. Note that $f^{-1}(1/t) \ge a$ for all $t \in (0,1/f(a))$. By the assumption, we get that for all $t \in (0, 1/f(a))$,
\begin{align*}
\int_a^\infty \exp(-c_2tf(\xi))d\xi &\le \int_a^{2f^{-1}(1/t)} \exp(-c_2tf(\xi))d\xi + \int_{2f^{-1}(1/t)}^\infty \exp\big(-c_2\frac{f(\xi)}{f(2f^{-1}(1/t))}\big)d\xi \\
& \le 2f^{-1}(1/t) + \int_{2f^{-1}(1/t)}^\infty \exp\big(-c_1c_2\big(\frac{\xi}{2f^{-1}(1/t)}\big)^\beta\big)d\xi\\
& = 2f^{-1}(1/t)\left(1 + \int_{1}^\infty \exp\big(-c_1c_2u^\beta\big)du \right)=c_3f^{-1}(1/t).
\end{align*}

On the other hand, assume that \eqref{assu} holds for $0 \le r \le R \le a$. If $a \le 2f^{-1}(1/t)$, then there is nothing to prove. Hence, assume that $a>2f^{-1}(1/t)$. Then, for all $t \in [1/f(a), \infty)$,
\begin{align*}
\int_0^a \exp(-c_2tf(\xi))d\xi &\le \int_0^{2f^{-1}(1/t)} \exp(-c_2tf(\xi))d\xi + \int_{2f^{-1}(1/t)}^a \exp\big(-c_2\frac{f(\xi)}{f(2f^{-1}(1/t))}\big)d\xi \\
& \le 2f^{-1}(1/t) + \int_{2f^{-1}(1/t)}^a \exp\big(-c_1c_2\big(\frac{\xi}{2f^{-1}(1/t)}\big)^\beta\big)d\xi\\
& \le 2f^{-1}(1/t)\left(1 + \int_{1}^\infty \exp\big(-c_1c_2u^\beta\big)du \right)=c_3f^{-1}(1/t).
\end{align*}
\qed

\begin{prop}\label{hkesup}{\rm(cf. \cite[Theorems 3.1 and 3.10]{GS}.)}
	
	\noindent {\rm (1)} Suppose that the condition {\bf (S-1)} holds.  Then, for every $T>0$, there exists a constant $c>0$ such that for all $t \in (0, T]$, 
	\begin{equation}\label{e:hkesup}
	\sup_{x \in \R} p(t,x) \le c H^{-1}(1/t).
	\end{equation}
	
	\noindent {\rm (2)} Suppose that the condition {\bf (L-1)} holds.  Then, for every $T>T_0$, where $T_0$ is the constant in {\bf (E)}, there exists a constant $c>0$ such that \eqref{e:hkesup} holds for all $t \in [T, \infty)$.
\end{prop}
\proof
(1) By Lemma \ref{l:smallH}(4) (and Lemma \ref{l:largeH}(4) as well if $R_1=\infty$) and the Fourier inversion theorem, for every $t>0$ and $x \in \R$,
\begin{align*}
p(t,x) &= \frac{1}{2 \pi} \int_{\R} e^{-i \xi x}e^{-t \phi(-i\xi)} d\xi \le \frac{1}{2 \pi} \int_{\R} |e^{-i \xi x}e^{-t \phi(-i\xi)}| d\xi \\
& \le \frac{1}{ \pi} \int_{0}^\infty \exp \big( -t\int_0^\infty (1-\cos(\xi s)) \nu(s) ds \big)d\xi \\
& \le \frac{2R_1^{-1}}{\pi}+\frac{1}{ \pi} \int_{2R_1^{-1}}^{\infty} \exp \big( -\frac{\cos 1}{2}t \xi^2\int_0^{1/\xi} s^2 \nu(s) ds \big)d\xi \\
&\le \frac{2R_1^{-1}}{\pi}+\frac{1}{ \pi} \int_{2R_1^{-1}}^{\infty} \exp \big( -c_2t H(\xi) \big)d\xi.
\end{align*}
We used the fact that $1-\cos x \ge 2^{-1} (\cos 1)x^2$ for all $|x|\le 1$ in the third inequality.

By Lemmas \ref{l:smallH}(3) and \ref{fourier}, (and Lemma \ref{l:largeH}(3) as well if $R_1=\infty$), there is a constant $c_3>0$ such that
\begin{equation*}
 \int_{2R_1^{-1}}^{\infty} \exp \big( -c_2t H(\xi) \big)d\xi \le \int_{2R_1^{-1} \wedge H^{-1}(1/T)}^{\infty} \exp \big( -c_2t H(\xi) \big)d\xi \le c_3H^{-1}(1/t) \;\; \text{for all} \;\; t \in (0,T].
\end{equation*}
Since $H^{-1}(1/t) \ge H^{-1}(1/T)$ for $t \in (0,T]$, we see that \eqref{e:hkesup} holds.

\smallskip
\noindent
(2) Fix any $T' \in (T_0, T)$. By the proof of Proposition \ref{p:existence},  $\int_{\R} |e^{-T' \phi(-i\xi)}| d\xi \le c_4<\infty$. On the other hand, by the condition {\bf (E)}, there exists a constant $s_0>0$ such that $\nu(s) \ge 1/(2T_0s)$ for all $s \in (0, s_0]$.
Then, by the similar arguments as the ones given in the proof of (1) and using Lemma \ref{l:largeH} instead of Lemma \ref{l:smallH}, we get
\begin{align*}
p(t,x) & \le \frac{1}{2 \pi} \int_{|\xi| \le 1/s_0} |e^{-t \phi(-i\xi)}| d\xi +  \frac{1}{2 \pi} \int_{|\xi| >1/s_0} |e^{-(t-T') \phi(-i\xi)}||e^{-T' \phi(-i\xi)}| d\xi  \\
&\le
\frac{1}{ \pi} \int_0^{1/s_0} \exp \big( -c_5t H(\xi) \big)d\xi + \frac{c_4}{2\pi} \sup_{|\xi|>1/s_0}\big|\exp \big(-(t-T') 
\phi(-i\xi)\big)\big| =: I_1+I_2.
\end{align*}
By Lemmas \ref{l:largeH}(3) and \eqref{fourier},
$$I_1 \le \frac{1}{\pi} \int_0^{1/s_0 \vee H^{-1}(1/T)} \exp \big(-c_5 t H(\xi) \big) d \xi \le c_6 H^{-1}(1/t) \quad \text{for all} \;\; t \ge T.$$
On the other hand, we also have
\begin{align*}
I_2 &\le \frac{c_4}{2\pi} \sup_{|\xi|>1/s_0}\exp \big(- \frac{(T-T')\cos 1}{2T}t \xi^2
\int_0^{1/\xi} s^2 \nu(s)ds \big)\\
& \le \frac{c_4}{2\pi} \sup_{|\xi|>1/s_0}\exp \big(-c_7t \xi^2
\int_0^{1/\xi} s^2s^{-1}ds \big) = \frac{c_4}{2\pi} \exp(-\frac{c_7}{2}t).
\end{align*}

Note that the lower bound in Lemma \ref{l:largeH}(3) implies that there exists $c_8>0$ such that $H^{-1}(1/t) \ge c_8t^{-1/(\alpha_3 \wedge (3/2))}$ for all $t \ge T$. Since it also holds that $\exp(-c_7t/2) \le c_{9}t^{-1/(\alpha_3 \wedge (3/2))}$ for all $t \ge T$, for some constant $c_{9}>0$, we finish the proof.
\qed

Now, we find a range of $x$ which achieves the maximum value of $p(t,x)$. One of the important points in the following proposition is that $N$ can be arbitrarily big number. This point allows us to remove the constant $M_0$ in estimates in Corollary \ref{c:lefttail}. 

A similar result to the following proposition was established in \cite[Theorem 5.3]{GS} which considers a class of L\'evy processes whose L\'evy measure dominates some symmetric measure. Note that since the support of the L\'evy measure of a subordinator is one-sided, that is always contained in $(0,\infty)$, we can only push the $y$-variable to the positive direction in the following unlike \cite[Theorem 5.3]{GS}.

\begin{prop}\label{arbi}
\noindent {\rm (1)}
	Suppose that the condition {\bf (S-1)} holds. Then, for every $T>0$ and $N>0$, there exists a constant $c_1>1$ such that for all  $t \in (0, T]$, 
	\begin{equation}\label{ondiag}
	c_1^{-1} H^{-1}(1/t) \le p(t, tb(t/(2M_0))+y) \le c_1 H^{-1}(1/t) \quad \text{for all} \;\; 0\le y \le NH^{-1}(1/t)^{-1},
	\end{equation}
	where $M_0$ is the constant in Proposition \ref{lefttail}.
	
\noindent {\rm (2)} Suppose that the conditions {\bf (L-1)} and {\bf (L-3)} hold. Then, for every $N>0$, there is a comparison constant such that \eqref{ondiag} holds for all $t \in [2T_1, \infty)$ and $y \in [0,NH^{-1}(1/t)^{-1}]$ with the constants $T_1, M_0$ in Proposition \ref{lefttail:L}.

 Moreover, if $T_0=0$ in the condition {\bf (E)}, then for every $T>0$ and $N>0$, there is a comparison constant such that \eqref{ondiag} holds for all $t \in [T, \infty)$ and $y \in [0,NH^{-1}(1/t)^{-1}]$.
\end{prop}
\proof
By Proposition \ref{hkesup}, it remains to prove the lower bound.
Below, we give the full proof for (1) and explain main differences in the proof of (2) in the end.

 For $p \in [1,4]$, we observe that
$$b(t/(pM_0)) \le b(t/M_0) \qquad \text{and} \qquad ((\phi')^{-1} \circ b)(t/(pM_0)) = H^{-1}(pM_0/t).$$
 Hence, by Proposition \ref{lefttail} and Lemma \ref{l:smallH}(3) and (4),  for all $p \in [1,4]$,
\begin{equation}\label{gap}
p(t, tb(t/(pM_0))) \ge \frac{ce^{-pM_0}}{\sqrt{t|(\phi''\circ H^{-1})(pM_0/t)|}} \ge \frac{ce^{-4M_0} H^{-1}(pM_0/t)}{\sqrt{pM_0}} \ge c_2 H^{-1}(1/t).
\end{equation}
According to Lemma \ref{l:b}, there is a constant $c_3>0$ such that 
\begin{equation}\label{gapb}
tb(t/M_0) - tb(t/(4M_0)) \ge c_3H^{-1}(1/t)^{-1} \qquad \text{for all} \;\; t \in (0,T].
\end{equation}
Then, by the intermediate value theorem, for all $t \in (0,T]$ and $u \in [0,c_3H^{-1}(1/t)^{-1}]$, there exists $p \in [1,4]$ such that $tb(t/M_0)-u=tb(t/(pM_0))$. Hence, by \eqref{gap}, we get
\begin{equation}\label{max}
p(t, tb(t/M_0)-u) \ge c_2H^{-1}(1/t) \qquad \text{for all} \;\; t \in (0,T], \;\; u \in [0,c_3H^{-1}(1/t)^{-1}].
\end{equation}

\smallskip

By the semigroup property and \eqref{max}, we have that for any $t \in (0, T]$ and $y\ge 0$,
\begin{align*}
p(2t, 2tb(t/M_0)+y) &= \int_\R p(t, tb(t/M_0)-u)p(t, tb(t/M_0)+y+u)du\\
&\ge c_2H^{-1}(1/t) \P\big(y \le S_t - tb(t/M_0) \le  y+c_3H^{-1}(1/t)^{-1}\big).
\end{align*}
Thus, since $H^{-1}(1/t) \asymp H^{-1}(1/(2t))$ for all $t \in (0,T]$ by Lemma \ref{l:smallH}(3), it suffices to show that for every fixed $N>0$, it holds that
\begin{equation}\label{e:drag}
\inf_{t \in (0, T]}\inf_{y \in [0, NH^{-1}(1/t)^{-1}]}\P\big(y \le S_t-tb(t/M_0) \le  y+c_3H^{-1}(1/t)^{-1}\big) >0.
\end{equation}

Let $(t_n : n \ge 1)$ be a sequence of time variables realizing the infimum in \eqref{e:drag}. Since $(0,T]$ is a bounded interval, after taking a subsequence, we can assume that $t_n$ converges to $t_* \in [0, T]$. If $t_* \in (0, T]$, then since the support of the distribution of $S_{t_*}$ is $(0, \infty)$, we obtain \eqref{e:drag}. Hence, we assume that $t_* =0$ and all $t_n$ are sufficiently small.

\smallskip

Define $\nu_n(s):=\nu(s)\1_{(0, H^{-1}(1/t_n)^{-1}]}(s)$ and let $\tilde{S}_u$ be a subordinator without drift, whose L\'evy measure is given by $\nu_n(s)ds$. Then, for all $u>0$, $S_{u}=\tilde{S}_u + P_u$, $\P$-a.s. where $P$ is a compounded Poisson process whose L\'evy measure is given by $\nu(s)\1_{(H^{-1}(1/t_n)^{-1}, \infty)}(s)ds$. Thus, by \eqref{basicH2},
\begin{align*}
\P(\tilde{S}_{t_n} = S_{t_n}) &= \P(P_{t_n}=0) = \exp\big(-t_nw(\frac{1}{H^{-1}(1/t_n)})\big)\\
& \ge \exp \big(- \frac{e}{e-2}t_n (H\circ H^{-1})(1/t_n)\big) = \exp \big(-\frac{e}{e-2}\big).
\end{align*}
Hence, to prove \eqref{e:drag}, it is enough to show that
\begin{equation}\label{e:dragcut}
\liminf_{n\to \infty}\inf_{y \in [0, NH^{-1}(1/t_n)^{-1}]}\P\big(y \le \tilde{S}_{t_n}-t_nb(t_n/M_0) \le  y+c_3H^{-1}(1/t_n)^{-1}\big) >0.
\end{equation}

\smallskip

Define $Z_n=H^{-1}(1/t_n)(\tilde{S}_{t_n}-t_nb(t_n/M_0))$. Then, we have that for $\xi \in \R$,
\begin{align*}
\E[\exp(i\xi Z_n)] &= \exp\bigg(-i \xi t_nH^{-1}(1/t_n)b(t_n/M_0)\big)\E\big[\exp\big(i\xi H^{-1}(1/t_n)\tilde{S}_{t_n} \big)\big]\bigg)\\
& = \exp\left(-i \xi t_nH^{-1}(1/t_n)b(t_n/M_0) + t_n\int_0^{\infty} \big( \exp\big(i \xi H^{-1}(1/t_n) s\big) -1 \big) \nu_n(s)ds \right).
\end{align*}
Therefore, we get $\E[\exp(i \xi Z_n)] = \exp \big(\Psi_n(\xi)\big)$ for all $\xi \in \R$ and $n \ge 1$ where
\begin{align*}
\Psi_n(\xi) &= \int_0^\infty\left(e^{i\xi s}-1- \frac{i \xi s}{1+s^2} \right) \lambda_n(s) ds - i \xi \gamma_n,\\
\lambda_n(s) &= t_nH^{-1}(1/t_n)^{-1}\nu_n(H^{-1}(1/t_n)^{-1}s),\\
\gamma_n &= t_nH^{-1}(1/t_n)b(t_n/M_0)- \int_0^\infty \frac{s}{1+s^2} \lambda_n(s)ds,
\end{align*}
by the change of the variables.
We claim that the family of random variables $\{Z_n:n \ge 1\}$ is tight. Indeed, according to \cite[(3.2)]{JP}, it holds that for all $n \ge 1$ and $R>1$,
\begin{align*}
\P(Z_n \ge R) &\le c_4 \left(\int_0^\infty \left(\frac{s^2}{R^2} \wedge 1\right) \lambda_n(s)ds+ \frac{1}{R}\left| \gamma_n + \int_R^\infty \frac{s}{1+s^2}\lambda_n(s)ds - \int_0^R \frac{s^3}{1+s^2} \lambda_n(s)ds \right|\right) \\
&=: c_4(I_1+I_2).
\end{align*}

First, by the change of variables and \eqref{basicH1}, we have
\begin{align*}
I_1 &= t_n\int_0^\infty \left(\frac{H^{-1}(1/t_n)^{2}u^2}{R^2} \wedge 1 \right)  \nu_n(u)du = R^{-2}t_nH^{-1}(1/t_n)^2\int_0^{H^{-1}(1/t_n)^{-1}}u^2 \nu_n(u)du \\
& \le 2eR^{-2}t_nH(H^{-1}(1/t_n))= 2eR^{-2}.
\end{align*}

On the other hand, by the change of variables, Lemma \ref{l:smallH}(3) and (4) and \eqref{basicH2},
\begin{align*}
RI_2&=\left|t_nH^{-1}(1/t_n)b(t_n/M_0) - \int_0^R s\lambda_n(s)ds\right|  \\
&=t_nH^{-1}(1/t_n)\left|\int_0^\infty s\exp \big(-H^{-1}(M_0/t_n)s\big) \nu(s)ds - \int_0^{RH^{-1}(1/t_n)^{-1}} s \nu_n(s)ds \right| \\
&\le t_nH^{-1}(1/t_n) \bigg(\int_0^{H^{-1}(1/t_n)^{-1}} s\left(1-\exp \big(-H^{-1}(M_0/t_n)s\big)\right) \nu(s)ds \\
& \qquad \qquad \qquad \quad\;\;  + \int_{H^{-1}(1/t_n)^{-1}}^\infty s\exp\big(-H^{-1}(M_0/t_n)s\big) \nu(s)ds\bigg) \\
& \le t_nH^{-1}(1/t_n)\left(H^{-1}(M_0/t_n)\int_0^{H^{-1}(1/t_n)^{-1}} s^2 \nu(s)ds +H^{-1}(M_0/t_n)^{-1}w\big(H^{-1}(1/t_n)^{-1}\big)\right) \\
& \le c_5t_nH^{-1}(1/t_n) \left(H^{-1}(1/t_n)^{-1} (H\circ H^{-1})(1/t_n) \right) = c_5.
\end{align*}
We used the fact that the support of $\nu_n$ is contained in $(0, H^{-1}(1/t_n)^{-1}]$ in the first inequality, and the mean value theorem and the fact that for every $a>0$, $\sup_{x \in (0,\infty)} xe^{-ax}=e^{-1}a^{-1}$ in the second inequality.

Therefore, we deduce that $\P(Z_n \ge R) \le c_4(2e+c_5)R^{-1}$ for all $n \ge 1$ and $R>1$, 
which yields that the family $\{Z_n:n \ge 1\}$ is tight. Then by the Prokhorov's
theorem, by taking a subsequence, we can assume that $Z_n$ is weakly convergent to the random variable $Z_*$.

\smallskip

Now, from the weak convergence, we can prove \eqref{e:dragcut} by showing the following:
\begin{equation}\label{scaledrag}
\inf_{z \in [0,N]} \P(z \le Z_* \le z+c_3) >0.
\end{equation}

According to \cite[Theorem 8.7]{Sa}, $Z_*$ is a infinitely divisible random variable with the characteristic function
$$
\Psi_*(\xi) = -\frac{1}{2}A_*\xi^2 - i \xi \gamma_* + \int_0^\infty \left(e^{i \xi s} - 1 - \frac{i \xi s}{1+s^2} \lambda_*(s) ds \right),
$$
where the triplet $(A_*, \gamma_*, \lambda_*)$ is characterized by

(i) $\lim_{\eps \to 0} \limsup_{n \to \infty} \left| \int_0^\eps s^2 \lambda_n(s) ds - A_* \right| = 0$;

(ii) $\gamma_*= \lim_{n \to \infty} \gamma_n$;

(iii) $\int_0^\infty f(s) \lambda_*(s)ds = \lim_{n \to \infty} \int_0^\infty f(s) \lambda_n(s)ds$ for any bounded continuous function $f$ vanishing in a neighborhood of $0$.

\smallskip

If $A_*>0$, then it is evident that the support of $Z_*$ is $\R$ and hence \eqref{scaledrag} holds. Hence, we assume that $A_*=0$. Then, by (i) and (iii) in the above characterization, for every $\eta\in (0,1)$,
\begin{align*}
\int_0^\eta s^2 \lambda_*(s) ds &= \lim_{\eps \to 0+} \int_\eps^\eta s^2 \lambda_*(s)ds = \lim_{\eps \to 0+} \lim_{n \to \infty} \left(\int_0^\eta s^2 \lambda_n(s)ds - \int_0^\eps s^2 \lambda_n(s)ds \right) \\
& = \lim_{\eps \to 0+} \lim_{n \to \infty} \int_0^\eta s^2 \lambda_n(s)ds = \lim_{n \to \infty} t_n H^{-1}(1/t_n)^2 \int_0^{\eta H^{-1}(1/t_n)^{-1}} u^2 \nu(u) du \\
& \ge \lim_{n \to \infty} c_7\eta^2 t_nH(\eta^{-1} H^{-1}(1/t_n)) \ge c_7 \eta^2>0.
\end{align*}
We used Lemma \ref{l:smallH}(4) in the first inequality and the monotonicity of $H$ in the second inequality. It follows that according to \cite[Lemma 2.5]{Pi}, if $\int_0^1 s \lambda_*(s)ds = \infty$, then the support of $Z_*$ is $\R$ so that \eqref{scaledrag} holds. Assume that $\int_0^1 s \lambda_*(s)ds < \infty$. Then we see from (iii) in the characterization that $\limsup_{n \to \infty} \int_0^1 s \lambda_n(s)ds<\infty$. Hence, again by \cite[Lemma 2.5]{Pi}, the support of $Z_*$ is $[-\chi, \infty)$ where $\chi = \lim_{n \to \infty} t_nH^{-1}(1/t_n)b(t_n/M_0) \ge 0$. 

Since the support of $Z_*$ includes $(0, \infty)$ in any cases, we see that \eqref{scaledrag} holds. This finishes the proof of the proposition under the condition {\bf (S-1)}.

\smallskip

Hereafter, we assume that the conditions {\bf (L-1)} and {\bf (L-3)} hold instead of {\bf (S-1)}. By Proposition \ref{lefttail:L} and Lemmas \ref{l:largeH}(3$\&$4) and \ref{l:b}, we can follow \eqref{gap}, \eqref{gapb} and \eqref{max}, and hence it suffice to show that for every $N>0$, there exists a constant $c_8>0$ such that
\begin{equation}\label{e:dragL}
\inf_{t \in [T, \infty)}\inf_{y \in [0, NH^{-1}(1/t)^{-1}]}\P\big(y \le S_t-tb(t/M_0) \le  y+c_3H^{-1}(1/t)^{-1}\big) \ge c_8.
\end{equation}
(The constant $c_3$ may differ.) For convenience, we still denote a sequence of time variables realizing the infimum in \eqref{e:dragL} by $(t_n:n\ge 1)$. Then, after taking a subsequence, we can assume that either $t_n$ converges to $t^* \in [T, \infty)$ or $\lim_{n \to \infty}t_n =\infty$. If $t_n$ converges, then we are done. Hence, assume that $\lim_{n \to \infty}t_n =\infty$ and all $t_n$ are sufficiently large. 
Then, by using Lemma \ref{l:largeH} instead of Lemma \ref{l:smallH}, we can follow the proof under the condition {\bf (S-1)} and deduce the desired result. 

Furthermore, note that the restriction that $t \ge 2T_1$ is only required to obtain \eqref{gap}. Hence, in view of the second statement of Proposition \ref{lefttail:L}, we can see that the later assertion holds. This completes the proof.
\qed

Recall that $\sigma=\sigma(t,x)=(\phi')^{-1}(x/t)$. As a corollary to the above proposition, we can erase the constant $M_0$ in Propositions \ref{lefttail} and \ref{lefttail:L}.

\begin{cor}\label{c:lefttail}

\noindent {\rm (1)}
	Suppose that the condition {\bf (S-1)} holds. Then, for every fixed $T>0$ and $N>0$, it holds that for all $t \in (0, T]$,
\begin{align}\label{slefttail}
&p(t,x) \asymp \frac{1}{\sqrt{t(-\phi''(\sigma))}}\exp\big(-tH(\sigma)\big),  \qquad \text{for} \quad x \in (0, tb(t)],\nn\\
&p(t, tb(t)+y) \asymp H^{-1}(1/t), \quad\;\;\; \quad \qquad\qquad\text{for} \quad y \in [0, NH^{-1}(1/t)^{-1}].
\end{align}

\noindent {\rm (2)} 
Suppose that the conditions {\bf (L-1)} and {\bf (L-3)} hold. Then, for every $N>0$, there are comparison constants such that  for all $t \in [2T_1,\infty)$, \eqref{slefttail} holds for $x \in (0, tb(t)]$ and $y \in [0,NH^{-1}(1/t)^{-1}]$ where $T_1$ is the constant in Proposition \ref{lefttail:L}.

 Moreover, if $T_0=0$ in the condition {\bf (E)}, then for every $T>0$ and $N>0$, there is a comparison constant such that \eqref{slefttail} holds for all $t \in [T, \infty)$, $x \in (0, tb(t)]$ and $y \in [0,NH^{-1}(1/t)^{-1}]$.
\end{cor}
\proof
(1) The second comparison follows from Proposition \ref{arbi} and Lemma \ref{l:b}.
To prove the first one, in view of Proposition \ref{lefttail}, it suffices to consider for $x \in [tb(t/M_0), tb(t)]$. For those $x$, we see from Lemma \ref{l:b} and Proposition \ref{arbi} that $p(t,x) \asymp H^{-1}(1/t)$. On the other hand, 
 observe that for $x \in [tb(t/M_0), tb(t)]$,
$
\sigma = (\phi')^{-1}(x/t) \ge H^{-1}(1/t) \ge H^{-1}(1/T),
$
and hence by Lemma \ref{l:smallH}(3) and (4), 
\begin{align}\label{e:ondiag}
\frac{\exp\big(-tH(\sigma)\big)}{\sqrt{t(-\phi''(\sigma))}} \asymp \frac{\sigma}{\sqrt{tH(\sigma)}} \asymp \sigma \asymp H^{-1}(1/t) \asymp p(t,x).
\end{align}
This proves the corollary under the condition {\bf (S-1)}.

 \smallskip
 
 \noindent (2) Similarly, the second comparison follows from Proposition \ref{arbi} and Lemma \ref{l:b}. Moreover, in this case, for $x \in [tb(t/M_0), tb(t)]$, 
 $
 \sigma \le H^{-1}(M_0/t) \le H^{-1}(M_0/(2T_1)).
 $
 Therefore,  by Lemmas \ref{l:largeH}(3$\&$4) and \ref{l:b} and Proposition \ref{lefttail:L}, we can deduce that \eqref{e:ondiag} holds for all $t \in [2T_1, \infty)$ and $x \in [tb(t/M_0), tb(t)]$. Hence, we obtain the result by the same argument as the one for (1).
 
 Besides, in view of the second statement in Proposition \ref{lefttail:L}, the proof for the assertion under the assumption $T_0=0$ is exactly the same.
 This completes the proof.
\qed

\vspace{3mm}

\subsection{Estimates on right tail probabilities}
In this subsection, we establish estimates on $p(t,x)$ when $x \ge tb(t)$. Then, by combining with the results established in Section 3.2, we obtain full estimates on  $p(t,x)$.

\vspace{3mm}

Recall that we have $w^{-1}(2e/t) \le H^{-1}(1/t)^{-1}$ for all $t>0$,
$$
D(t):= t\max_{s \in [w^{-1}(2e/t), H^{-1}(1/t)^{-1}]} sH(s^{-1})
$$
and
\begin{equation*}
\theta(t,y):=
\begin{cases}
H^{-1}(1/t)^{-1}, \; &\text{if} \;\; y \in [0, H^{-1}(1/t)^{-1}),\\
\min\left\{s \in \big[w^{-1}(2e/t), H^{-1}(1/t)^{-1}\big] : t s H(s^{-1}) = y\right\}, \; &\text{if} \;\; y \in [H^{-1}\left(1/t\right)^{-1}, D(t)], \\
w^{-1}(2e/t), \quad &\text{if} \;\; y \in (D(t), \infty).
\end{cases}
\end{equation*}

\smallskip

We also recall that $\theta(t,y) \in [w^{-1}(2e/t), H^{-1}(1/t)^{-1}]$ for all $t>0$ and $y \ge 0$, and hence $\lim_{t \to 0} \theta(t,y)=0$ and $\lim_{t \to \infty} \theta(t,y) = \infty$  for each fixed $y \ge 0$.

\begin{lemma}\label{l:theta}
     It holds that
	$$t \theta(t,y) H(\theta(t,y)^{-1}) \le y \vee H^{-1}(1/t)^{-1} \qquad \text{for all} \;\; t>0,\; y \ge 0.$$
	In particular, for $y \in [H^{-1}(1/t)^{-1}, D(t)]$, it holds that
	$y =
	t \theta(t,y) H(\theta(t,y)^{-1}).
	$
\end{lemma}
\proof If $y \in [0, H^{-1}(1/t)^{-1}]$, then $t \theta(t,y) H(\theta(t,y)^{-1}) = t H^{-1}(1/t)^{-1} (H \circ H^{-1})(1/t) = H^{-1}(1/t)^{-1}$. Else if $y \in [H^{-1}(1/t)^{-1}, D(t)]$, then $t \theta(t,y) H(\theta(t,y)^{-1}) = y$. Otherwise, if $y > D(t)$, then $t \theta(t,y) H(\theta(t,y)^{-1}) < y$. \qed

The following theorem is the main result in this subsection.

\begin{thm}\label{righttail}

\noindent
{\rm (1)} Suppose that the condition {\bf (S)} holds. Then, for every $T>0$, there exist constants $c_1>1, c_2>0$ such that for all $t \in (0, T]$ and $y \in [0, R_1/2)$,
	\begin{align}\label{e:righttail}
	&c_1^{-1} H^{-1}(1/t) \min\left\{1, \frac{t \nu(y)}{H^{-1}(1/t)} + \exp \big(-\frac{c_2y}{\theta(t, y/(8e^2))} \big) \right\} \nn\\
	&\qquad \le p(t, tb(t)+y) \le c_1 H^{-1}(1/t) \min\left\{1, \frac{t \nu(y)}{H^{-1}(1/t)} + \exp \big(-\frac{y}{8\theta(t, y/(8e^2))} \big) \right\}.
	\end{align}

	Moreover, if the condition {\bf (S-3*)} further hold, then \eqref{e:righttail} holds true for all $t \in (0,T]$ and $y \in [0, \infty)$.

\noindent
{\rm (2)} Suppose that the condition {\bf (L)} holds. Then, there exist constants $c_1>1,c_2>0$ such that \eqref{e:righttail} holds for all $t \in [2T_1, \infty)$ and $y \in [0, \infty)$ where $T_1$ is the constant in Proposition \ref{lefttail:L}. Moreover, if $T_0=0$ in the condition {\bf (E)}, then for every $T>0$, there are comparison constants such that \eqref{e:righttail} holds for all $t \in [T, \infty)$ and $y \in [0, \infty)$.
\end{thm}
\proof
The result follows from Propositions \ref{mainupper} and \ref{p:mixedlow}.
\qed

\begin{prop}\label{mainupper}
Under the settings of Theorem \ref{righttail}, the upper bound in \eqref{e:righttail} holds.
\end{prop}
\proof
For convenience of notation, we let $\delta:=1/(8e^2)$. 

\noindent (1) We first assume that only the condition {\bf (S)} holds.
Fix $T>0$ and $t \in (0, T]$. If $\delta y \le H^{-1}(1/t)^{-1}$,
then
$
\exp \big(-y/\theta(t, \delta y) \big) \ge \exp \big(- 1/\delta \big)
$
and hence we obtain the upper bound in \eqref{e:righttail} from Proposition \ref{hkesup}. 
Therefore, 
for the remainder part of the proof of (1), we assume that $\delta y>H^{-1}(1/t)^{-1}$. 

 Define
\begin{align*}
&\nu_1(s):=  \1_{(0, \theta(t,\delta y)]} (s)\nu(s) \qquad \text{and} \qquad \nu_2(s):= \nu(s)-\nu_1(s) = \1_{(\theta(t,\delta y),\infty)}(s) \nu(s).
\end{align*}

Denote by $S^i$ and $H_i$ the corresponding subordinator and $H$-function with respect to the L\'evy measure $\nu_i$ for $i=1,2$, respectively. Since $\liminf_{s \to 0} s\nu_1(s) = \liminf_{s \to 0} s\nu(s) = 1/T_0$, by Proposition \ref{p:existence}, $S^1_u$ has a transition density function $p^1(u, \cdot )$ for every $u>T_0$. Recall that $T_0=0$ under the condition {\bf (S)}. Since $S_t = S^1_t + S^2_t$, we see that for $t>0$,
\begin{align}\label{semigroup}
&p(t,tb(t)+y) = \int_{\R} p^1(t, tb(t)+y-z) \P(S^2_t \in dz)\nn \\
&\quad= \int_{\{z \le y/4\}} p^1(t, tb(t)+y-z) \P(S^2_t \in dz) + \int_{\{z>y/4\}} p^1(t, tb(t)+y-z) \P(S^2_t \in dz) \nn\\
&\quad\le \sup_{z \ge 3y/4} p^1(t, tb(t)+z) + \sup_{z>y/4} \frac{\P(S_t^2 \in dz)}{dz}=:A_1 + A_2.
\end{align}

\vspace{2mm}

\noindent
\textit{Step1.}
First, we estimate $A_1$. By the semigroup property, for every $z \ge 3y/4$,
\begin{align}\label{semigroup1}
p^1(t, tb(t)+z) & = \left(\int_{u < z/2} + \int_{u \ge z/2}\right)  p^1\big(t/2, tb(t)/2 + u\big) p^1\big(t/2, tb(t)/2 + z- u\big) du \nn\\
& \le 2 \P\big(S^1_{t/2} \ge \frac{t}{2}b(t) + \frac{3y}{8}\big) \sup_{u \in \R} p^1(t/2,u).
\end{align}

We claim that there exists a constant $c_2>0$ such that for all $t \in (0,T]$,
\begin{equation}\label{p1supp}
\sup_{u \in \R} p^{1}(t/2,u) \le c_2 H^{-1}(1/t).
\end{equation}
Indeed, by \cite[Lemma 7.2]{CKW}, Proposition \ref{hkesup} and Lemma \ref{l:smallH}(3), we see that
$$
\sup_{u \in \R} p^{1}(t/2,u) \le \exp\big(2^{-1}t w(\theta(t, \delta y)) \big)\sup_{u \in \R} p(t/2,u) \le c\exp(e)H^{-1}(2/t) \le cH^{-1}(1/t).
$$

\smallskip

On the other hand, by the Chebychev's inequality and \cite[Lemma 2.5]{CK}, for every $\lambda>0$,
\begin{align*}
\P(S^1_{t/2} \ge \frac{t}{2}b(t)&+\frac{3y}{8}) \le \E\big[\exp(\lambda S^1_{t/2} - \frac{\lambda t}{2}b(t) - \frac{3 \lambda y}{8})\big] \\
&=e^{-3\lambda y/8} \exp \left(\frac{t}{2} \int_0^{\theta(t,\delta y)} (e^{\lambda s} - 1)\nu(s) ds - \frac{t}{2} \int_0^{\infty} \lambda se^{-H^{-1}(1/t)s} \nu(s) ds \right)\\
&\le e^{-3\lambda y/8} \exp \left(\frac{t}{2} \int_0^{\theta(t,\delta y)}(e^{\lambda s}- e^{-H^{-1}(1/t)s})  \lambda s\nu(s) ds \right) \\
&\le e^{-3\lambda y/8} \exp \left(\frac{t}{2} \int_0^{\theta(t,\delta y)} (\lambda + H^{-1}(1/t))  \lambda e^{\lambda s} s^2\nu(s) ds \right).
\end{align*}
We used the mean value theorem in the second and third inequalities.
Thus, by letting $\lambda = \theta(t,\delta y)^{-1} \ge H^{-1}(1/t)$, we see from \eqref{basicH1} and Lemma \ref{l:theta} that
\begin{align}\label{cheby1}
\P(S^1_{t/2} \ge \frac{t}{2}b(t)+\frac{3y}{8}) &\le \exp \bigg(-\frac{3y}{8\theta(t, \delta y)} + \frac{et}{\theta(t, \delta y)^2}\ \int_0^{\theta(t, \delta y)} s^2 \nu(s) ds\bigg) \nn\\
&\le \exp \bigg(-\frac{1}{8\theta(t, \delta y)} \big(3y - 16e^2 t \theta(t, \delta y)H(\theta(t, \delta y)^{-1})\big) \bigg)\nn\\
& \le \exp \bigg(-\frac{1}{8\theta(t, \delta y)} \big(3y - 16e^2 \delta y\big) \bigg) = \exp \bigg(-\frac{y}{8\theta(t, \delta y)} \bigg).
\end{align}

Consequently, from \eqref{semigroup1}, \eqref{p1supp} and \eqref{cheby1}, we deduce that
\begin{align}\label{e:estiA1}
A_1 &\le 2c_2H^{-1}(1/t)\exp \bigg(-\frac{y}{8\theta(t, \delta y)} \bigg).
\end{align}
Note that \eqref{semigroup} and \eqref{e:estiA1} hold for all $y>\delta^{-1}  H^{-1}(1/t)^{-1}$ and we have not assumed $y <R_1/2$ yet.

\noindent
\textit{Step2.}
Next, we assume  $y \in [0, R_1/2)$ and estimate $A_2$. Since $S^2$ is a compounded Poisson process, for every $z>0$ and $\rho>0$, we have that 
\begin{equation}\label{compoundP}
\P(S^2_{t/2} \in (z, z+\rho)) = \sum_{n=1}^\infty e^{-tw(\theta(t, \delta y))/2}\frac{t^n\nu_2^{n*}(z,z+\rho)}{n!}  \le \sum_{n=1}^\infty \frac{t^n\nu_2^{n*}(z,z+\rho)}{n!},
\end{equation}
where $\nu_2^{n*}$ is the $n$-fold convolution of the measure $\nu_2$.
Define
$$
f(r):= \begin{cases}
\sup_{u \ge r} \nu(u) \qquad &\text{if} \quad r<R_1/2, \\
\sup_{u \ge R_1/2} \nu(u) \qquad &\text{if} \quad r\ge R_1/2.
\end{cases}
$$
Then $f$ is a non-increasing function on $(0,\infty)$ and $\nu(r) \le f(r)$ for all $r>0$. Moreover, we see from the conditions {\bf (S-1)} and {\bf (S-3)} that $\nu(r) \asymp f(r)$ for $r \in (0,R_1/2)$. Also, we see that there exists a constant $c_5>1$ such that
$$
f(r) \le c_5f(2r) \qquad \text{for all} \;\; r>0.
$$
Indeed, if $r<R_1/4$, then $f(r) \le cf(2r)$ for some constant $c>0$ by {\bf (S-2)}. Else if $R_1/4 \le r <R_1/2$, then by {\bf (S-2)} and {\bf (S-3)}, $f(r) \le c\nu(R_1/4) \le c(\sup_{u \ge R_1} \nu(u))^{-1}f(2r) \le cf(2r)$. Otherwise, if $r\ge R_1/2$, then $f(r) = f(2r)$.

Now, we prove by induction that for every $n \ge 1$, $z>0$ and $\rho>0$, it holds that (cf. \cite[Lemma 9 and Corollary 10]{KS},)
\begin{equation}\label{comp}
\nu_2^{n*}(z, z+\rho) \le (4ec_5)^n t^{1-n} f(z)\rho.
\end{equation}

First, we see that $\nu_2(z, z+ \rho) \le f(z) \rho$. Assume that \eqref{comp} is true for $n$. Then we have
\begin{align*}
\nu_2^{(n+1)*}(z, z+\rho) &= \left(\int_{u<z/2} + \int_{u  \ge z/2} \right)\nu_2^{n*}(z-u,z-u+\rho)  \nu_2(du) \\
& \le \int_{u<z/2} (4ec_5)^n t^{1-n}f(z-u) \rho \nu_2(du) + \int_0^{z/2+ \rho} \int_{(z-v) \vee (z/2)}^{z-v+\rho} \nu_2(du) \nu_2^{n*}(dv)\\
& \le (4ec_5)^n t^{1-n} f(z/2) \rho \nu_2(\R) + f(z/2) \rho \nu_2(\R)^n \\
& \le (2e (4ec_5)^n + (2e)^n) c_5t^{-n} f(z) \rho  \le (4ec_5)^{n+1}t^{-n} f(z) \rho.  
\end{align*}
We used the Fubini's theorem in the first inequality and the fact that $\nu_2(\R) = w(\theta(t, \delta y)) \le w(w^{-1}(2e/t))=2e/t$ in the third inequality. Hence, we conclude that \eqref{comp} holds.

It follows from \eqref{compoundP} and \eqref{comp} that since $y \in [0, R_1/2)$,
\begin{align}\label{e:estiA2}
A_2 = \sup_{z>y/4} \lim_{\rho \to 0} \rho^{-1}\P(S^2_{t/2} \in (z, z+\rho)) \le  \sup_{z>y/4} tf(z) \sum_{n=1}^\infty \frac{(4ec_5)^n}{n!} \le e^{4ec_5} tf(y/4) \le c_6 t\nu(y).
\end{align}

Finally, we get the desired result from \eqref{semigroup}, \eqref{e:estiA1} and \eqref{e:estiA2}.

\vspace{2mm}

Now, we further assume that {\bf (S-3*)} holds and assume that $y > (R_1/2) \vee  (\delta^{-1}  H^{-1}(1/t)^{-1})$.  Recall that   \eqref{semigroup} and \eqref{e:estiA1} still hold for those values of $y$. Define
$
f_*(r):= \sup_{u \ge r} \nu(u).
$
Then, we see from the conditions {\bf (S)} and {\bf (S-3*)} that $\nu(r) \asymp f(r)$ for all $r \in (0,\infty)$ and there exists a constant $c_7>1$ such that $f_*(r) \le c_7f_*(2r)$ for all $r>0$. By following the above proof given in \textit{Step2.}, we get $A_2 \le e^{4ec_7}tf_*(y/4) \le c_8\nu(y)$. Thus, \eqref{e:estiA2} still hold for those values of $y$ and this completes the proof.

\vspace{2mm}

\noindent (2) We follow the proof of (1). Since $T_1 >T_0$, we see that $S^1_u$ has a transition density function $p^1(u, \cdot)$ for every $u \ge T_1$ by Proposition \ref{p:existence}. Hence, \eqref{semigroup} still hold. Also, by using Proposition \ref{hkesup} and Lemma \ref{l:largeH}(3), \eqref{p1supp} holds for all $t \in [2T_1, \infty)$.
We can prove \eqref{cheby1} by exactly the same way. Moreover,
by using Lemma \ref{l:largeH} instead of Lemma \ref{l:smallH} and the function
$
f_*(r):= \sup_{u \ge r} \nu(u)
$
instead of the function $f$, we can follow the proof in \textit{Step2.} This proves the proposition under the condition {\bf (L)}. 

Furthermore, if $T_0=0$, then for every fixed $T>0$, by Proposition \ref{lefttail:L}, \eqref{semigroup} holds for all $t \ge T$. Then, there is no difference in the proof for the last assertion. 
\qed

\smallskip

Now, we begin to prove the lower bound in Theorem \ref{righttail}. We first establish a preliminary jump type estimates for $p(t,x)$.

\begin{prop}\label{p:jumplow}

\noindent {\rm (1)}
	Suppose that the conditions {\bf (S-1)} and {\bf (S-3)} hold. Then, for every $T>0$, there exist constants $c_1,c_2>0$ such that for all $t \in (0, T]$ and $y \in [0, R_1/2)$,
	\begin{align}\label{e:lowjump}
	p(t, tb(t)+y) \ge c_1 H^{-1}(1/t) \min\left\{1, \frac{t \nu(y)}{H^{-1}(1/t)}\right\}.
	\end{align}

	Moreover, if the condition {\bf (S-3*)} further hold, then \eqref{e:lowjump} holds true for all $t \in (0,T]$ and $y \in [0, \infty)$.

\noindent {\rm (2)} Suppose that the conditions {\bf (L-1)} and {\bf (L-3)} hold. Then, there exist constants $c_1,c_2>0$ such that \eqref{e:lowjump} holds for all $t \in [2T_1, \infty)$ and $y \in [0, \infty)$ where $T_1$ is the constant in Proposition \ref{lefttail:L}.

 Moreover, if $T_0=0$ in the condition {\bf (E)}, then for every $T>0$, there exist $c_1,c_2>0$ such that \eqref{e:lowjump} holds for all $t \in [T, \infty)$ and $y \in [0, \infty)$.
\end{prop}
\proof
(1) According to Corollary \ref{c:lefttail}, it suffices to prove \eqref{e:lowjump} for $y>2H^{-1}(1/t)^{-1}$. Hence, we assume $y>2H^{-1}(1/t)^{-1}$. 

Let $\eps \in (0,1/2)$ be a small constant which will be chosen later and define
\begin{align*}
&\mu_1(s):= (1 - \eps  \1_{[ H^{-1}(1/t)^{-1}, \infty)}(s)) \nu(s), \\
&\mu_2(s):= \eps \1_{[ H^{-1}(1/t)^{-1}, \infty)}(s) \nu(s).
\end{align*}
We denote by $T^i$ the corresponding subordinator with respect to the L\'evy measure $\mu_i$ for $i=1,2$, respectively.
Since $\liminf_{s \to 0} s\mu_1(s) = \liminf_{s \to 0} s\nu(s) = 1/T_0$, by Proposition \ref{p:existence}, $T^1_u$ has a transition density function $q^1(u, \cdot )$ for every $u>T_0$. 

We claim that there exists a constant $c_3>0$ such that for all $t \in (0, T]$,
\begin{equation}\label{smallNDL}
q^1(t, tb(t) + z) \ge c_3 H^{-1}(1/t) \qquad \text{for all} \;\; z \in [0, H^{-1}(1/t)^{-1}].
\end{equation}
Indeed, we see from the conditions {\bf (S-1)} and {\bf (S-3)}, Lemma \ref{l:smallH}(1) and \eqref{basicH2} that
\begin{align*}
\sup_{s >0} |\mu_2(s)| &\le \eps c(\nu(H^{-1}(1/t)^{-1}) + 1) \le \eps c (H^{-1}(1/t)w(H^{-1}(1/t)^{-1}/2)+1)\\
& \le \eps c ( H^{-1}(1/t) H(2H^{-1}(1/t))+1) \le \eps c_4 t^{-1}H^{-1}(1/t).
\end{align*}
On the other hand, by Corollary \ref{c:lefttail}, there exists $c_5>0$ such that
$$
p(t, tb(t) + z) \ge c_5H^{-1}(1/t) \qquad \text{for all} \;\; z \in [0, H^{-1}(1/t)^{-1}].
$$
Hence, by \cite[Lemma 3.1(c)]{BGK}, we get that for all $z \in [0, H^{-1}(1/t)^{-1}]$,
\begin{align*}
q^1(t, tb(t)+z) \ge p(t, tb(t)+z) - t\sup_{s >0} |\mu_2(s)| \ge (c_5- \eps c_4)H^{-1}(1/t). 
\end{align*}
Therefore, by taking $\eps=c_5/(2c_4)$, we obtain \eqref{smallNDL}.

\smallskip

Then, since $S_t = T^1_t + T^2_t$ and $T^2$ is a compounded Poisson process, by \eqref{smallNDL} and \eqref{basicH2}, for all $t \in (0, T]$ and $y \in [0, \infty)$, 
\begin{align*}
&p(t, tb(t) + y) = \int_{\R}  q^1(t, tb(t) + y-z) \P(T^2_t \in dz) \\
&\quad\ge c_3H^{-1}(1/t) \P\big(T^2_t \in [y-H^{-1}(1/t)^{-1}, y]\big) \\
&\quad\ge c_3H^{-1}(1/t) \eps t \nu([y- H^{-1}(1/t)^{-1},y]) \exp \big(- \eps t w(H^{-1}(1/t)^{-1})\big) \\
&\quad\ge  c_3H^{-1}(1/t) \eps t H^{-1}(1/t)^{-1}  \exp \big(- 4e\eps t (H \circ H^{-1})(1/t)\big) \inf_{u \in [y- H^{-1}(1/t)^{-1},y]} \nu(u) \\
&\quad \ge c_6 t\inf_{u \in [y/2,y]} \nu(u).
\end{align*}

We see from the condition {\bf (S-1)} that $\inf_{u \in [y/2,y]} \nu(u) \asymp \nu(y)$ for all $y \in (2H^{-1}(1/t)^{-1}, R_1/2)$. Moreover, if the condition {\bf (S-3*)} further hold, then  $\inf_{u \in [y/2,y]} \nu(u) \asymp \nu(y)$ for all $y \in (2H^{-1}(1/t)^{-1}, \infty)$. Hence, we get the results.
\vspace{2mm}

\noindent(2) Fix any $N>2$ such that $NH^{-1}(1/T)^{-1} \ge R_2$. In view of Corollary \ref{c:lefttail}, we can assume that $y>NH^{-1}(1/t)^{-1} \ge R_2$. Then, by repeating the proof for (1), we get the desired result. The proof for the second assertion is exactly the same.
\qed

\vspace{2mm}

\begin{lemma}\label{boundary}
	(1) Suppose that the condition {\bf (S-1)} holds. Then, for every $a>0$ and $T>0$, there exists $c_1>0$ such that for all $t \in (0,T]$ and $y \in [H^{-1}(1/t)^{-1},R_1/2)$, 
	\begin{equation}\label{exppoly}
	\exp \big(-\frac{ay}{w^{-1}(2e/t)} \big) \le c_1\frac{t \nu(y)}{H^{-1}(1/t)}.
	\end{equation}
	
	Moreover, if the condition {\bf (S-3*)} further hold, then \eqref{exppoly} holds true for all $t \in (0,T]$ and $y \in [H^{-1}(1/t)^{-1}, \infty)$.

	\noindent
	(2) Suppose that the condition {\bf (L-1)} holds. Then, for every $a>0$ and $T>0$, there exists $c_1>0$ such that \eqref{exppoly} holds for all $t \in [T, \infty)$ and $y \in [H^{-1}(1/t)^{-1} \vee R_2, \infty)$.
\end{lemma}
\proof
Since the proofs are similar, we only give the proof for (1). By \eqref{wHinverse} and Lemma \ref{l:smallH}(1) and (2), we have that for all $y \in [H^{-1}(1/t)^{-1},R_1/2)$,
\begin{align}\label{exppoly0}
\exp \big(-\frac{ay}{w^{-1}(2e/t)} \big) &\le c_0 \left(\frac{w^{-1}(2e/t)}{y}\right)^{\alpha_1+1} \le  c\frac{w(y)}{w(w^{-1}(2e/t))}  \frac{w^{-1}(2e/t)}{y}\le c\frac{t \nu(y)}{ H^{-1}(1/t)}.
\end{align}
 In the first inequality above, we used the fact that for every $p>0$, there exists a constant $c(p)>0$ such that $e^x \ge c(p)x^p$ for all $x>0$.

Next, we further assume that the condition {\bf (S-3*)} holds. Since both conditions {\bf (S-1)} and {\bf (S-3*)} holds, there exist constants $c_2,c_3 \in (0,1)$ such that
\begin{align}\label{e:globaldouble}
c_2\sup_{u \ge r} \nu(u) \le \nu(r) \quad \text{and} \quad c_3\nu(r) \le \nu(2r) \quad \text{for all} \;\; r >0.
\end{align}
Let $r_0= (\log(c_2^{-1}c_3^{-1})+1)(a^{-1}+1)H^{-1}(1/T)^{-1}>0$. Then, $r_0>H^{-1}(1/T)^{-1} \ge H^{-1}(1/t)^{-1} $ for all $t \in (0,T]$ and by \eqref{basicH2},
\begin{align}\label{e:chooser0}
\exp \big(-ar_0/w^{-1}(2e/t) \big) \le \exp \big(-ar_0H^{-1}(1/T) \big)\le c_2c_3.
\end{align}
We first note that, using the condition {\bf (S-3*)},  we can see that the condition {\bf (S-1)} holds with $R_1 = 9r_0$  (after changing the constant $c_1$ therein). Therefore, we see that \eqref{exppoly0} holds for all $t \in (0,T]$ and $y \in [H^{-1}(1/t)^{-1}, 4r_0]$.

Now, assume that $y \in (4r_0, \infty)$. Choose $n \in \N$ such that $2^{n-1}r_0 < y \le 2^{n}r_0$. Then,  by \eqref{e:globaldouble} and \eqref{e:chooser0}, and using 
\eqref{exppoly0} for $y=r_0$, it holds that for all $t \in (0, T]$,
\begin{align*}
\frac{t \nu(y)}{ H^{-1}(1/t)} &\ge c_2\frac{t \nu(2^n r_0)}{ H^{-1}(1/t)} \ge c_2c_3^n\frac{t \nu( r_0)}{ H^{-1}(1/t)} \ge c_0^{-1}c_2c_3^n \exp \big(-\frac{ar_0}{w^{-1}(2e/t)} \big) \\
& \ge c_0^{-1}  \exp \big(-\frac{(n+1)ar_0}{w^{-1}(2e/t)} \big)\ge c_0^{-1}   \exp \big(-\frac{2^{n-1}ar_0}{w^{-1}(2e/t)} \big)\ge c_0^{-1}   \exp \big(-\frac{ay}{w^{-1}(2e/t)} \big).
\end{align*}
The fifth inequality above holds since $n \ge 3$. This completes the proof.
\qed

\begin{prop}\label{p:mixedlow}
	\noindent
	{\rm (1)} Suppose that the conditions {\bf (S-1)} and {\bf (S-3)} hold. Then, for every $T>0$, there exist constants $c_1, c_2>0$ such that for all $t \in (0, T]$ and $y \in [0, R_1/2)$,
	\begin{align}\label{e:righttai}
	p(t, tb(t)+y)  \ge c_1 H^{-1}(1/t) \min\left\{1, \frac{t \nu(y)}{H^{-1}(1/t)} + \exp \big(-\frac{c_2y}{\theta(t, y/(8e^2))} \big) \right\}
	\end{align}
	
	Moreover, if the condition {\bf (S-3*)} further hold, then \eqref{e:righttai} holds true for all $t \in (0,T]$ and $y \in [0, \infty)$.

	\noindent
	{\rm (2)} Suppose that the conditions {\bf (L-1)} and {\bf (L-3)} hold. Then, there exist constants $c_1,c_2>0$ such that \eqref{e:righttai} holds for all $t \in [2T_1, \infty)$ and $y \in [0, \infty)$ where $T_1$ is the constant in Proposition \ref{lefttail:L}.
	
	 Moreover, if $T_0=0$ in the condition {\bf (E)}, then for every $T>0$, there exist $c_1,c_2>0$ such that \eqref{e:righttai} holds for all $t \in [T, \infty)$ and $y \in [0, \infty)$.
\end{prop}
\proof
We first give the proof for (2).
Suppose that the conditions {\bf (L-1)} and {\bf (L-3)} hold. Since the proof for the case when $T_0=0$ is easier, we only give the proof for the case when $T_0>0$.

Let $\rho = (16e^2 T_1 H(w^{-1}(e/T_1)^{-1}))^{-1} \wedge (4e^2)^{-1}$. Then, by the monotonicities of $H$ and $w$,
\begin{equation}\label{enough}
\frac{1}{ 8e^2 \rho H(w^{-1}(2e/t)^{-1})} \ge 2T_1 \qquad \text{for all} \;\; t \ge 2T_1.
\end{equation}

By Corollary \ref{c:lefttail}, Lemma \ref{boundary} and Proposition \ref{p:jumplow}, it remains to prove that there are constants $c_1, c_2>0$ such that for all $t \in [2T_1, \infty)$ and $y \in [2\rho^{-1} H^{-1}(1/t)^{-1}, 8e^2D(t))$,
$$
p(t, tb(t)+y) \ge c_1 H^{-1}(1/t) \exp \big(-\frac{c_2y}{\theta(t, y/(8e^2))} \big).
$$

Fix  $t \in [2T_1, \infty)$, $y \in [2\rho^{-1}H^{-1}(1/t)^{-1}, 8e^2D(t))$ and we simply denote $\theta:=\theta(t,y/(8e^2))$. Then, since $2\rho^{-1} \ge 8e^2$, by Lemma \ref{l:theta}, we have
\begin{equation}\label{deftheta}
8e^2t\theta H(\theta^{-1})=y.
\end{equation}

Let $n=\lfloor \rho y/\theta \rfloor:=\max\{m \in \mathbb{Z}: m \le \rho y/\theta\}$. Then, since $\theta \le H^{-1}(1/t)^{-1}$, we have $n \ge \rho y H^{-1}(1/t)-1 \ge 1$. 

We claim that there exist constants $\kappa_1 \in (0,1)$ and $\kappa_2 \in (1, \infty)$ independent of $t$ and $y$ such that
\begin{equation}\label{chaining}
\kappa_1H^{-1}(n/t)^{-1} \le y/n \le \kappa_2 H^{-1}(n/t)^{-1}.
\end{equation}
Indeed, first note that  \eqref{chaining} is equivalent to $H(\kappa_1 n/y) \le n/t \le H(\kappa_2 n/y).$

Since $\rho/\theta \le \rho w^{-1}(e/T_1)^{-1}$, by Lemma \ref{l:largeH}(3) and \eqref{deftheta}, there exists a constant $c_3 \in (0,1)$ independent of $t$ and $y$ such that for every $\kappa \in [1, y/n]$,
\begin{align*}
H(\kappa n/y) \ge c_3 \kappa^{\alpha_3'} H(\rho \theta^{-1})  \ge c_3 \rho^2 \kappa^{\alpha_3'} H(\theta^{-1}) = c_3 \rho^2\kappa^{\alpha_3'} y/(8e^2t \theta) \ge c_3\rho \kappa^{\alpha_3'} n/(8e^2t)
\end{align*}
where $\alpha_3'=\alpha_3 \wedge (3/2)$.
Hence, if $y/n \ge  (8e^2 c_3^{-1}\rho^{-1})^{1/\alpha_3'}$, then by choosing $\kappa_2$ bigger than $(8e^2 c_3^{-1}\rho^{-1})^{1/\alpha_3'}$, we get the upper bound in \eqref{chaining}. Otherwise, if $y/n <  (8e^2 c_3^{-1}\rho^{-1})^{1/\alpha_3'}$, then we have
$$(8e^2 c_3^{-1}\rho^{-1})^{1/\alpha_3'} \ge y/n \ge \rho^{-1} \theta \ge \rho^{-1}w^{-1}(2e/t)$$
which implies that $t \asymp 1$. It follows that $y \asymp \theta \asymp n \asymp 1$. Therefore, by choosing $\kappa_2>(8e^2 c_3^{-1}\rho^{-1})^{1/\alpha_3'}$ sufficiently large, we obtain the upper bound in \eqref{chaining}.

On the other hand, we also have that by Lemma \ref{l:largeH}(3) and \eqref{deftheta},
\begin{align*}
H(\kappa_1 n/y) \le c_4 \kappa_1^{\alpha_3'}\rho^{\alpha_3'} H(\theta^{-1}) =  c_4 \kappa_1^{\alpha_3'}\rho^{\alpha_3'} y/(8e^2t \theta).
\end{align*}
Therefore, by choosing $\kappa_1 = (4e^2c_4^{-1}\rho^{1-\alpha_3'})^{1/\alpha_3'} \wedge (1/2)$, we get $H(\kappa_1 n/y) \le \rho y/(2t \theta) \le n/t$,
which proves the lower bound in \eqref{chaining}.

\smallskip

Define $z=y+tb(t)-tb(t/n)$ and $z_j=jz/n$ for $j=1,...,n-1$. Then, according to Lemma \ref{l:b}, \eqref{chaining} and Lemma \ref{l:largeH}(2),
\begin{align*}
y \le z &\le y + 2enH^{-1}(n/t)^{-1} + e^{-1}tH^{-1}(1/t)^{-1}w(H^{-1}(n/t)^{-1})\\
& \le (2e \kappa_1^{-1}+1)y + e^{-1}tH^{-1}(1/t)^{-1} w(\theta/(\rho \kappa_2)) \\
& \le (2e \kappa_1^{-1}+1)y + c_5  \kappa^{\alpha_4} tH^{-1}(1/t)^{-1}w(\theta)  \\
& \le (2e \kappa^{-1}+1)y + 2ec_5\kappa^{\alpha_4} H^{-1}(1/t)^{-1} \le c_6y.
\end{align*}
We used the definition that $\theta \ge w^{-1}(2e/t)$ in the fourth inequality and the assumption that $y \ge 8e^2H^{-1}(1/t)^{-1}$ in the last inequality.
Then, by \eqref{chaining}, for any $u \in (z_j-z/(2n),z_j+z/(2n))$ and $v \in (z_{j+1}- z/(2n), z_{j+1}+z/(2n))$ for some $j=1,..., n-2$, we have 
$$|u-v| \le 2z/n \le 2c_6y/n \le 2c_6 \kappa_2 H^{-1}(n/t)^{-1}.$$ Moreover, we see from \eqref{enough} and \eqref{deftheta} that 
$$
\frac{t}{n} \ge \frac{t \theta}{\rho y} = \frac{1}{8e^2 \rho H(\theta^{-1}) } \ge  \frac{1}{8e^2 \rho H(w^{-1}(2e/t)^{-1}) } \ge 2T_1.
$$
Thus, by Corollary \ref{c:lefttail}, there exists a constant $c_7 \in (0, 1)$ independent of $t$ and $y$ such that
$$
p\left(t/n, (t/n)b(t/n)+v-u\right) \ge c_7 H^{-1}(n/t),
$$
 for any $u \in (z_j-z/(2n),z_j+z/(2n))$ and $v \in (z_{j+1}- z/(2n), z_{j+1}+z/(2n))$ for some $j=1,..., n-2$.
Then, by the semigroup property and \eqref{chaining}, we get
\begin{align*}
&p(t,tb(t)+y)\\
& \ge \int_{z_{n-1}-z/(2n)}^{z_{n-1}+z/(2n)}...\int_{z_1 -z/(2n)}^{z_1 + z/(2n)} p(t/n, (t/n)b(t/n)+u_1)p(t/n, (t/n)b(t/n)+u_2-u_1)...\\
& \qquad \quad \times p(t/n, (t/n)b(t/n)+u_{n-1}-u_{n-2}) p(t/n, (t/n)b(t/n)+z-u_{n-1}) du_1...du_{n-1} \\
& \ge (c_7H^{-1}(n/t))^{n} (z/n)^{n-1} \ge c_7^n H^{-1}(n/t) \ge H^{-1}(1/t)\exp \big(-n \log c_7^{-1}\big).
\end{align*}
Since $n \asymp y/\theta$,  we obtain the desired result.
\smallskip

Now, we assume that the conditions {\bf (S-1)} and {\bf (S-3)} hold and follow the above proof. In this case, we simply let $\rho = (4e^2)^{-1}$. Since $\theta^{-1} \ge H^{-1}(1/t)^{-1} \ge H^{-1}(1/T)^{-1}$ in this case, by using Lemma \ref{l:smallH} instead of Lemma \ref{l:largeH}, we obtain \eqref{chaining}. Then, we get the result by exactly the same proof as the one given in the above.
We note that there is no difference in the proof for the second assertion in (1).
\qed

\subsection{Proofs of Theorems \ref{t:main1} and \ref{t:main2} and Corollaries \ref{c:main1},\ref{cor1}, \ref{cor2} and \ref{cor3}}

\

\begin{lemma}\label{l:simpleleft}
	(1) Suppose that the condition {\bf (S-1)} holds. Then, for every fixed $T>0$, 
	there exist constants $c_1,c_2>0$ such that for all $t \in (0,T]$ and $x \in (0,tb(t)]$,
	\begin{align}\label{e:simpleleft}
	c_1H^{-1}(1/t) \exp \big(-2tH(\sigma)\big) \le \frac{\exp \big(-t H(\sigma)\big)}{\sqrt{t(-\phi''(\sigma))}} \le c_2H^{-1}(1/t) \exp \big(-\frac{t}{2}H(\sigma)\big).
	\end{align}
	In particular, if the condition {\bf (S-1)} holds with $R_1= \infty$, then \eqref{e:simpleleft} holds for all $t \in (0, \infty)$ and $x \in (0, tb(t)]$.
	
	\noindent (2) Suppose that the condition {\bf (L-1)} holds. Then, for every fixed $T >0$, 
	there exist constants $c_1,c_2>0$ such that \eqref{e:simpleleft} holds for all $t \in [T,\infty)$ and $x \in [tb(T),tb(t)]$,
\end{lemma}
\proof (1) Observe that for all
$t \in (0,T]$ and $x\in (0, tb(t)]$, we have $\sigma \ge H^{-1}(1/t) \ge H^{-1}(1/T)$ and $tH(\sigma) \ge 1$. Hence, by Lemma \ref{l:smallH}(3) and (4), we get 
\begin{align}\label{e:s1}
c_3\sigma \exp(-2tH(\sigma)) \le \frac{\exp(-tH(\sigma))}{\sqrt{t(-\phi''(\sigma))}} &\asymp \frac{ \exp(-tH(\sigma))}{\sigma^{-1}\sqrt{tH(\sigma)}} \le \sigma \exp(-tH(\sigma)).
\end{align}
Moreover, by applying Lemma \ref{l:smallH}(4) again, 
\begin{align}\label{e:s2}
&\sigma \exp(-tH(\sigma)) = H^{-1}(1/t)\frac{\sigma}{H^{-1}(1/t)} \exp(-tH(\sigma)) \nn\\
&\quad \le c_4H^{-1}(1/t)\left(\frac{H(\sigma)}{1/t}\right)^{1/\alpha_1} \exp(-tH(\sigma)) \le  c_5H^{-1}(1/t) \exp \big(-\frac{t}{2}H(\sigma)\big)
\end{align}
and
\begin{align}\label{e:s3}
&\sigma \exp(-2tH(\sigma))= H^{-1}(1/t)\frac{\sigma}{H^{-1}(1/t)} \exp(-2tH(\sigma))  \ge H^{-1}(1/t)\exp(-2tH(\sigma)).
\end{align}
This proves the first assertion. If we further assume that $R_1=\infty$, then by combining Lemmas \ref{l:smallH} and \ref{l:largeH}, we can see that \eqref{e:s1}, \eqref{e:s2} and \eqref{e:s3} holds for all $t \in (0, \infty)$ and $x \in (0, tb(t)]$ since $tH(\sigma) \ge 1$ on those values of $t$ and $x$. 

(2) Note that for all $t \in [T, \infty)$ and $x \in [tb(T), tb(t)]$, we have $\sigma \le H^{-1}(1/T)$ and $tH(\sigma) \ge 1$. Hence, by using Lemma \ref{l:largeH} instead of Lemma \ref{l:smallH}, we can follow the proof for (1) and conclude that (2) also holds.\qed

\textit{Proof of Theorems \ref{t:main1} and \ref{t:main2}.} The results follow from Corollary \ref{c:lefttail}, Theorem \ref{righttail} and Lemmas \ref{boundary} and \ref{l:simpleleft}. \qed

\textit{Proof of Corollary \ref{c:main1}.} The results follow from Theorems \ref{t:main1} and \ref{t:main2}. \qed 

\textit{Proof of Corollaries \ref{cor1} and \ref{cor2}.} Since the proofs are similar, we only give the proof for Corollary \ref{cor1}. Under the condition {\bf (S.Pure)}, by \cite[Lemma 2.1(iii)]{CK}, we have that $w^{-1}(2e/t) \asymp H^{-1}(1/t)^{-1}$ for all $t \in (0,T]$ and hence $D(t) \asymp H^{-1}(1/t)^{-1}$ for all $t \in (0,T]$. Then, by Theorem \ref{t:main1}(2), \eqref{pure} and Corollary \ref{c:lefttail}, we get that \eqref{e:pj} holds. 

On the other hand, note that by the condition {\bf (E)}, $\nu\big((x-tb(t))_+\big) = \nu(0) = \infty$ for all $ x\le tb(t)$. Thus, by joining \eqref{e:pj} and \eqref{otherform} together, we also deduce \eqref{e:pj'}.   \qed

\textit{Proof of Corollary \ref{cor3}.}
Since the proofs for the case $T_0=0$ and the case $T_0>0$ are similar, we give the proof for the case $T_0>0$ only.
Let $T_1>0$ is the constant in Theorem \ref{t:main2}(1).
Let $\alpha_3'=\alpha_3 \wedge (3/2)$. Since $\alpha_3>1$, we also have that $\alpha_3'>1$.
Observe that for every $t \ge T_1$, by Lemma \ref{l:largeH}(4) and {\bf (L.Mixed)},
\begin{align}\label{push}
0 \le t\phi'(0)- t b(t) &= t \int_0^{H^{-1}(1/t)} (-\phi''(\lambda)) d\lambda \le c  \int_0^{H^{-1}(1/t)} \lambda^{-2}\frac{H(\lambda)}{H(H^{-1}(1/t))} d \lambda \nn\\
& \le c H^{-1}(1/t)^{-\alpha_3'}\int_0^{H^{-1}(1/t)} \lambda^{-2+\alpha_3'} d \lambda \le c_1 H^{-1}(1/t)^{-1},
\end{align}
for some constant $c_1>1$.
Moreover, by Lemma \ref{l:largeH}(3), there exists $c_2 \in (0,1)$ such that
\begin{equation}\label{e:scaling}
H(\kappa \lambda) \ge c_2\kappa^{\alpha_3'} H(\lambda) \qquad \text{for all} \;\; \kappa \ge 1, \;\; 0<\lambda \le \kappa^{-1}.
\end{equation}

 Let $y_t = y + t\phi'(0) - tb(t)$. Note that $y \le y_t \le  y+ c_1 H^{-1}(1/t)^{-1}$. By Theorem \ref{t:main2}(1), we have that for every $t \ge T_1$,
\begin{align}\label{absorb}
p(t, t\phi'(0)+y) = p(t, tb(t) + y_t) \simeq H^{-1}(1/t) \min\left\{1, \frac{t \nu(y_t)}{H^{-1}(1/t)} + \exp \big(-\frac{cy_t}{\theta(t, y_t/(8e^2))} \big) \right\}.
\end{align}

Define
$$F(t,y)=\min\left\{1, \frac{t \nu(y)}{H^{-1}(1/t)} + \exp \big(-\frac{y}{\mathscr{H}^{-1}(t/y)} \big) \right\}.$$
Then, in view of \eqref{absorb}, it remains to prove that for all $ t \ge T_1$ and $ y \ge 0,$
\begin{align}\label{FvsG}
F(t,y) \simeq \min\left\{1, \frac{t \nu(y_t)}{H^{-1}(1/t)} + \exp \big(-\frac{cy_t}{\theta(t, y_t/(8e^2))} \big) \right\}=:G(t,y,c).
\end{align}
We prove \eqref{FvsG} by considering several cases. We use the following notations below.
$$
\eps_1:=(c_2/(8e^2))^{1/(\alpha_3'-1)} \in (0,1), \quad   \kappa_1:=c_2^{-1/(\alpha_3'-1)}>1, \quad  \theta:=\theta(t, y_t/(8e^2)), 
$$

\smallskip
	\noindent
(i) Suppose that $0 \le y_t < 8e^2H^{-1}(1/t)^{-1}$. Then, we have $\theta= H^{-1}(1/t)^{-1} \ge y_t/(8e^2)$ and hence $G(t,y,1) \asymp 1$.

 We claim that it also holds that $F(t,y) \asymp 1$ which yields the desired result in this case. To prove this claim, we consider the following two cases:
 
 \smallskip
 
 (a) If $t \ge 1/H(\eps_1)$, then we see from \eqref{e:scaling} that
 \begin{align*}
 \mathscr{H}(\eps_1 H^{-1}(1/t)^{-1}) \le \frac{H^{-1}(1/t)}{\eps_1  H(\eps_1^{-1}H^{-1}(1/t))} \le \frac{tH^{-1}(1/t)}{c_2 \eps_1^{1-\alpha_3'} } =\frac{tH^{-1}(1/t)}{8e^2} \le \frac{t}{y_t} \le \frac{t}{y}.
 \end{align*}
 Thus, $\mathscr{H}^{-1}(t/y) \ge \eps_1 H^{-1}(1/t)^{-1} \ge \eps_1 y/(8e^2)$ and hence $F(t,y) \asymp 1$.
 
 (b) If $T_1 \le t \le 1/H(\eps_1)$, then $y \le y_t<8e^2/\eps_1$ and hence from the monotonicity, we get $\mathscr{H}^{-1}(t/y) \ge \mathscr{H}^{-1}(\eps_1 T_1/(8e^2)) \ge  \eps_1\mathscr{H}^{-1}(\eps_1 T_1/(8e^2))y/(8e^2)$ and hence $F(t,y) \asymp 1$.

\vspace{3mm}
	\noindent
(ii) Suppose that $8e^2H^{-1}(1/t)^{-1} \le y_t < 8e^2D(t)$. Then, by Lemma \ref{l:theta}, we have $y_t = 8e^2 t \theta H(\theta^{-1})$. Denote by $\eps_2=\eps_2(t,y) = \theta H^{-1}(1/t) \in (0,1]$ so that $\theta = \eps_2 H^{-1}(1/t)^{-1}$.
 
 \smallskip
 
 (a) Assume that $y <c_1H^{-1}(1/t)^{-1}$. Then we see from \eqref{push} and \eqref{e:scaling} that
 \begin{align*}
 2c_1H^{-1}(1/t)^{-1} > y_t = 8e^2t \theta H(\theta^{-1}) \ge 8e^2 c_2 \eps_2^{1-\alpha_3'}H^{-1}(1/t)^{-1},
 \end{align*} 
 provided that $\theta \ge 1$. Hence, if $\theta \ge 1$, then $\eps_2 \ge (4e^2c_2 /c_1)^{1/(\alpha_3'-1)}>0$ and hence $y_t \asymp \theta \asymp H^{-1}(1/t)^{-1}$. Therefore, combining with the results in (i)(a) and (i)(b), we can deduce that $F(t,y) \asymp G(t,y,1) \asymp 1$ in this case. Otherwise, if $\theta<1$, then $w^{-1}(2e/t) \le \theta <1$ and hence $t <2ew(1)$. By the similar proof to the one given in (i)(b), we can also deduce that $F(t,y) \asymp G(t,y,1) \asymp 1$.
 
 \smallskip
 
 (b) Assume that $y \ge c_1H^{-1}(1/t)^{-1}$. By the proof given in (i)(b), we may assume that $H^{-1}(1/t)^{-1} \ge R_3$ and $w^{-1}(2e/t) \ge \eps_1^{-1}$. Since \eqref{push} implies that $y \le y_t \le 2y$ in this case, by the condition {\bf (L)}, we get $\nu(y) \asymp \nu(y_t)$. Hence, it remains to prove that $\mathscr{H}^{-1}(t/y) \asymp \theta$ in this case. First, we see that by \eqref{e:scaling},
 $$
 \mathscr{H}( \eps_1 \theta) \le \frac{1}{\eps_1   \theta H(\eps_1^{-1}\theta^{-1})} \le \frac{1}{c_2 \eps_1^{1-\alpha_3'} \theta H(\theta^{-1})} = \frac{1}{8e^2 \theta H(\theta^{-1})} = \frac{t}{y_t} \le \frac{t}{y}
 $$
 and hence $\mathscr{H}^{-1}(t/y) \ge \eps_1\theta$. On the other hand, since $\kappa_1=c_2^{-1/(\alpha_3'-1)}$, by \eqref{e:scaling},
 $$
 \mathscr{H}(\kappa_1 \theta) = \inf_{ \kappa \ge \kappa_1} \frac{1}{\kappa \theta H(\kappa^{-1} \theta^{-1})} \ge \frac{c_2\kappa_1^{\alpha_3'-1}}{  \theta H(\theta^{-1})} = \frac{8e^2 t}{y_t} > \frac{t}{y}
 $$
 and hence $\mathscr{H}^{-1}(t/y) \le \kappa_1 \theta$. Therefore, we obtain $\mathscr{H}^{-1}(t/y) \asymp \theta$.

 \vspace{3mm}
 	\noindent
 (iii) Suppose that $y_t>8e^2D(t)$. If $y<c_1H^{-1}(1/t)^{-1}$, then by the proof given in (ii)(a), we get the result. Hence, suppose that $y \ge c_1H^{-1}(1/t)^{-1}$ and hence $y_t \le 2y$. By the proof given in (ii)(b), we may assume that $H^{-1}(1/t)^{-1} \ge R_2$ and $\nu(y) \asymp \nu(y_t)$.
 Moreover, by Lemma \ref{l:largeH}(1) and \eqref{basicH2}, we see that $tH^{-1}(1/t)^{-1}\nu(y_t) \le c t y_t \nu(y_t) \le c t w(y_t) \le ct H(y_t^{-1}) \le c$. Moreover, by Lemma \ref{boundary} and the condition {\bf (L)}, for any fixed $a>0$,
 $$
 \exp \big(-\frac{ay_t}{\theta(t, y_t/(8e^2))}\big) \le \exp \big(-ac_1\frac{c_1^{-1}y}{w^{-1}(2e/t)}\big)  \le c\frac{t\nu(c_1^{-1}y)}{H^{-1}(1/t)} c\le \frac{t\nu(y)}{H^{-1}(1/t)}.
 $$ 
 Thus, $G(t,y,1) \asymp tH^{-1}(1/t)^{-1}\nu(y)$ in this case. Therefore, it remains to prove that there exists $c_3>0$ such that
 \begin{equation}\label{last} 
 \exp \big(-\frac{y}{\mathscr{H}^{-1}(t/y)}\big) \le c_3 \frac{t \nu(y)}{H^{-1}(1/t)}.
 \end{equation}
 
 As before, we may assume that $w^{-1}(2e/t) \ge 1$. Then, since $\theta  = w^{-1}(2e/t)$ in this case, by \eqref{e:scaling} and Lemma \ref{l:theta},
 \begin{align*}
 \mathscr{H}(\kappa_1 w^{-1}(2e/t)) = \inf_{ \kappa \ge \kappa_1} \frac{1}{\kappa \theta H(\kappa^{-1} \theta^{-1})} \ge \frac{c_2\kappa_1^{\alpha_3'-1}}{ \theta H(\theta^{-1})} =  \frac{8e^2 t}{y_t} > \frac{t}{y},
 \end{align*}
 which implies that $\mathscr{H}^{-1}(t/y) \le \kappa_1 w^{-1}(2e/t)$. Hence, we get \eqref{last} from Lemma \ref{boundary}. This completes the proof.
\qed

\section{Examples}

In this section, we provide non-trivial and concrete  examples of subordinators which our main results can be applied to. 

\smallskip

Recall that $b(t)=(\phi' \circ H^{-1})(1/t)$ and $\sigma=(\phi')^{-1}(x/t)$.

\subsection{Polynomially decaying L\'evy measure perturbed by a logarithmic function}
	In  this subsection, we use the notation $\log^p x := (\log x)^p$ for $p \in \R$ and $x>1$. 	
	Suppose that $\gamma_1 \in [0,1]$ and $p \in \R$.
	If $\gamma_1=0$ we further assume that  $p>0$ and, if $\gamma_1 =1$ then we further assume that $p<-1$.

Let $f:(1,\infty) \to (0,\infty)$ be a measurable function satisfying
\begin{align}\label{largenu}
\int_1^\infty f(r)dr<\infty, \quad c_1\sup_{u \ge r} f(u) \le f(r) \quad \text{and} \quad c_2f(r) \le f(2r) \quad \text{for all} \;\; r \ge 1.
\end{align}
for some constants $c_1,c_2>0$ and
	$$\nu(s):=\1_{\{0<s \le 1\}}s^{-1-\gamma_1}\log^p (1+1/s)+\1_{\{s > 1\}}f(s).$$
Let $S$ be a subordinator with the Laplace exponent
	$
	\phi(\lambda) = \int_0^\infty (1-e^{-\lambda t})\nu(s)ds.
	$
	 We see that $S$ satisfies the conditions {\bf (E)} with $T_0=0$. Thus, by Proposition \ref{p:existence}, for every $t>0$, the transition density $p(t,x)$ of $S_t$ exists and is a continuous bounded function. In the following, 
	 we show that, using  our main results, we can get the precise two-sided estimates on $p(t,x)$. 
	 
\subsubsection{Small time estimates.} Below, we assume that $\gamma_1>0$ and obtain estimates on $p(t,x)$ for $t \in (0,2]$. Note that the conditions  {\bf (S.Pure)} and  {\bf (S-3*)}   hold. Using Lemmas \ref{l:smallH}(1$\&$4) and \ref{l:L1diff}, and \eqref{estib}, for every fixed $\lambda_0>0$, we get that for all $\lambda \ge \lambda_0$,
\begin{align}\label{e:smallest}
H(\lambda) &\asymp  \lambda^{\gamma_1}\log^p (1+ \lambda), \quad \qquad H^{-1}(\lambda) \asymp \lambda^{1/\gamma_1}\log^{-p/\gamma_1}(1+\lambda),\nn\\
 \phi'(\lambda) &\asymp \lambda^{-1}\phi(\lambda) \asymp \1_{\{\gamma_1 \in (0,1)\}}\lambda^{\gamma_1-1} \log^p (1+ \lambda)  + \1_{\{\gamma_1 =1 \}}\log^{p+1} (1+ \lambda),\nn\\
 (\phi')^{-1}(1/\lambda) & \simeq \1_{\{\gamma_1 \in (0,1)\}}\lambda^{1/(1-\gamma_1)} \log^{p/(1-\gamma_1)} (1+ \lambda)  + \1_{\{\gamma_1 =1 \}} \exp \big(c\lambda^{-1/(p+1)}\big),\nn\\
 \lambda^{-1}b(1/\lambda) &\asymp  \1_{\{\gamma_1 \in (0,1)\}}\lambda^{-1/\gamma_1} \log^{p/\gamma_1} (1+ \lambda) + \1_{\{\gamma_1 =1 \}}\lambda^{-1}\log^{p+1} (1+ \lambda).
\end{align}
In particular, $tb(t) \asymp H^{-1}(1/t)^{-1}$ for $t \in (0,2]$ unless $\gamma_1=1$.

\smallskip
 	\noindent
	 (i) Suppose that $\gamma_1 \in (0,1)$. Then, for all $t \in (0,2]$ and $x \in (0, tb(t)]$, by \eqref{e:smallest},
	 \begin{align*}
	 tb(t)&\asymp H^{-1}(1/t)^{-1} \asymp t^{1/\gamma_1} \log^{p/\gamma_1} (1+1/t), \\
	 \sigma&\asymp (t/x)^{1/(1-\gamma_1)} \log^{p /(1-\gamma_1)}(1+t/x),\\
	 H(\sigma) &\asymp \sigma^{\gamma_1}\log^p(1+\sigma) \asymp \left(\displaystyle t/x\right)^{\gamma_1/(1-\gamma_1)}\displaystyle\log^{p/(1-\gamma_1)}(1+ t/x).
	 \end{align*} 
	Hence, by Theorem \ref{t:main1} and Corollary \ref{cor1}, we have that for all $t \in (0,2]$, 
	\begin{align}\label{e:logsmall}
	p(t,x) \begin{cases}
	\simeq t^{-1/\gamma_1} \log^{-p/\gamma_1} (1+1/t)  \exp \left(-c t\left(\displaystyle\frac{t}{x}\right)^{\frac{\gamma_1}{1-\gamma_1}}\displaystyle\log^{\frac{p}{1-\gamma_1}}\big(1+ \frac{t}{x}\big)  \right), \\
		\vspace{1mm}
	\qquad \qquad \qquad \qquad \qquad \qquad \qquad \quad  \qquad \quad \text{if} \;\; x \in (0,  2t^{1/\gamma_1} \log^{p/\gamma_1} (1+1/t)],\\
		\vspace{1mm}
	\asymp t x^{-1-\gamma_1} \log^p (1+1/x), \;\qquad \quad  \qquad \qquad\text{if} \;\; x \in (2t^{1/\gamma_1} \log^{p/\gamma_1} (1+1/t), 1],\\
	\asymp t f(x), \qquad\quad\qquad\qquad\qquad \qquad  \quad \qquad\text{if} \;\; x \in (1, \infty).
	\end{cases}
	\end{align}
	In the first and second comparison in \eqref{e:logsmall}, we used the following observation: In this case, for every fixed $a>0$, by \eqref{e:smallest}, it holds that for all $t \in (0,2]$, $tb(t) \asymp H^{-1}(1/t)^{-1}$ and 
	\begin{align*}
	t H \circ\sigma |_{x= atb(t)}= t(H \circ (\phi')^{-1})(ab(t)) \asymp t(H \circ (\phi')^{-1})(b(t))=1.
	\end{align*}
	Hence, according to Corollary \ref{c:lefttail} and Lemma \ref{l:simpleleft}, for every fixed $a>0$, we get $p(t, atb(t)) \asymp t^{-1/\gamma_1}\log^{-p/\gamma_1}(1+1/t) \asymp H^{-1}(1/t)$ for $t\in (0,2]$. Therefore, we can use $(0,  2t^{1/\gamma_1} \log^{p/\gamma_1} (1+1/t)]$ instead of $(0, tb(t)+H^{-1}(1/t)^{-1}]$, and use $(2t^{1/\gamma_1} \log^{p/\gamma_1} (1+1/t), 1]$ instead of $(tb(t)+H^{-1}(1/t)^{-1},1]$ in \eqref{e:logsmall}.
	
	\vspace{2mm}
		\noindent
		 (ii) Suppose that $\gamma_1 = 1$ and $p<-1$. Then,  for all $t \in (0,2]$ and $x \in (0, tb(t)]$, by \eqref{e:smallest},
		 \begin{align}\label{asymp1}
		 tb(t)&\asymp t \log^{p+1}(1+1/t), \quad H^{-1}(1/t)^{-1} \asymp t \log^{p}(1+1/t), \quad \sigma\simeq \exp\left(c\left(t/x\right)^{-1/(p+1)}\right).
		 \end{align} 
		Moreover, for all $t \in (0,2]$ and $x \in (0,tb(t)]$, we get
		$$
		\exp(-ct H(\sigma)) \simeq \exp(-ct \sigma^c) \simeq \exp\left(-ct \exp \big(c \left(t/x\right)^{-1/(p+1)}\big)\right).
		$$
		  Thus, by 
	 Theorem \ref{t:main1}, and  Corollaries \ref{cor1} and \ref{c:lefttail}, we have that for all $t \in (0,2]$, 
	\begin{align*}
	p(t,x) \begin{cases}
		\vspace{1mm}
	\simeq  t^{-1} \log^{-p} (1+1/t)\exp \left(\displaystyle - ct \exp \left(c\left(\frac{t}{x}\right)^{\frac{-1}{p+1}} \right) \right), &\text{if} \; x \in (0,  tb(t) + t \log^{p} (1+1/t)],\\
	\asymp t (x-tb(t))^{-2} \log^p \big(1+1/(x-tb(t))\big),   &\text{if} \; x \in (tb(t) + t \log^{p} (1+1/t), 1],\\
	\asymp t f(x),&\text{if} \; x \in (1, \infty).
	\end{cases}
	\end{align*}
	In particular, for $t \in (0,2]$ and $x \in (tb(t),  tb(t) + t \log^{p} (1+1/t)]$, by \eqref{asymp1} and Corollary \ref{c:lefttail}, 
	$$1 \simeq \exp \left( - ct \exp \left(c\left(t/x\right)^{-1/(p+1)} \right) \right) \quad \text{and} \quad p(t,x) \asymp t^{-1} \log^{-p} (1+1/t) \asymp H^{-1}(1/t)^{-1}.$$
	We note that for $t \in (0,2]$,
	$$
	p(t,tb(t)) \asymp t^{-1}\log^{-p}(1+1/t) \not\asymp t^{-1}\log^{-p-2}(1+1/t) \asymp p(t,2tb(t)).
	$$

	\subsubsection{Large time estimates.} Next, we further assume that 
	$
	f(s)=s^{-1-\gamma_2} \log^q(1+s)
	$ so that
	$$\nu(s)=\1_{\{0<s \le 1\}}s^{-1-\gamma_1}\log^p (1+1/s)+\1_{\{s > 1\}}s^{-1-\gamma_2} \log^q(1+s),$$
	for some $\gamma_2 \in (1, \infty)$ and $q \in \R$,
	and obtain estimates on $p(t,x)$ for $t\in[2, \infty)$. Clearly, \eqref{largenu} is satisfied. Since the condition {\bf (L.Mixed)} holds, by Remark \ref{r:remark2}(2) and Lemma \ref{l:largeH}(4), for every fixed $\lambda_0>0$, we get that for all $\lambda \in (0, \lambda_0]$,
	\begin{align}\label{e:HHin}
	\qquad  \phi(\lambda) &\asymp \lambda, \qquad \phi'(\lambda) \asymp \phi'(0) \asymp 1,\nn\\
	H(\lambda) &\asymp  \begin{cases}
	\lambda^{\gamma_2} \log^q (1+ 1/\lambda), \qquad &\text{ if} \;\; \gamma_2<2, \\
	\lambda^2 \log^{q+1}(1+1/\lambda), \qquad &\text{ if} \;\; \gamma_2 = 2, \; q>-1, \\
	\lambda^2 \log\log(1+1/\lambda), \qquad &\text{ if} \;\; \gamma_2 = 2, \; q=-1, \\
	\lambda^2, \qquad &\text{  if} \;\; \gamma_2=2\;\; \text{and}\;\; q<-1,\;\; \text{or}\;\; \gamma_2>2,
	\end{cases}\nn\\
		H^{-1}(\lambda) &\asymp  \begin{cases}
	\lambda^{1/\gamma_2} \log^{-q/\gamma_2} (1+ 1/\lambda),  &\text{if} \;\; \gamma_2<2, \\
	\lambda^{1/2} \log^{-(q+1)/2}(1+1/\lambda),  &\text{if} \;\; \gamma_2 = 2, \; q>-1, \\
	\lambda^{1/2} \log^{-1}\log(1+1/\lambda),  &\text{if} \;\; \gamma_2 = 2, \; q=-1, \\
	\lambda^{1/2}, \qquad &\text{if} \;\; \gamma_2=2\;\; \text{and}\;\; q<-1,\;\; \text{or}\;\; \gamma_2>2.
	\end{cases}
		\end{align}
	We also have that for all $s  \ge \lambda_0$,
		\begin{align}\label{e:SH}
	\mathscr{H}^{-1}(s) &\asymp  \begin{cases}
s \log^{q+1}(1+s),  &\text{if} \;\; \gamma_2 = 2, \; q>-1, \\
	s \log\log(1+s),  &\text{if} \;\; \gamma_2 = 2, \; q=-1, \\
	s, \qquad &\text{if} \;\; \gamma_2=2\;\; \text{and}\;\; q<-1,\;\; \text{or}\;\; \gamma_2>2,
	\end{cases}
	\end{align}
	
	\vspace{2mm}

	Note that even though $t$ is large, $\sigma$ can be arbitrary big.
Hence, to obtain large time estimates on $p(t,x)$, we still need estimates for $\phi''(\lambda)$ and $H(\lambda)$ for large $\lambda$ to calculate the function in \eqref{e:main1}.
	To cover the case $\gamma_1=0$, we observe that by Lemma \ref{l:L1diff} and \eqref{e:CK21}, if $\gamma_1=0$ and $\lambda_0>0$, we get that for all $\lambda \in [\lambda_0, \infty)$,
	\begin{align*}
	\phi'(\lambda) &\asymp \lambda^{-1} \log^p(1+\lambda), \quad (\phi')^{-1}(1/\lambda) \asymp \lambda \log^p(1+\lambda),\quad |\phi''(\lambda)| \asymp \lambda^{-2} \log^p(1+\lambda), \\
	 w(1/\lambda) &\asymp \int_{1/\lambda}^{2/\lambda_0}s^{-1}\log^p(1+1/s)ds + 1\asymp \log^{p+1}(1+\lambda),\\
	H(\lambda) &\asymp \lambda^2\int_0^{1/\lambda} sw(s)ds \asymp \log^{p+1}(1+\lambda).
	\end{align*}
	Thus, if $\gamma_1=0$ and $p>0$, then for all $t \in [2,\infty)$ and $x \in (0,t(\phi' \circ H^{-1})(1)]$, since $\sigma \ge H^{-1}(1)$,
	\begin{align*}
 \sigma \asymp \frac{t}{x} \log^p (1+\frac{t}{x}), \qquad |\phi''(\sigma)| \asymp \frac{x^2}{t^2} \log^{-p}(1+\frac{t}{x}),\qquad H(\sigma)  \asymp \log^{p+1}(1+\frac{t}{x}),
	\end{align*}
	and hence
	\begin{align}\label{e:PS} \frac{\exp\big(-tH(\sigma)\big)}{\sqrt{t(-\phi''(\sigma))}} &\simeq 	t^{-1/2} \left(\displaystyle\frac{t}{x}\right) \displaystyle\log^{\frac{p}{2}}\big(1+ \frac{t}{x}\big)\exp \left(-(c t+1)\displaystyle\log^{p+1}\big(1+ \frac{t}{x}\big)  \right)\nn\\
	& \simeq t^{1/(2p+2)}x^{-1}\exp \left(-c t\displaystyle\log^{p+1}\big(1+ \frac{t}{x}\big)  \right)\exp \left(-\displaystyle\log\big(1+ \frac{t}{x}\big)  \right)\nn\\
	& \asymp t^{-(2p+1)/(2p+2)}\exp \left(-c t\displaystyle\log^{p+1}\big(1+ \frac{t}{x}\big)  \right).
	\end{align}

	Define 
	\begin{align*}
	& p_S(t,x,c):= \\
	&\begin{cases}
t^{-(2p+1)/(2p+2)}\exp \left(-c t\displaystyle\log^{p+1}\big(1+ \frac{t}{x}\big)  \right), \; &\text{if} \; \gamma_1=0,\\
	t^{-1} \log^{-p} (1+1/t)\exp \left(  \displaystyle - ct \exp \left(\left(\frac{t}{x}\right)^{\frac{-1}{p+1}} \right) \right), \; &\text{if} \; \gamma_1=1,\\
	t^{-1/\gamma_1} \log^{-p/\gamma_1} (1+1/t) \exp \left(-c t\left(\displaystyle\frac{t}{x}\right)^{\frac{\gamma_1}{1-\gamma_1}}\displaystyle\log^{\frac{p}{1-\gamma_1}}\big(1+ \frac{t}{x}\big)  \right), \; &\text{otherwise}.
	\end{cases}
	\end{align*}
 The function $p_S(t,x,c)$ with $\gamma_1>0$ appears in small time left tail estimates on $p(t,x)$ in  Section 4.1.1.
 We will see that the function $p_S(t,x,c)$ also appears in large time left tail estimates on $p(t,x)$.
	
	\smallskip
		\noindent
	(i) Suppose that $\gamma_2 < 2$. Since both conditions {\bf (L.Pure)} and {\bf (L.Mixed)} hold in this case, by \eqref{e:HHin},  \eqref{e:PS}, Theorem \ref{t:main2}, and Corollaries \ref{cor2} and \ref{cor3}, it holds that for all $t \in [2,\infty)$, 
	\begin{align*}
	p(t,x) \begin{cases}
	\simeq p_S(t,x,c),   \qquad\qquad\qquad\qquad\quad\quad	\qquad\qquad\qquad\qquad\;\; \text{if} \;\; x \in (0, tb(1)],\\
	\simeq t^{-1/\gamma_2} \log^{-q/\gamma_2}(1+t)\exp \big(-ct \sigma^{\gamma_2} \displaystyle\log^q\big(1+1/\sigma \big)\big),\quad\; \text{if} \;\; x  \in (tb(1),t\phi'(0)),\\
	\asymp t^{-1/\gamma_2} \log^{-q/\gamma_2}(1+t), \;\;	\quad\qquad\; \text{if} \;\; x=t\phi'(0)+y, \; y \in [0,t^{1/\gamma_2} \log^{q/\gamma_2}(1+t)),\\
	\asymp t y^{-1-\gamma_2} \log^q (1+y), \qquad\qquad\;\; \text{if} \;\; x=t\phi'(0)+y, \; y \in [t^{1/\gamma_2} \log^{q/\gamma_2}(1+t),\infty).
	\end{cases}
	\end{align*}
 In the second comparison, we used the following observation: by Corollary \ref{c:lefttail}, \eqref{push} and \eqref{e:HHin}, we get that for all $t \in [2, \infty)$ and $x \in [tb(t), t\phi'(0))$,
 \begin{equation}\label{diagonal}
p(t,x) \asymp H^{-1}(1/t) \asymp  t^{-1/\gamma_2} \log^{-q/\gamma_2}(1+t) \quad \text{and} \quad t \sigma^{\gamma_2} \displaystyle\log^q\big(1+1/\sigma \big) \asymp tH(\sigma) \le 1.
 \end{equation}

 Note that by Lemma \ref{l:largeH}(4) and \eqref{e:HHin}, for $t \in [2, \infty)$ and  $x \in (tb(1),tb(t))$,  
	\begin{align}\label{estisigma}
	\phi'(0)-\frac{x}{t} = \int_0^{\sigma}(-\phi''(u))du \asymp \int_0^{\sigma} u^{\gamma_2-2}\log^q(1+1/u)du \asymp \sigma^{\gamma_2-1}\log^q(1+1/\sigma)
	\end{align}
and hence
$$
\sigma \asymp \left(\phi'(0)-\frac{x}{t}\right)^{1/(\gamma_2-1)} \log^{-q/(\gamma_2-1)}\left(1+ \frac{1}{\phi'(0)-x/t}\right).
$$
	\smallskip
		\noindent
		(ii) Suppose that $\gamma_2 = 2$, $q>-1$.
	By \eqref{e:HHin},  \eqref{e:PS}, Theorem \ref{t:main2} and Corollary \ref{cor3}, it holds that for all $t \in [2,\infty)$, 
	\begin{align*}
	p(t,x) \begin{cases}
	\simeq p_S(t,x,c),  \quad \qquad\qquad\qquad\qquad\qquad\quad	\qquad\qquad\qquad\qquad\; \text{if} \;\; x \in (0, tb(1)],\\
	\vspace{1mm}
	\simeq t^{-1/2} \log^{-(q+1)/2}(1+t)\exp \big(-ct \sigma^{2} \displaystyle\log^{q+1}\big(1+1/\sigma \big)\big), \;\;\quad \text{if} \;\;  x  \in (tb(1),t\phi'(0)),\\
	\asymp t^{-1/2} \log^{-(q+1)/2}(1+t), \;\;	\qquad\quad \text{if} \;\; x=t\phi'(0)+y, \; y \in [0,t^{1/2} \log^{(q+1)/2}(1+t)),\\
	\simeq t^{-1/2} \log^{-(q+1)/2}(1+t) \exp \left(-c \displaystyle\frac{y^2}{t \log^{q+1}(1+t)} \right) +  t y^{-3} \log^q (1+y),\\
	\qquad\qquad\qquad\qquad\qquad\qquad\qquad\;\; \text{if} \;\; x=t\phi'(0)+y, \; y \in [t^{1/2} \log^{(q+1)/2}(1+t),\infty).
	\end{cases}
	\end{align*}
	We used a similar argument to \eqref{diagonal} in the second comparison.
	In the last comparison, we used the facts that the exponential term is dominated by $t\nu(y)$ for all $y \ge D(t)$ and for all $t \in [2,\infty)$, we have $w^{-1}(2e/t) \asymp t^{1/2}\log^{q/2}(1+t) \ge ct^{1/3}$,
	$$
	D(t) \asymp t\max_{s \in [w^{-1}(2e/t),H^{-1}(1/t)^{-1}]}s^{-1}\log^{q+1}(1+s) \le c t^{2/3} \log^{q+1}(1+t)  \le ct^{3/4},
	$$
	and hence by \eqref{e:SH}, for all $t \in [2, \infty)$ and $y\in [0, D(t))$, 
	(cf. \cite[Corollary 6.3]{BKKL},)
	\begin{align}\label{slowexp}
	\exp \left(-c \frac{y}{\mathscr{H}^{-1}(t/y)} \right) \simeq \exp \left(-c \frac{y^2}{t \log^{q+1}(1+t/y)} \right) \simeq \exp \left(-c \frac{y^2}{t \log^{q+1}(1+t)} \right).
	\end{align}
	
	In particular, we can see that for every fixed $\eps>0$, there are comparison constants such that $p(t,t\phi'(0)+y) \asymp ty^{-3}\log^q(1+y)$ for all $y  \ge t^{1/2} \log^{(q+1+\eps)/2}(1+t)$.
	
	We also note that by a similar calculation to \eqref{estisigma}, for $t \in [2, \infty)$ and $x \in (tb(1), tb(t))$,
$$
\sigma \asymp \left(\phi'(0)-\frac{x}{t}\right) \log^{-(q+1)}\left(1+ \frac{1}{\phi'(0)-x/t}\right).
$$

	\smallskip
	\noindent
	(iii) Suppose that $\gamma_2 = 2$, $q=-1$.
By \eqref{e:HHin}, \eqref{e:PS}, Theorem \ref{t:main2} and Corollary \ref{cor3}, it holds that for all $t \in [2,\infty)$, 
\begin{align*}
p(t,x) \begin{cases}
\simeq p_S(t,x,c),  \quad\qquad \qquad\qquad\qquad\qquad	\qquad\qquad\qquad\qquad\; \text{if} \;\; x \in (0, tb(1)],\\
\vspace{1mm}
\simeq t^{-1/2}  \log^{-1}\log(1+t)\exp \big(-ct \sigma^{2} \displaystyle\log\log\big(1+1/\sigma \big)\big),\quad\; \text{if} \;\;  x  \in (tb(1),t\phi'(0)),\\
\asymp t^{-1/2} \log^{-1}\log(1+t), \qquad	\qquad\;\; \text{if} \;\; x=t\phi'(0)+y, \; y \in [0,t^{1/2} \log\log(1+t)),\\
\simeq t^{-1/2} \log^{-1}\log(1+t) \exp \left(-c \displaystyle\frac{y^2}{t \log\log(1+t)} \right) +  t y^{-3} \log^{-1} (1+y),\\
\qquad\qquad\qquad\qquad\qquad\qquad\qquad\;\;\; \text{if} \;\; x=t\phi'(0)+y, \; y \in [t^{1/2} \log\log(1+t),\infty).
\end{cases}
\end{align*}
We used a similar argument to \eqref{diagonal} in the second comparison. Also, the last comparison holds by a similar argument to the one which is used to obtain \eqref{slowexp}.

In particular, we can see that for every fixed $\eps>0$, there are comparison constants such that $p(t,t\phi'(0)+y) \asymp ty^{-3}\log^{-1}(1+y)$ for all $y  \ge t^{1/2} \log^{\eps}(1+t)$.

	We also note that by a similar calculation to \eqref{estisigma}, for $t \in [2, \infty)$ and $x \in (tb(1), tb(t))$,
$$
\sigma \asymp \left(\phi'(0)-\frac{x}{t}\right) \log^{-1}\log\left(1+ \frac{1}{\phi'(0)-x/t}\right).
$$

\smallskip

\noindent
(iv) Suppose that either $\gamma_2 = 2$, $q<-1$ or $\gamma_2>2$.
By \eqref{e:HHin}, \eqref{e:PS}, Theorem \ref{t:main2} and Corollary \ref{cor3}, it holds that  for all $t \in [2,\infty)$, 
\begin{align*}
&p(t,x) \begin{cases}
\simeq p_S(t,x,c),  &\text{ if} \;\; x \in (0, tb(1)],\\
\vspace{1mm}
\simeq t^{-1/2} \exp \big(-ct \sigma^{2} \big), &\text{ if} \;\;  x  \in (tb(1),t\phi'(0)),\\
\asymp t^{-1/2},\;\; & \text{ if} \;\; x=t\phi'(0)+y, \; y \in [0,t^{1/2}),\\
\simeq t^{-1/2} \exp \left(-c \displaystyle\frac{y^2}{t} \right) +  t y^{-1-\gamma_2} \log^q (1+y),\;\; &\text{ if} \;\; x=t\phi'(0)+y, \; y \in [t^{1/2},\infty).
\end{cases}
\end{align*}
We used a similar argument to \eqref{diagonal} in the second comparison.

In particular, we can see that for every fixed $\eps>0$, there are comparison constants such that $p(t,t\phi'(0)+y) \asymp ty^{-1-\gamma_2}\log^q(1+y)$ for all $y  \ge t^{1/2}\log^{1/2+\eps}(1+t)$. Indeed, for all $t \in [2, \infty)$ and $y  \ge t^{1/2}\log^{1/2+\eps}(1+t)$,
\begin{align*}
&t^{-1/2} \exp\left(-c_1 \displaystyle\frac{y^2}{t} \right) \le t^{-1/2} \exp\left(-\frac{c_1}{2} \log^{1+2\eps}(1+t) \right) \exp\left(-\frac{c_1}{2} \displaystyle\frac{y^2}{t} \right) \\
&\;\;\; \le c_2 t^{-1/2-\gamma_2} \left(\frac{t}{y^2}\right)^{1+\gamma_2} =  \frac{c_2}{t^{1/2}y}\frac{t \log^q(1+y)}{y^{1+\gamma_2}} \le \frac{c_2}{2}\frac{t \log^q(1+y)}{y^{1+\gamma_2}}.
\end{align*}

	Note that $\phi''(0)<\infty$ in this case. Hence, we get that for $t \in [2, \infty)$ and $x \in (tb(1), tb(t))$, by the mean value theorem,
	$
	\phi''(H^{-1}(1/2)) \le \sigma^{-1}(\phi'(0)-x/t) \le \phi''(0)
	$
	and hence
$$
\frac{\phi'(0)-x/t}{\phi''(0)} \le \sigma \le \frac{\phi'(0)-x/t}{\phi''(H^{-1}(1/2))}.
$$

\bigskip

\subsection{An example of varying transition density estimates}
	 In this subsection, we give an example of subordinator whose   transition density has the estimates given in Theorem \ref{t:main2} and the exponential term in estimates only appears at specific time ranges.
	 
	 \smallskip
	Define an increasing sequence $(a_n)_{n \ge 0}$ as follows:
	\begin{equation}\label{defan}
	a_0:=0, \qquad a_1:=3, \qquad a_{n+1}:=\exp(a_n^{3/2}) \quad \text{for} \;\; n \ge 1.
	\end{equation}
Using this $(a_n)_{n \ge 0}$, we define a non-decreasing function $\psi:(0, \infty) \to (0, \infty)$ as follows:
	$$
	\psi(r)=\begin{cases}
	(4/3)r^{1/2} &\qquad \text{for} \;\; r \in (0,a_1],\\
	r^4+\psi(a_{2n-1})-a_{2n-1}^{4} &\qquad \text{for} \;\; r \in (a_{2n-1}, a_{2n}],\\
	(4/3)r^{1/2}+\psi(a_{2n})-(4/3)a_{2n}^{1/2} &\qquad \text{for} \;\; r \in (a_{2n}, a_{2n+1}].
	\end{cases}
	$$
One can easily check that there exist $c_2 \ge c_1>0$ such that
	\begin{equation}\label{scalepsi}
	c_1 \left(\frac{R}{r}\right)^{1/2} \le \frac{\psi(R)}{\psi(r)} \le c_2 \left(\frac{R}{r}\right)^4 \qquad \text{for all} \;\; 0<r\le R.
	\end{equation}

	Let 
	$$
	\Phi(r):=\frac{r^2}{2\int_0^r s\psi(s)^{-1}ds}.
	$$
	Then by \cite[Lemma 2.4]{BKKL} and \eqref{scalepsi}, there exists a constant $c_3>0$ such that
	\begin{equation}\label{scalephi}
	c_3 \left(\frac{R}{r}\right)^{1/2} \le \frac{\Phi(R)}{\Phi(r)} \le  \left(\frac{R}{r}\right)^2 \qquad \text{for all} \;\; 0<r\le R.
	\end{equation}

	\subsubsection{Preliminary calculations}
	
	\begin{lemma}\label{basic}
		
		For every $\eps \in (0,1)$, there exists $N \in \N$ such that for every $n \ge N$, the following estimates hold:

		\noindent {\rm (1)} For every $r \in [a_{2n+1}^{1-\eps}, a_{2n+1}]$,
		\begin{equation}\label{odd}
		\frac{4}{3}r^{1/2} \le \psi(r) \le 2r^{1/2} \qquad \text{and} \qquad r^{1/2} \le \Phi(r) \le 2r^{1/2}.
		\end{equation}
		
		\noindent {\rm (2)} For every $r \in [a_{2n}^{1-\eps}, a_{2n}]$,
		\begin{equation}\label{even}
		\frac{1}{2}r^4 \le \psi(r) \le r^4 \qquad \text{and} \qquad \frac{2(1-\eps)r^2}{3\log r} \le \Phi(r) \le \frac{2r^2}{\log r}.
		\end{equation}

	\end{lemma}
	\proof
	From the definition \eqref{defan} of the sequence $(a_n)$,
	by choosing $N$ sufficiently large, we can assume that $[a_{2n+1}^{1-\eps},a_{2n+1}]\subset (a_{2n},a_{2n+1}]$ and $[a_{2n+1}^{1-\eps},a_{2n+1}]\subset (a_{2n},a_{2n+1}]$ for all $n \ge N$.
	
	First, we prove the assertions for the function $\psi$. From the construction, we have
	\begin{equation}\label{basicpsi} 
	\frac{4}{3}r^{1/2} \le \psi(r) \le r^4 \qquad \text{for all} \;\; r \ge 1.
	\end{equation}
	Moreover, for all sufficiently large $n$ and $r \in [a_{2n+1}^{1-\eps}, a_{2n+1}]$, by \eqref{defan},
	$$
	\psi(r) \le \bigg(1+ \frac{a_{2n}^4}{(4/3)a_{2n+1}^{(1-\eps)/2}}\bigg) \frac{4}{3}r^{1/2} \le \left(1+ a_{2n}^4 \exp\big(- 2^{-1}(1-\eps)a_{2n}^{3/2}\big) \right)\frac{4}{3}r^{1/2}.
	$$ 
	Similarly, for all sufficiently large $n$ and $r \in [a_{2n}^{1-\eps}, a_{2n}]$, 
	$$
	\psi(r) \ge \bigg(1- \frac{a_{2n-1}^4}{a_{2n}^{4(1-\eps)}}\bigg)r^4 \ge \left(1- a_{2n-1}^4\exp\big(-4(1-\eps)a_{2n-1}^{3/2}\big)\right)r^4.
	$$
	Since $\lim_{x \to \infty}x^4e^{-4(1-\eps)x^{3/2}}= \lim_{x \to \infty} x^4 e^{-2^{-1}(1-\eps)x^{3/2}}=0$, we deduce the results for $\psi$.
	
	Now, we prove the assertions for the function $\Phi$. Fix $\eps' \in (0, 1- \eps)$. By using the results for $\psi$ and \eqref{basicpsi}, we can see that for all sufficiently large $n$, it holds that for $r \in [a_{2n+1}^{1-\eps}, a_{2n+1}]$,
	\begin{align*}
	&\frac{2}{3}r^{3/2}(1-a_{2n+1}^{-3\eps'/2}) \le \frac{2}{3}r^{3/2}(1-(a_{2n+1}^{(1-\eps-\eps')}/r)^{3/2})=\frac{2}{3}(r^{3/2}-a_{2n+1}^{3(1-\eps-\eps')/2})\nn\\
	& =\int_{a_{2n+1}^{1- \eps - \eps'}}^r s^{1/2}ds   \le 2\int_0^r s\psi(s)^{-1} ds \le \frac{3}{2}\int_0^r  s^{1/2} ds = r^{3/2}.
	\end{align*}
	Since $\lim_{n \to \infty} a_{2n+1}^{-3\eps'/2} = 0$, it follows that for all sufficiently large $n$ and $r \in [a_{2n+1}^{1-\eps}, a_{2n+1}]$,
	\begin{equation}\label{intodd}
	\frac{1}{2}r^{3/2} \le 2 \int_0^r s \psi(s)^{-1}ds \le r^{3/2} \quad \text{and hence} \quad r^{1/2} \le \Phi(r) \le 2r^{1/2}.
	\end{equation}
	
	Next, by \eqref{intodd}, for all sufficiently large $n$ and $r \in [a_{2n}^{1-\eps}, a_{2n}]$, we get
	\begin{align}\label{inteven}
	&\frac{1}{2}a_{2n-1}^{3/2} \le 2\int_0^{a_{2n-1}} s \psi(s)^{-1} ds  \le 2\int_0^r s\psi(s)^{-1} ds \nn\\
	&\quad = 2\int_0^{a_{2n-1}} s \psi(s)^{-1} ds +2 \int_{a_{2n-1}}^r \frac{s}{s^4+\psi(a_{2n-1}) - a_{2n-1}^4} ds.
	\end{align}
	Note that for all sufficiently large $n$, by \eqref{basicpsi},
	\begin{align*}
	&\int_{a_{2n-1}}^r \frac{s}{s^4+\psi(a_{2n-1}) - a_{2n-1}^4} ds \le
	\int_{a_{2n-1}}^r \frac{s}{(s-a_{2n-1}) s^3+\psi(a_{2n-1})} ds\\
	& \le  \psi(a_{2n-1})^{-1}\int_{a_{2n-1}}^{a_{2n-1}+1} sds + \int_{a_{2n-1}+1}^r s^{-2} ds \le \frac{3}{4}a_{2n-1}^{-1/2}(a_{2n-1}+1) + 1 \le a_{2n-1}^{1/2}.
	\end{align*}
	Thus, by combining the above inequality with \eqref{intodd} and \eqref{inteven}, we have that for all sufficiently large $n$ and $r \in [a_{2n}^{1-\eps}, a_{2n}]$,
	$$
	\frac{1}{2} \log r \le \frac{1}{2}a_{2n-1}^{3/2} \le 2\int_0^r s\psi(s)^{-1} ds \le (1+a_{2n-1}^{-1})a_{2n-1}^{3/2} \le \frac{3}{2}a_{2n-1}^{3/2} \le \frac{3}{2(1-\eps)} \log r,
	$$
	and hence
	$$
	\frac{2(1-\eps)r^2}{3 \log r}\le \Phi(r) \le \frac{2r^2}{\log r}.
	$$
	This completes the proof. 	\qed

	Let $t_{2n}=a_{2n}^2/(\log a_{2n})$ and $t_{2n+1}=a_{2n+1}^{1/2}$ for $n \ge 1$. Since $\exp(x^{3/2}) \ge 4x^4$ for $x \ge 10$, we can check that $t_{n+1} \ge 4t_n$ for all $n \ge 2$. As a corollary to Lemma \ref{basic}, we obtain the following estimates for the inverse functions of $\Phi$ and $\psi$, respectively.

	\begin{lemma}\label{inverse}
		There exists $N \in \N$ such that for every $n \ge N$, the following estimates hold: 

		\noindent {\rm (1)} For every $t \in [t_{2n+1}/2,t_{2n+1}]$, it holds that
		$$
		\Phi^{-1}(t) \asymp \psi^{-1}(t) \asymp t^{2}.
		$$
		
		\noindent {\rm (2)} For every $t \in [t_{2n}/2,t_{2n}]$, it holds that
		$$
		\Phi^{-1}(t) \asymp t^{1/2} (\log t)^{1/2} \quad \text{and} \quad \psi^{-1}(t) \asymp t^{1/4}.
		$$

	\end{lemma}
	\proof
		(1) For all sufficiently large $n$ and $t \in [t_{2n+1}/2, t_{2n+1}]$, by \eqref{scalepsi}, \eqref{scalephi} and Lemma \ref{basic}(1), we have $\Phi(t^2) \asymp \Phi(a_{2n+1}) \asymp \psi(t^2)\asymp \psi(a_{2n+1}) \asymp t_{2n+1} \asymp t$. Then, we get the result from \eqref{scalepsi} and \eqref{scalephi}.
		
	(2) For all sufficiently large $n$ and $t \in [t_{2n}/2, t_{2n}]$, we have $t^{1/2}(\log t )^{1/2} \asymp a_{2n}$ and $t^{1/4} \asymp a_{2n}^{1/2}(\log a_{2n})^{-1/4}$. Since for all sufficiently large $n$, $\Phi(a_{2n}) \asymp \psi(a_{2n}^{1/2}(\log a_{2n})^{-1/4}) \asymp t_{2n} \asymp t$ by Lemma \ref{basic}(2) with $\eps=2/3$, we obtain the results.	\qed

	\vspace{2mm}
	
	\subsubsection{Construction of subordinator and its transition density estimates}
	
	Let $S$ be a subordinator without drift whose L\'evy measure $\nu(dr)$ is given by
	$$
	\nu(dr)=\frac{1}{r\psi(r)}dr,
	$$
	i.e., the Laplace exponent is given by $\phi(\lambda)=\int_0^\infty (1-e^{-\lambda s})\nu(ds)$. Since $\nu$ satisfies the condition {\bf (E)} with $T_0=0$, we see that the subordinator $S_t$ has a transition density function $p(t,x)$ for all $t>0$. The following theorem is the main result in this example.
	
	Recall that $b(t)=(\phi' \circ H^{-1})(1/t)$ for $t>0$.
	
	\begin{thm}
		\noindent {\rm (1)} For every $n \ge 2$ and $t \in [(1/2)t_{2n+1}, t_{2n+1}]$, it holds that
		$$
		p(t,tb(t)+y) \asymp t^{-2} \wedge \frac{t}{y\psi(y)} \;\; \text{for all} \;\; y \ge 0.
		$$
		
		\noindent {\rm (2)} For every $n \ge 2$ and $t \in [(1/2)t_{2n}, t_{2n}]$, it holds that
		$$
		p(t,tb(t)+y) \simeq t^{-1/2}(\log t)^{-1/2} \wedge \left(\frac{t}{y\psi(y)} + t^{-1/2}(\log t)^{-1/2} \exp \big(-c \frac{y^2}{t \log t}\big) \right)  \;\; \text{for all} \;\; y \ge 0.
		$$
	\end{thm}
	
	\begin{remark}
		{\rm
			Note that for every $n \ge 2$, $t \in [(1/2)t_{2n}, t_{2n}]$ and $y \in [a_{2n}, a_{2n}(\log a_{2n})^{1/3}]$, since $\lim_{n \to \infty}a_{2n} = \infty$, we have
			\begin{align*}
			&t^{-1/2}(\log t)^{-1/2}\exp \big(-c_1 \frac{y^2}{t \log t}\big) \ge c_2a_{2n}^{-1} \exp \big(-c_3 \frac{y^2}{a_{2n}^2}\big) \\
			&\;\; \ge c_2a_{2n}^{-1} \exp \big(-c_3 (\log a_{2n})^{-1/3}\log a_{2n}\big) =c_2a_{2n}^{-1-c_3(\log a_{2n})^{-1/3}}  \ge c_4a_{2n}^{-2},
			\end{align*}
			while
			$$
			\frac{t}{y\psi(y)} \le c_5\frac{a_{2n}^2}{a_{2n}^{1+4}\log a_{2n}} = c_5a_{2n}^{-3} (\log a_{2n})^{-1}.
			$$
			Hence, we see that the exponential term is the dominating factor in heat kernel estimates at those intervals. Therefore, we deduce that although the lower index $\alpha_1<1$, the exponential term in \eqref{e:main2} is indispensable in heat kernel estimates. (cf. \cite[Theorem 1.5]{BKKL}.)
		}
	\end{remark}

	{\it Proof of Theorem 4.3.} Since $S$ satisfies the condition {\bf (G)}, by Lemmas \ref{l:smallH} and \ref{l:largeH}, 
	\begin{align}\label{Hphi}
	H(r^{-1}) \asymp 2r^{-2} \int_0^r \frac{s}{\psi(s)}ds = \Phi(r)^{-1}, \qquad \;\; w(r) \asymp r\nu(r) = \psi(r)^{-1} \quad &\text{for all} \;\; r>0.
	\end{align}
Hence, by the scaling property of the function $w$, we get
\begin{equation}\label{Hw}
H^{-1}(1/t) \asymp \Phi^{-1}(t)^{-1} \qquad \text{and} \qquad w^{-1}(2e/t) \asymp \psi^{-1}(t) \quad \text{for all} \;\; t>0.
\end{equation}
We simply denote by $\theta=\theta(t, y/(8e^2)) \in [w^{-1}(2e/t), H^{-1}(1/t)^{-1}]$.

\smallskip

	(1) By Lemma \ref{inverse}(1) and \eqref{Hw}, there is $N \in \N$ such that for all $n \ge N$,
	\begin{equation}\label{pj}
	H^{-1}(1/t)^{-1} \asymp w^{-1}(2e/t) \asymp t^2 \quad \text{for all} \;\; t \in [(1/2)t_{2n+1}, t_{2n+1}].
	\end{equation} 
	Since $H^{-1}(1/t)^{-1} \asymp w^{-1}(2e/t) \asymp 1$ for $2 \le n <N$ and $t \in [(1/2)t_{2n+1}, t_{2n+1}]$, after taking comparison constants larger, we can see that \eqref{pj} holds for all $n \ge 2$.  
	It follows that $\theta(t,y)\asymp t^2$ for all $t \in [(1/2)t_{2n+1}, t_{2n+1}]$ and $y \ge 0$. Hence, by Lemma \ref{boundary}, for fixed $a>0$ and all $t \in [(1/2)t_{2n+1}, t_{2n+1}]$,
	\begin{align*}
	\exp \big(- \frac{ay}{\theta}\big) \asymp 1 \quad \text{for} \; y \in [0,H^{-1}(1/t)^{-1}]
	\end{align*}
	and
	\begin{align*}
	\exp \big(-\frac{ay}{\theta}\big) \le \exp \big(- \frac{ay}{w^{-1}(2e/t)}\big) \le c_1 \frac{t \nu(y)}{H^{-1}(1/t)} \quad \text{for} \; y \in (H^{-1}(1/t)^{-1},\infty).
	\end{align*}
	Therefore, we get the result from Corollary \ref{c:main1}.
	
	\smallskip

	(2) By Lemma \ref{inverse}(2), \eqref{Hw} and using the same argument as the one given in the proof of (1), we get that for all $n \ge 2$,
	\begin{equation}\label{mix}
	H^{-1}(1/t)^{-1} \asymp  t^{1/2}(\log t)^{1/2} \;\; \text{and} \;\; w^{-1}(2e/t) \asymp t^{1/4}   \quad \text{for all} \;\;  t \in [(1/2)t_{2n}, t_{2n}].
	\end{equation} 
	Then, by \eqref{Hphi}, \eqref{mix} and Lemma \ref{basic}(2) with $\eps=4/5$, we have that for all $t \in [(1/2)t_{2n}, t_{2n}]$,
	$$
	D(t) \asymp t\max_{s \in [w^{-1}(2e/t),H^{-1}(1/t)^{-1}]} \frac{s}{\Phi(s)} \asymp t\max_{s \in [w^{-1}(2e/t),H^{-1}(1/t)^{-1}]} \frac{\log s}{s} \asymp t^{3/4} \log t.
	$$
	From \eqref{basicpsi}, we see that for all $t \in  [(1/2)t_{2n}, t_{2n}]$ and $y \ge D(t) \ge ct^{3/4} \log t$,
	\begin{align*}
	t^{-1/2}(\log t)^{-1/2} \exp \big(-c \frac{y^2}{t \log t}\big) &\le ct^{-1/2}(\log t)^{-1/2} \left(\frac{t \log t}{y^2}\right)^{11/2} = c \frac{t}{y^5}\frac{t^4 \log^5 t}{y^6}\le c \frac{t}{y\psi(y)}.
	\end{align*}

	Hence, by \eqref{mix}, Corollaries \ref{c:main1} and \ref{c:lefttail} and Lemma \ref{boundary}, it suffices to show that  
	\begin{equation}\label{thetacal}
	\frac{y}{\theta} \asymp \frac{y^2}{t\log t} \quad \text{for all} \;\; t \in [(1/2)t_{2n},t_{2n}], \; y \in [H^{-1}(1/t)^{-1}, D(t)].
	\end{equation}

	 Since $\theta \in [w^{-1}(2e/t), H^{-1}(1/t)^{-1}]$, by \eqref{mix}, there are constants $c_1,c_2>0$ such that $c_1t^{1/4} \le \theta \le c_2t^{1/2}(\log t)^{1/2}$. Since $t\theta H(\theta^{-1}) = y$, by \eqref{Hphi} and Lemma \ref{basic}(2) with $\eps=4/5$, (as before, it suffices to consider large $t$ only,)
	\begin{align*}
	y \theta = t\theta^2H(\theta^{-1}) \asymp t \frac{\theta^2}{\Phi(\theta)} \asymp t  \log \theta \asymp t \log t.
	\end{align*}
	This proves \eqref{thetacal}. \qed

\bigskip

\section{Relationship between subordinator and symmetric jump processes}\label{s:R}
In this short section, we discuss a resemblance between transtion density estimates on subordinators and heat kernel estimates on symmetric jump processes.

 Let $S$ be a subordinator without drift whose L\'evy density is $\nu(r) = 1/(r\psi(r))$.
 We assume that $\psi$ is non-decreasing and that $S$  satisfies the condition {\bf (G)}.
Define
$$
\Phi(r)= \frac{r^2}{2\int_0^r s\psi(s)^{-1}ds}.
$$
Then, by Lemmas \ref{l:smallH} and \ref{l:largeH}, we have
\begin{equation}\label{HPhi}
 H(r) \asymp \Phi(r^{-1})^{-1}  \qquad  \text{and} \qquad w(r) \asymp \psi(r)^{-1} \quad \text{for all} \;\; r>0.
\end{equation}

\smallskip
Following \cite[(1.16)]{BKKL}, we define
$$\mathscr{K}_{\infty}(r):= \begin{cases}
\sup_{1\le s\le r} s^{-1}\Phi(s) \quad &\text{if}\quad r \ge 1, \\
\Phi(r) \quad &\text{if} \quad 0<r<1.
\end{cases}$$
If we further assume that {\bf (L.Mixed)} holds, that is $\alpha_3>1$, then we can see from \eqref{HPhi} that
\begin{equation}\label{HPhi1}
\mathscr{K}_{\infty}(r) \asymp \mathscr{H}(r)  =\inf_{s \ge r} \frac{1}{sH(s^{-1})}\qquad \text{for all} \;\; r \ge 1.
\end{equation}

\smallskip

Let $X=(X_t, x \in \R, t \ge 0)$ be a pure jump symmetric Markov process on $\R$ whose jumping kernel $J(x,y)$ satisfies
$$
J(x,y) \asymp \frac{1}{|x-y| \psi(|x-y|)}, \quad x,y \in \R,
$$
that is, its associated Dirichlet form $(\sE, \sF)$ in $L^2(\R)$ is given by
\begin{align*} 
\sE(f,g)&= \int_{\R \times \R \setminus \text{diag}} (f(x)-f(y))(g(x)-g(y))J(x,y)dxdy, \quad f,g \in \sF,\\
\sF&= \{f \in L^2(\R): \sE(f,f)<\infty\}.
\end{align*}

According to \cite[Theorem 1.5(2)]{BKKL}, under the conditions {\bf (G)} and  {\bf (L.Mixed)}, the process $X$ admits a transition density $p^X(t,x,y)$ enjoying the following estimates: for all $(t,x,y) \in [1, \infty) \times \R \times \R$,
\begin{align*}
p^X(t,x,y) &\simeq \Phi^{-1}(t)^{-1} \wedge \left(\frac{t}{|x-y|\psi(|x-y|)} + \Phi^{-1}(t)^{-1}\exp \big(-c \frac{|x-y|}{\mathscr{K}_\infty^{-1}(t/|x-y|)} \big) \right).
\end{align*}
Hence, in view of \eqref{HPhi}, \eqref{HPhi1} and Corollary \ref{cor3}, we see that if the L\'evy density for  a subordinator and a symmetric jump process  are decaying in the same order and the conditions {\bf (G)} and  {\bf (L.Mixed)} hold, then right tail estimates on the transition density for the subordinator and off-diagonal estimates on the one for the symmetric jump process on $\R$ are the same.

\vspace{.1in}

%
%\noindent
%{\bf Acknowledgements:} 
%

\end{document}